\documentclass[reqno]{amsart}
\usepackage{amssymb}
\usepackage{amscd}

\newtheorem{theorem}{Theorem}[section]

\newtheorem{lemma}[theorem]{Lemma}
\newtheorem{proposition}[theorem]{Proposition}

\newtheorem{definition}[theorem]{Definition}

\theoremstyle{definition}
\theoremstyle{remark}
\numberwithin{equation}{section}

\theoremstyle{plain}

\begin{document}
\title[Half-sided modular inclusions]{Extension of the structure
theorem of Borchers and its application to half-sided modular
inclusions}
\author{Huzihiro Araki}
\address{Research Institute for Mathematical Sciences, Kyoto University,
Sakyoku, Kyoto 606-8205, Japan}
\email{araki@kurims.kyoto-u.ac.jp}
\author{L\'aszl\'o Zsid\'o}
\address{Department of Mathematics, University of Rome ``Tor Vergata'',
Via della Ricerca Scientifica, 00133 Rome, Italy}
\email{zsido@axp.mat.uniroma2.it}
\date{November 7, 2004}
\renewcommand{\subjclassname}{\textup{2000} 
 Mathematics Subject Classification}
\thanks{The second author was supported by MIUR, INDAM and EU}
\dedicatory{Dedicated to Professor D. Buchholz on his
$\, 60^{\text{th}}$ birthday}

\begin{abstract}
A result of H.-W. Wiesbrock is extended from the case of a common
cyclic and separating vector for the half-sided modular inclusion
$N\subset M$ of von Neumann algebras to the case of a common
faithful normal semifinite weight and at the same time a gap in
Wiesbrock's proof is filled in.
\end{abstract}
\maketitle

\section{Introduction}

J. Bisognano and E. Wichmann (\cite{B-W}) made a discovery about the
connection of the modular operator and the modular conjugation for the
von Neumann algebra generated by quantum fields in a wedge region of
the Minkowski space-time with kinematical transformations, namely
pure Lorentz transformation and the TCP$_1$ operator.

H. J. Borchers (\cite{Bo2}) formulates an important feature of this
connection in the abstract setting of a pair of von Neumann algebras
$N\subset M$ with a common cyclic and separating vector $\Omega\,$,
and a one-parameter group of unitaries $U(\lambda )$ having a positive
generator, which induces a semigroup of endomorphisms of $(M\, ,\,
\Omega )\,$, obtaining a commutation relation of $U(\lambda )$ with
the modular operator and the modular conjugation for $(M\, ,\,
\Omega )\,$, which reproduces the kinematical commutation relations
in the Bisognano-Wichmann situation.

A further development has been achieved by H.-W. Wiesbrock (\cite{Wie1},
\cite{Wie2}, \cite{Wie3}, \cite{Wie4}, \cite{Wie5}, \cite{Wie6}), who
introduces the notion of the half-sided modular inclusion and obtains
an underlying group structure (cf. \cite{Sch}), as well as an
imbedding of the canonical endomorphisms of the subfactor theory into
a one-parameter semigroup of endomorphisms in this specific situation.
Thus he gets a correspondence between 2-dimensional chiral conformal
field theories and a class of type III$_1$ subfactors.

Unfortunately, there is a gap in Wiesbrock's proof of his basic
theorem (\cite{Wie1}, Theorem 3, Corollary 6 and Corollary 7,
\cite{Wie6}). We will fill in this gap in Wiesbrock's proof and
further generalize the result to the case of a common normal
semifinite faithful weight.

As a basic tool to prove general half-sided modular inclusion results,
we generalize a structure theorem of H. J. Borchers (\cite{Bo3},
Theorem B) considerably, making Borchers' proof at the same time more
transparent.

The extension from the state case to the weights turns out to be not
straightforward. For this purpose we introduce as a basic tool the
notion of a Hermitian map by using modular structure.

It seems that before the summer of 1995, when we independently
noticed the gap in the proof of Wiesbrock's half-sided modular
inclusion theorem, this gap was generally overlooked. Thanks to
Prof. Detlev Buchholz, who has been visiting the first named author
in the fall of 1995, we learned about each others insights and
started to collaborate on this paper. The first version of the
paper, containing a complete proof of the General Half-sided
Modular Inclusion Theorem, Theorem \ref{hsmi}, was available
already at the end of 1995. It had a restricted circulation,
but it was presented at several conferences. Other topics, like
Theorem \ref{type} on the structure and type of the involved
von Neumann algebras and Proposition \ref{counterex.} on
pathologies of the analytic extension of orbits of one-parameter
automorphism groups, are of more recent date.

We notice that since 1995 a number of papers appeared, containing
proposals for a complete proof of the half-sided modular inclusion
theorem (see, for example, \cite{Fl}, Section 3 and \cite{Bo-Y},
pp. 608 and 609). None the less, till now we have no knowledge of
a completely elaborated proof, even in the state case.

\bigskip
\section{Main Results}

\bigskip
\noindent (a) {\sc Notations and facts from the Modular Theory
of von Neumann}
\\ \phantom{xx} {\sc Algebras} (see, for example, \cite{S-Z}, Chapter 10).
\bigskip

For two Hilbert spaces $H$ and $K$ we denote by $B(K,H)$ the Banach
space of all bounded linear maps from $K$ to $H\,$. $B(H,H)$ will
be denoted simply by $B(H)\,$. If $T$ is a not necessarily
everywhere defined linear operator from $K$ to $H\,$, then Dom$\, (T)$
will stand for the domain of $T\,$.

We denote the weak and the strong operator topology on $B(K,H)$
respectively by $wo$ and $so\,$. The weak topology defined on $B(K,H)$
by all linear functionals belonging to the norm-closure of the
$wo$-continuous linear functionals in the dual of $B(K,H)$ will be
denoted by $w\,$. Further, the locally convex vector space topology
defined on $B(K,H)$ by the seminorms
\begin{equation*}
B(K,H)\ni T\longmapsto\varphi (T^*T)^{1/2}\, ,
\end{equation*}
where $\varphi$ ranges over all $w$-continuous positive linear
functionals on $B(K)\,$, will be denoted by $s\,$. We notice that
on any bounded subset of $B(K,H)$ $wo=w$ and $so=s\,$.

For a weight $\varphi$ on a von Neumann algebra $M$ we use the
standard notations:
\begin{align*}
\mathfrak N_\varphi &=\{ x\in M\, ;\,\varphi (x^*x)<+\infty\}
\hspace{0.2 cm} (\text{left ideal})\, , \\
\mathfrak M_\varphi &=(\mathfrak N_\varphi )^*\mathfrak N_\varphi =
\,\text{the linear span of}\,\{ a\in M^+\, ;\,\varphi (a)<+\infty\} \\
&\hspace{4.45 cm} (\text{hereditary $*$-subalgebra})\, ,\notag \\
\mathfrak A_\varphi &= (\mathfrak N_\varphi )^*\cap
\mathfrak N_\varphi \supset \mathfrak M_\varphi
\hspace{1.15 cm} (*\text{-subalgebra})\, .
\end{align*}
We notice that for $\, a\in M^+$ we have $\, a\in\mathfrak M_\varphi
\,\Longleftrightarrow\, \varphi (a)<+\infty\,$.

A von Neumann algebra $M$ on a Hilbert space $H$ is in \emph{standard
form with respect to a normal semifinite faithful weight} $\varphi$
on $M$ if there is a linear map with dense range
\begin{center}
$\mathfrak N_\varphi\ni x\longmapsto x_\varphi\in H$
\end{center}
such that
\begin{center}
$\varphi (x^*x)=\| x_\varphi\|^2\, ,\, (a\, x)_\varphi =a\, x_\varphi
\, ,\qquad x\in \mathfrak N_\varphi\, ,\, a\in M\, .$
\end{center}
In particular, by the faithfulness of $\varphi\,$, the map
$x\longmapsto x_\varphi$ is injective. We notice also that the above
map $x\longmapsto x_\varphi$ is unique up to natural unitary
equivalence. If $\varphi$ is bounded, then $\xi_\varphi =(1_H)_\varphi$
is a cyclic and separating vector for $M$ and $x_\varphi =
x\,\xi_\varphi$ for all $x\in M=\mathfrak N_\varphi\,$. In this case
$\varphi$ is the vector form $M\ni x\longmapsto
\omega_{\xi_\varphi} (x)=
(x\,\xi_\varphi\, |\,\xi_\varphi )\,$.

Let $M$ be a von Neumann algebra on a Hilbert space $H\,$, in
standard form with respect to a normal semifinite faithful weight
$\varphi$ on $M\,$. Then the antilinear operator
\begin{center}
$H\supset \{ x_\varphi\, ;\, x\in \mathfrak A_\varphi\}\ni x_\varphi
\longmapsto (x^*)_\varphi$
\end{center}
has closure $S_\varphi$ and the invertible positive selfadjoint
operator $\Delta_\varphi =S_\varphi^{\, *}S_\varphi$ is called the
\emph{modular operator} of $\varphi\,$. If
\begin{equation*}
S_\varphi =J_\varphi\,\Delta_\varphi^{1/2}
\end{equation*}
is the polar decomposition of $S_\varphi\,$, then $J_\varphi$
is an involutive antiunitary operator (antilinear surjective isometry
with $J_\varphi^{\, 2}=1_H$), called the \emph{modular conjugation} of
$\varphi\,$. The operators $\Delta_\varphi$ and $J_\varphi$ satisfy
the commutation relation
\begin{equation}\label{mod.comm}
J_\varphi\,\Delta_\varphi^z =\Delta_\varphi^{-\,\overline z}\,
J_\varphi\, ,\qquad z\in\mathbb C\, ,
\end{equation}
in particular,
\begin{equation}\label{mod.comm-}
S_\varphi =J_\varphi\,\Delta_\varphi^{1/2} =\Delta_\varphi^{-1/2}\,
J_\varphi\, ,\qquad J_\varphi\,\Delta_\varphi^{it}=
\Delta_\varphi^{it}\,J_\varphi\, ,\; t\in\mathbb R\, .
\end{equation}
If $\varphi$ is bounded and $\,\xi_\varphi =(1_H)_\varphi$ is the
associated cyclic and separating vector, then
\begin{equation*}
S_\varphi\,\xi_\varphi =\xi_\varphi\, ,\quad
\Delta_\varphi\,\xi_\varphi =\xi_\varphi\, ,\quad
J_\varphi\,\xi_\varphi =\xi_\varphi\, .
\end{equation*}

The fundamental result of the modular theory claims that
\begin{equation}\label{delta}
x\in \mathfrak N_\varphi\, ,\, t\in \mathbb R\;\Longrightarrow\;
\Delta_\varphi^{it}\, x\,\Delta_\varphi^{-it}\in \mathfrak N_\varphi
\, ,\; (\Delta_\varphi^{it}\, x\,\Delta_\varphi^{-it})_\varphi
=\Delta_\varphi^{it}\, x_\varphi\, ,
\end{equation}
so that
\begin{equation}\label{sigma}
M\ni x\longmapsto \sigma^\varphi_t (x) =\Delta_\varphi^{it}\, x\,
\Delta_\varphi^{-it}\in M\, ,\qquad t\in \mathbb R
\end{equation}
defines an $so$-continuous one-parameter group of automorphisms
$(\sigma^\varphi_t)_{\substack{{}\\ t\in \mathbb R}}$ of $M\,$,
called the \emph{modular automorphism group} of $\varphi\,$, and
\begin{equation}\label{j}
J_\varphi\, M\, J_\varphi =M'\, ,\qquad
x\, ,\, y\in \mathfrak N_\varphi \Longrightarrow
x\, J_\varphi\, y_\varphi =J_\varphi\, y\, J_\varphi\, x_\varphi\, ,
\end{equation}
so that $M\ni x\longmapsto J_\varphi\, x^* J_\varphi\in M'$ is
a $*$-antiisomorphism. Moreover, the weight $\varphi$ is invariant
under the action of the modular automorphism group :
\begin{equation}\label{inv}
\varphi\big( \sigma^\varphi_t (a)\big) =\varphi (a)\, ,\qquad
a\in M^+\, ,\, t\in\mathbb R\, .
\end{equation}

The center $Z(M)$ of $M$ is contained in the fixed point von
Neumann subalgebra $\{ x\in M\, ;\,\sigma^\varphi_t (x) =
x\, ,\, t\in\mathbb R \}\subset M\,$, which is usually denoted
by $M^\varphi\,$. On the other hand, $J_\varphi\, z\, J_\varphi
=z^*$ for all $z\in Z(M)\,$. We recall also (see the proof of
\cite{P-T}, Lemma 5.2 or \cite{Z3}, Corollary 1.2) :
\begin{equation}\label{right-mult}
x\in\mathfrak N_\varphi\, ,\, y\in M^\varphi\;\Longrightarrow\;
x\, y\in\mathfrak N_\varphi\, ,\, (x\, y)_\varphi =J_\varphi\,
y^* J_\varphi\, x_\varphi\, .
\end{equation}

Let $e\in M^\varphi$ be a projection and let $\varphi_e$ denote
the restriction of $\varphi$ to $e M e\,$.
By \cite{P-T}, Proposition 4.1 and Theorem 4.6
(see also \cite{S}, Propositions 4.5 and 4.7), $\varphi_e$ is a
normal semifinite faithful weight and its modular group is the
restriction of the modular group of $\varphi$ to $e M e\,$.
Thus, if $\pi_e : e M e\longrightarrow B(e H)$ is the
faithful normal $*$-representation which associates to every
$x\in e M e$ the restriction $x\, |\, e H$ considered as a linear
operator $e H\rightarrow e H\,$, then the modular group of the
weight $\varphi_e\circ\pi_e^{-1}$ on $\pi_e (e M e)\,$, that
is $(\pi_e\circ\sigma^{\varphi_e}_t\circ\pi_e^{-1})_{
\substack{{}\\ t\in\mathbb R}}\,$, is implemented by the unitary
group $(\Delta_\varphi^{i t}\, |\, e H)_{
\substack{{}\\ t\in\mathbb R}}$ on $e H\,$. Nevertheless,
$\pi_e (e M e)$ is not always in standard form with respect to
$\varphi_e\circ\pi_e^{-1}$ (indeed, if $M\subset B(\mathbb C^4 )$
is a type I$_2$ factor, in standard form with respect to its
trace, and $e\in M$ is a minimal projection, then $\pi_e (M)$ is
one-dimensional, while its commutant $\pi_e (M)'$ is
four-dimensional, so $\pi_e (M)$ and $\pi_e (M)'$ are not
antiisomorphic).

However, for any projection $e\in M^\varphi\,$,
\begin{equation}\label{stand.red}
\pi =\pi_{\substack{{}\\ e J_\varphi e J_\varphi}} : e M e\ni
x\longmapsto x\, |\, e J_\varphi e J_\varphi H \in B(e J_\varphi
e J_\varphi H)
\end{equation}
is a faithful normal $*$-representation, such that the von Neumann
algebra $\pi (e M e)$ is in standard form with respect to
$\varphi_e\circ\pi^{-1}$ (cf. \cite{Haa1}, Lemma 2.6). Moreover,
$\Delta_\varphi$ and $J_\varphi$ commute with $e J_\varphi e
J_\varphi$ and we have the identifications
\begin{equation}\label{red.mod}
\Delta_{\varphi_e\circ\pi^{-1}} =\Delta_\varphi\, |\, e J_\varphi
e J_\varphi H\, ,\qquad J_{\varphi_e\circ\pi^{-1}} =
J_\varphi\, |\, e J_\varphi e J_\varphi H\, .
\end{equation}
For the convenience of the reader, let us outline the proof of
(\ref{red.mod}).

For the faithfulness of $\pi\,$, let $x\in e M e$ be such that
$x\, eJ_\varphi eJ_\varphi =0\,$. Then
\smallskip
\begin{center}
$x\, J_\varphi M eJ_\varphi =x\, eJ_\varphi MJ_\varphi J_\varphi
eJ_\varphi =J_\varphi MJ_\varphi\, x\, eJ_\varphi eJ_\varphi =0\,$,
\end{center}
\smallskip
so $x\, J_\varphi$ vanishes on $M e H\,$, hence on the range of
the central support $z(e)\in Z(M)$ of $e\,$. Thus $x=x\, z(e)=
x\, J_\varphi z(e) J_\varphi =0\,$.

To see that $\pi (e M e)$ is in standard form with respect to
$\psi =\varphi_e\circ\pi^{-1}\,$, first we notice that, according
to (\ref{right-mult}), $\mathfrak N_\psi =
\pi (\mathfrak N_{\varphi_e}) =e\,\mathfrak N_\varphi e\,$.
Next, the linear map
\begin{center}
$\mathfrak N_\psi =\pi (\mathfrak N_{\varphi_e})\ni\pi (x)
\longmapsto x_\varphi =(e x e)_\varphi
\overset{(\ref{right-mult})}{=} eJ_\varphi eJ_\varphi\, x_\varphi
\in e J_\varphi e J_\varphi H$
\end{center}
\smallskip
has dense range. Indeed, every vector in $eJ_\varphi eJ_\varphi H$
belongs to the closure of
\begin{center}
$eJ_\varphi eJ_\varphi\,\{ x_\varphi\, ;\, x\in\mathfrak N_\varphi\}
\overset{(\ref{right-mult})}{=} \{ (e x e)_\varphi\, ;\,
x\in\mathfrak N_\varphi\} = \{ x_\varphi\, ;\, x\in
\mathfrak N_{\varphi_e} \}\, .$
\end{center}
\smallskip
Finally, for every $\pi (x)\in\pi (\mathfrak N_{\varphi_e})
=\mathfrak N_\psi$ and $\pi (a)\in\pi (e\, M\, e)$ hold true:
\begin{equation*}
\psi\big( \pi (x)^*\pi (x)\big) =\varphi_e (x^*x) =\varphi (x^*x)
=\| x_\varphi\|^2\, ,\qquad (a\, x)_\varphi =a\, x_\varphi =
\pi (a) x_\varphi\, .
\end{equation*}

The commutation of $J_\varphi$ with $e J_\varphi e J_\varphi$
follows immediately from the commutation of $e$ with $J_\varphi
e J_\varphi\,$. Let $J_e$ denote the involutive antiunitary
operator
\begin{equation*}
e J_\varphi e J_\varphi H\ni\xi\longmapsto J_\varphi\,\xi\in
e J_\varphi e J_\varphi H\, .
\end{equation*}
Further, using (\ref{mod.comm-}) and
$e\in M^\varphi\,$, we obtain for every $t\in\mathbb R\,$:
\begin{equation*}
e J_\varphi e J_\varphi\,\Delta_\varphi^{i t} =
e J_\varphi e\,\Delta_\varphi^{i t} J_\varphi =
e J_\varphi\,\Delta_\varphi^{i t} e J_\varphi =
e\,\Delta_\varphi^{i t} J_\varphi e J_\varphi =
\Delta_\varphi^{i t} e J_\varphi e J_\varphi\, .
\end{equation*}
Thus also $\,\Delta_\varphi$ commutes with $e J_\varphi e
J_\varphi\,$, so
\begin{equation*}
e J_\varphi e J_\varphi H\supset \text{Dom}\, (\Delta_\varphi )
\cap (e J_\varphi e J_\varphi H)\ni\xi\longmapsto \Delta_\varphi
\,\xi\in e J_\varphi e J_\varphi H
\end{equation*}
is an invertible positive selfadjoint operator $\Delta_e\,$,
whose positive selfadjoint square root is
\begin{equation*}
e J_\varphi e J_\varphi H\supset\text{Dom}\, (\Delta_\varphi^{1/2})
\cap (e J_\varphi e J_\varphi H)\ni\xi\longmapsto
\Delta_\varphi^{1/2}\,\xi\in e J_\varphi e J_\varphi H\, .
\end{equation*}

Since, for every $\pi (x)\in\pi (\mathfrak A_{\varphi_e}) =
\mathfrak A_\psi\,$,
\begin{equation*}
S_\psi x_\varphi =(x^*)_\varphi = S_\varphi x_\varphi =J_\varphi\,
\Delta_\varphi^{1/2} x_\varphi =J_e\,\Delta_\varphi^{1/2} x_\varphi
=J_e\,\Delta_e^{1/2}x_\varphi\, ,
\end{equation*}
we deduce that $S_\psi\subset J_e\,\Delta_e^{1/2}\,$. For the
equality $S_\psi =J_e\,\Delta_e^{1/2}$, which will imply
(\ref{red.mod}), let $\xi\in\text{Dom}\, (\Delta_e^{1/2}) =
\text{Dom}\, (\Delta_\varphi^{1/2})\cap (eJ_\varphi eJ_\varphi H)
=\text{Dom}\, (S_\varphi )\cap (e J_\varphi e J_\varphi H)$
be arbitrary. Then there is a sequence $(x_n)_{n\geq 1}$ in
$\mathfrak A_\varphi$ such that $(x_n)_\varphi\longrightarrow\xi$
and $(x_n^{\, *})_\varphi\longrightarrow S_\varphi\,\xi\,$.
By (\ref{right-mult}) the sequence $(e x_n e)_{n\geq 1}$ belongs
to $\mathfrak A_{\varphi_e}$ and we have
\begin{align*}
(e x_n e)_\varphi &=eJ_\varphi eJ_\varphi\,(x_n)_\varphi
\longrightarrow eJ_\varphi eJ_\varphi\,\xi =\xi\, , \\
S_\psi\, (e x_n e)_\varphi &=\big( (e x_n e)^*\big)_\varphi =
eJ_\varphi eJ_\varphi\, (x_n^{\, *})_\varphi
\longrightarrow eJ_\varphi eJ_\varphi\, S_\varphi\,\xi\, .
\end{align*}
Now the closedness of the graph of $S_\psi$ yields $\xi\in
\text{Dom}\, (S_\psi )\,$.

For a projection $p\in Z(M)\subset M^\varphi$ we have $J_\varphi
p J_\varphi =p\,$, so (\ref{red.mod}) yields
\begin{equation}\label{centr.red.mod}
\Delta_{\varphi_p\circ\pi_p^{-1}} =\Delta_\varphi\, |\, p H\, ,
\qquad J_{\varphi_p\circ\pi_p^{-1}} = J_\varphi\, |\, p H\, .
\end{equation}

Let $M\neq\{ 0\}$ be a von Neumann algebra, in standard form
with respect to a normal semifinite faithful weight $\varphi$
on $M\,$. Then the \emph{Connes spectrum} $\Gamma (\sigma^\varphi )$
of the modular automorphism group $\sigma^\varphi$ of $\varphi$
is the intersection of the Arveson spectra of all modular
automorphism groups $\sigma^{\varphi_e}$, where $e$ ranges over
all non-zero projections $e\in M^\varphi\,$. By \cite{Co},
Lemme 1.2.2 and Th\'eor\`eme 2.2.4 (see also \cite{S}, Theorem 3.1
and Proposition 16.3), $\Gamma (\sigma^\varphi )$ is a closed
additive subgroup of $\mathbb R$ and it does not depend on
the choice of $\varphi\,$, so it can be denoted (like in \cite{P},
8.15) by $\Gamma (M)\,$. Furthermore, by \cite{Co}, Lemme 3.2.2
(see also \cite{S}, Proposition 28.1), $\lambda\in\Gamma (M)$
if and only if $e^\lambda$ belongs to the spectrum $\sigma
(\Delta_{\varphi_e})$ of $\Delta_{\varphi_e}$ for all non-zero
projections $e\in M^\varphi\,$.

According to \cite{Co}, page 28, the von Neumann algebra $M\neq
\{ 0\}$ is called to be \emph{of type} III$_1$ if $\Gamma (M) =
\mathbb R\,$, or equivalently, if $\sigma (\Delta_{\varphi_e})=
[\, 0,+\infty )$ for every non-zero projection $e\in M^\varphi\,$.
By (\ref{red.mod}) we have also:
\begin{equation}\label{typeIII_1}
M\text{ is of type III}_1\,\Longleftrightarrow\,
\left\{ \begin{array}{l} \sigma \big(\,\Delta_\varphi\, |\,
e J_\varphi e J_\varphi H\,\big) = [\, 0,+\infty ) \\
\text{for every projection }0\neq e\in M^\varphi
\end{array}\, . \right.
\end{equation}

\bigskip
\noindent (b) {\sc The general Half-sided Modular Inclusion Theorem}.
\bigskip

Let $M$ be a von Neumann algebra on a Hilbert space $H\,$, in
standard form with respect to a normal semifinite faithful weight
$\varphi$ on $M\,$. Let further $N\subset M$ be a von Neumann
subalgebra such that the restriction $\psi$ of $\varphi$ to $N$ is
semifinite. If $\{ y_\varphi\, ;\, y\in\mathfrak N_\psi\}$ is dense
in $H\,$, then $N$ is in standard form with respect to $\psi$ such
that $y_\psi =y_\varphi$ for all $y\in\mathfrak N_\psi\,$.
This happens, for example, if $N\subset M\subset B(H)$ are von
Neumann algebras having a common cyclic and separating vector
$\xi_o\,$, and $\varphi$ is the vector form $M\ni x\longmapsto
(x\,\xi_o\, |\,\xi_o )\,$.

In the above situation, owing to (\ref{j}), we have
\begin{equation*}
J_\psi\, J_\varphi\, M\, J_\varphi\, J_\psi =J_\psi\, M'\, J_\psi
\subset J_\psi\, N'\, J_\psi =N\, ,
\end{equation*}
so the unitary $J_\psi\, J_\varphi$ implements a unital
$*$-homomorphism
\begin{equation*}
M\ni x\longmapsto \text{Ad}\, (J_\psi\, J_\varphi)(x)=J_\psi\,
J_\varphi\, x\, J_\varphi\, J_\psi\in N\subset M\, ,
\end{equation*}
considered by R. Longo (\cite{Lo1}, \cite{Lo2}) and called the
\emph{canonical endomorphism} of the inclusion $N\subset M\,$. The
canonical endomorphism $\gamma\,$, in particular the tunnel
\begin{equation}\label{tunnel}
M\supset N\supset \gamma (M)\supset \gamma (N)\supset \gamma^2 (M)
\supset \gamma^2(N)\supset\,\ldots\, ,
\end{equation}
plays an important role in the Subfactor Theory (see \cite{Lo3}
and \cite{Kos}).

Let $\mathcal P_+^\uparrow (1)$ denote the two-dimensional
Lie group generated by the hyperbolic rotations
\begin{equation*}
L_t : \mathbb R^2\ni \left( \begin{array}{c} \xi_1\\ \xi_2
\end{array}\right)\longmapsto\left( \begin{array}{rr}
\cosh (2\pi t) &-\sinh (2\pi t)\\ -\sinh (2\pi t) &\cosh (2\pi t)
\end{array}\right) \left(\begin{array}{c} \xi_1\\ \xi_2
\end{array} \right)\in\mathbb R^2\, ,\qquad t\in\mathbb R
\end{equation*}
and the lightlike translations
\begin{equation*}
T_s : \mathbb R^2\ni \left( \begin{array}{c} \xi_1 \\ \xi_2
\end{array}\right)\longmapsto\left( \begin{array}{c} \xi_1 +s\\
\xi_2 +s\end{array}\right)\in\mathbb R^2\, ,\quad s\in\mathbb R
\end{equation*}
(cf. \cite{Ba-Ra}, Ch. 17, \S$\,$2, A), which is the \emph{Poincar\'e
group on the light-ray}. Furthermore, the commutation relation
\begin{equation*}
T_s\, L_t =L_t\, T_{e^{2\pi t} s}\, ,\qquad s\, ,\, t\in\mathbb R
\end{equation*}
implies that $(T_{s_1} L_{t_1})\, (T_{s_2}\, L_{t_2}) =T_{s_1 +
e^{-2\pi t_1} s_2}\, L_{t_1+t_2}\,$. Therefore, endowing $\mathbb R^2$
with the Lie group structure defined by the composition law
\begin{equation*}
(s_1 ,t_1)\cdot (s_2 ,t_2) = (s_1 +e^{-2\pi t_1} s_2\, ,\,
t_1+t_2 )\, ,
\end{equation*}
the mapping $\mathbb R^2\ni (s ,t)\longmapsto T_s\, L_t\in
\mathcal P_+^\uparrow (1)$ becomes a Lie group isomorphism. In
particular, $\mathcal P_+^\uparrow (1)$ is connected and simply
connected. On the other hand, the map
$\displaystyle \Big( \begin{array}{cc} e^{-2\pi t} & s\\ 0 & 1
\end{array}\Big)\longmapsto T_s\, L_t$ is a Lie group isomorphism
of the two-dimensional $2\times 2$ matrix group
$\displaystyle
\mathcal G =\Big\{ \Big( \begin{array}{cc} e^{-2\pi t} & s\\
0 & 1\end{array}\Big)\, ;\, s\, ,\, t\in\mathbb R\Big\}$ onto
$\mathcal P_+^\uparrow (1)\,$. If we identify
$\mathcal P_+^\uparrow (1)$ with $\mathcal G$ along the above
isomorphism, the Lie algebra
$\mathfrak p_+^\uparrow (1)$ of $\mathcal P_+^\uparrow (1)$ will
be identified with the Lie algebra $\mathfrak g$ of $\mathcal G\,$,
and the exponential map $\mathfrak p_+^\uparrow (1)\longrightarrow
\mathcal P_+^\uparrow (1)$ with the exponential map
$\mathfrak g\longrightarrow\mathcal G\,$, that is with the usual
exponentiation of the matrices belonging to $\mathfrak g\,$.
We notice that $\mathfrak g$ is the set of all $2\times 2$ real
matrices $X$ such that $\exp (t\, X)\in\mathcal G\,$,
$t\in\mathbb R\,$, and $[X,Y]=X Y-Y X$ for all $X\, ,\, Y\in
\mathfrak g\,$. The elements
\begin{equation}\label{gener}
X_1=\Big( \begin{array}{cc} -2\pi & 0\\ 0 & 0\end{array}\Big)\, ,
\quad X_2=\Big( \begin{array}{cc} -2\pi & 2\pi\\ 0 & 0\end{array}
\Big)\, ,\quad X_3=\Big( \begin{array}{cc} 0 & 1\\ 0 & 0\end{array}
\Big)
\end{equation}
of $\mathfrak g\equiv\mathfrak p_+^\uparrow (1)$ are of particular
interest: we have
\begin{equation}\label{commutation}
X_3=\frac 1{2\pi}\,\big( X_2 -X_1\big)\, ,\quad [\, X_2\, ,\, X_1]=
4\,\pi^2 X_3\, ,
\end{equation}
any two of $X_1\, ,X_2\, ,X_3$ is a basis for $\mathfrak g\equiv
\mathfrak p_+^\uparrow (1)$ and
\begin{equation}\label{gen.exp}
\exp (t\, X_j) =
\Big( \begin{array}{cc} e^{-2\pi t} & 0\\ 0 & 1
\end{array}\Big)\, ,\;\Big( \begin{array}{cc}
e^{-2\pi t} & 1-e^{-2\pi t}\\ 0 & 1\end{array}\Big)\, ,\;
\Big( \begin{array}{cc}1 & t\\ 0 & 1\end{array}\Big)
\end{equation}
for $j=1\, ,\, 2\, ,\, 3\,$, respectively.

According to the general theory of unitary
representations of Lie groups (see e.g. \cite{Ba-Ra}, Ch. 11, \S$\,$1,
B or \cite{Schm}, Section 10.1), if $\pi$ is an $so$-continuous
unitary representation of $\mathcal G\equiv\mathcal P_+^\uparrow (1)$
on a Hilbert space $H$ and $\mathcal D_{\text{G}}(\pi )$ denotes the
G{\aa}rding subspace of $H$ for $\pi\,$, then the formula
\begin{equation*}
\text{d}\pi (X)\,\xi =\frac {d\,}{d t}\, \pi\big( \exp (t\, X)\big)
\,\xi\,\bigg|_{\, t=0}\, ,\qquad X\in\mathfrak p_+^\uparrow (1)\,
,\,\xi\in \mathcal D_{\text{G}}(\pi )
\end{equation*}
defines a representation of the Lie algebra
$\mathfrak p_+^\uparrow (1)$ into the Lie algebra of all
skew-symmetric linear mappings $\mathcal D_{\text{G}}(\pi )
\rightarrow\mathcal D_{\text{G}}(\pi )\,$. Moreover, for every
$X\in\mathfrak p_+^\uparrow (1)\,$, the linear mapping
$\, i\,\text{d}\,\pi (X) : H\supset\mathcal D_{\text{G}}(\pi )
\longrightarrow \mathcal D_{\text{G}}(\pi )\subset H$ is essentially
selfadjoint (see e.g. \cite{Ba-Ra}, Ch. 11, \S$\,$2, Corollary 4 or
\cite{Schm}, Corollary 10.2.11). Therefore, if
$X\in\mathfrak p_+^\uparrow (1)$ and $A$ is the selfadjoint linear
operator in $H\,$, then
\begin{equation}\label{closure}
\pi\big( \exp (t\, X)\big) =\exp (i t A)\,\text{ for all }\,t\in
\mathbb R\quad\Longrightarrow\quad\overline{\text{d}\,\pi (X)} =
i\, A
\end{equation}
(the first $\exp$ is the exponential map of the Lie group
$\mathcal P_+^\uparrow (1)\,$, while the second one indicates
functional calculus).
Indeed, by the definition of $\text{d}\,\pi (X)$
we have $\text{d}\,\pi (X)\subset i\, A\,$, so the selfadjoint
operator $\,\overline{-i\,\text{d}\,\pi (X)}$ is contained in
the selfadjoint operator $A\,$, which implies their equality.

We notice for completeness that, according to \cite{Dix-Mall},
Theorem 3.3, the G{\aa}rding subspace $\mathcal D_{\text{G}}(\pi )$
is actually equal to the set of all $C^\infty$-vectors for $\pi\,$.

We notice also the following simple fact concerning the essential
selfadjointness of sums of symmetric operators:
If $H$ is a Hilbert space, $\mathcal D\subset H$ is a dense linear
subspace and $A\, ,\, B : \mathcal D\longrightarrow\mathcal D$ are
linear operators, then
\begin{equation}\label{sum}
A\, ,\, B\,\text{ symmetric },\;
A+B\,\text{ essentially selfadjoint }\,\Longrightarrow\;
\overline{A} +\overline{B}\subset\overline{A+B}\, .
\end{equation}
Consequently also $\overline{A} +\overline{B}$ is essentially
selfadjoint and $\overline{\overline{A} +\overline{B}} =
\overline{A+B}\,$.

To prove (\ref{sum}), let $\eta\in\text{Dom}\, (\overline{A})
\cap\text{Dom}\, (\overline{B})$ be arbitrary. Then
\begin{equation*}
\big(\eta\, |\, (A+B)\,\xi\big) =(\eta\, |\, A\,\xi )+
(\eta\, |\, B\,\xi )=
(\,\overline{A}\,\eta\, |\,\xi )+
(\,\overline{B}\,\eta\, |\,\xi )=\big(\, (\,\overline{A} +
\overline{B})\,\eta\, |\,\xi\big)\, ,\quad \xi\in\mathcal D\, ,
\end{equation*}
so $\eta$ is in the domain of $(A+B)^* =\overline{A+B}$
and $(\overline{A+B})\,\eta =(A+B)^*\eta =(\,\overline{A} +
\overline{B})\,\eta\,$.

(\ref{sum}) implies that, for any $so$-continuous unitary
representation $\pi$ of $\mathcal P_+^\uparrow (1)\,$,
\begin{equation}\label{sumLie}
\overline{\overline{\text{d}\,\pi (X)} +
\overline{\text{d}\,\pi (Y)}} =\overline{\text{d}\,
\pi (X+Y)}\, ,\qquad X\, ,\, Y\in\mathfrak p_+^\uparrow (1)\, .
\end{equation}

\begin{theorem}\label{hsmi}
{\rm (General Half-sided Modular Inclusion Theorem)}
Let $M$ be a von Neumann algebra on a Hilbert space $H\,$, in
standard form with respect to a normal semifinite faithful weight
$\varphi$ on $M\,$, and $N\subset M$ a von Neumann subalgebra such
that the restriction $\psi$ of $\varphi$ to $N$ is semifinite and
$N$ is in standard form with respect to $\psi\,$. Let us denote for
convenience
\begin{center}
$\Delta_M =\Delta_\varphi\, ,\, J_M=J_\varphi$ and
$\,\Delta_N =\Delta_\psi\, ,\, J_N=J_\psi$ 
\end{center}
and assume the following half-sided modular inclusion $:$
\begin{equation}\label{hsmicond}
\Delta_M^{it}\, N\,\Delta_M^{-it}\subset N\, ,\qquad t\leq 0\, .
\end{equation}
Then
\begin{equation}\label{generator}
\frac 1{2\,\pi}\, (\,\log \Delta_N -\log \Delta_M )\, ,
\end{equation}
defined on the intersection of the domains of $\,\log \Delta_N$ and
$\,\log \Delta_M\,$, is an essentially selfadjoint operator with
positive selfadjoint closure $P$ and, letting
\begin{equation}\label{generation}
U(s) =\exp (i s P)\, ,\qquad s\in \mathbb R\, ,
\end{equation}
we have the following $:$
\begin{itemize}
\item[(1)] $\quad\Delta_M^{-it}\, U(s)\,\Delta_M^{it} =
\Delta_N^{-it}\, U(s)\,\Delta_N^{it} =U(e^{2\,\pi\, t} s)\, ,\qquad
s\, ,\, t\in \mathbb R\, ;$
\item[(2)] $\quad J_M\, U(s)\, J_M =J_N\, U(s)\, J_N =U(-s)\, ,\qquad
s\in \mathbb R\, ;$
\item[(3)] $\quad U\big( 1-e^{2\,\pi\, t}\big) =
\Delta_N^{-it}\,\Delta_M^{it}\,\text{ and }\,\Delta_N^{it} =
U(1)\,\Delta_M^{it}\, U(1)^*\, ,\qquad t\in \mathbb R\, ;$
\item[(4)] $\quad U(2)=J_N\, J_M\,\text{ and }\, J_N=U(1)\, J_M\,
U(1)^*\, ;$
\item[(5)] $\quad N=U(1)\, M\, U(1)^*\, ;$
\item[(6)] $\quad U(s)\, M\, U(s)^*\subset M\, ,\qquad s\geq 0\, ;$
\item[(7)] $\gamma_s = \text{\rm Ad}\, U(s)\, ,\, s\geq 0$ is an
$so$-continuous one-parameter semigroup of $*$-endomorphisms of $M$
such that $\gamma_2$ is equal to the canonical endomorphism $\gamma
=\text{\rm Ad}\, (J_N\, J_M)$ of the inclusion $N\subset M\, ;$
thus $\gamma_s(M)\, ,\, s\geq 0\,$, provide a continuous
interpolation of the tunnel {\rm (\ref{tunnel}):}
\begin{align*} \phantom{xxxxx}
M\supset \gamma_1(M)=N &\supset\gamma_2(M)=\gamma (M)\supset
\gamma_3(M)=\gamma (N)\supset \\
&\supset\gamma_4(M)=\gamma^2(M)\supset\gamma_5(M)=\gamma^2(N)
\supset\,\ldots\; ;
\end{align*}
\item[(8)] For $x\in M$ and $\, s\geq 0$ we have
\begin{equation*} \phantom{xxxxx}
x\in \mathfrak N_\varphi\,\Longleftrightarrow\,
U(s)\, x\, U(s)^*\in \mathfrak N_\varphi\,\Longrightarrow\,
\big(U(s)\, x\,U(s)^*\big)_\varphi =U(s)\, x_\varphi\; ;
\end{equation*}
\item[(9)] The weight $\varphi$ on $M$ is invariant under $\gamma_s$
for every $s\geq 0\, ;$
\item[(10)] $\{ \Delta_M^{it}\, ,\,\Delta_N^{is}\, ;\, t\, ,\, s\in
\mathbb R\}$ generates a group of unitary operators on $H\,$, which is
the image of an $so$-continuous unitary representation $\pi$ of
$\mathcal P_+^\uparrow (1)$ on $H\,$, uniquely determined by any
two of the relations
\begin{equation*} \phantom{xxxxx}
\overline{\text{\rm d}\pi (X_1)}=i\,\log \Delta_M\, ,\quad
\overline{\text{\rm d}\pi (X_2)}=i\,\log \Delta_N\, ,\quad
\overline{\text{\rm d}\pi (X_3)}=i\, P\, ,
\end{equation*}
where $X_1\, , X_2\, , X_3\in\mathfrak p_+^\uparrow (1)$ are the
Lie algebra elements defined in $(\ref{gener})\,$.
\end{itemize}
\end{theorem}
\smallskip

\noindent {\bf Remark.} This theorem is a generalization of
Wiesbrock's statement, where $\varphi$ is assumed to be bounded.
\medskip

\noindent {\bf The strategy of our proof.} The proof will be given
in Section 7 in several steps.

$\,\,$(i) First we study $\Delta_N^{-it}\,\Delta_M^{it}$ by using
the Modular Extension Theorem which will be proved in Section 5.

$\,$(ii) Then we show that $\Delta_N^{-it}\,\Delta_M^{it}$ has a
strong operator limit $T$ for $t\longrightarrow -\infty\,$. For the
existence of the wave operator $T$ we use our generalization of the
Borchers Structure Theorem which will be proved in Section 6.

(iii) Using the above ingredients, we define an $so$-continuous
one-parameter family $U(s)\, ,\, s\in \mathbb R\,$, of unitaries.

$\,$(iv) We prove that the defined family $U(s)\, ,\, s\in
\mathbb R\,$, is a one-parameter group having positive generator $P$
and we verify that for it the statements (1) - (10) in Theorem
\ref{hsmi} hold.

$\,\,$(v) Using that the generator $P$ of $\,\mathbb R\ni s
\to U(s)$ satisfies (10), (\ref{sumLie}) will imply that $P$ is
the closure of the operator (\ref{generator}).

We also prove
the following result about the structure of half-sided modular
inclusions:
\smallskip

\begin{theorem}\label{type}
Under the assumptions and with the notation as in Theorem \ref{hsmi},
the following hold $:$
\begin{itemize}
\item[(1)] $\gamma_s (z)=z$ for $s\geq 0$ and $z\in Z(M)\,$, so
$Z\big(\gamma_s (M)\big) =Z(M)$ for all $s\geq 0\,$.
\smallskip
\item[(2)] There exists the greatest central projection $p$ of $M$
satisfying $M p =N p\,$. For a projection $e\in M$ we have
\smallskip

\hspace{2 cm}$e\leq p\;\Longleftrightarrow\; U(s)\, e =e$ for all
$s\in\mathbb R\,$,
\smallskip

\noindent while for a projection $e\in M^\varphi$ with $\gamma_s
(e)=e\, ,\, s\geq 0\,$, we have even
\smallskip

\hspace{2 cm}$e\leq p\;\Longleftrightarrow\; U(s)\, e\, J_M\, e\, J_M =
e\, J_M\, e\, J_M$ for all $s\in\mathbb R\,$.
\smallskip
\item[(3)] $\displaystyle \qquad M^\varphi\subset
\bigcap\limits_{s\geq 0} \gamma_s (M)\;\Longrightarrow\;
\bigcap\limits_{s\geq 0}\gamma_s (M) =\big\{ x\in M\, ;\,\gamma_s (x)
= x\, ,\, s\geq 0\big\}\,$.
\smallskip
\item[(4)] $\hspace{3.8 mm}\overline{\mathfrak M_\varphi\cap
M^\varphi}^{\, so}
= M^\varphi\;\Longrightarrow\; M^\varphi\subset
\big\{ x\in M\, ;\,\gamma_s (x)=x\, ,\, s\geq 0\big\}$
\smallskip

$\hspace{3.09 cm}\Longrightarrow\; M (1_H-p)$ and $N (1_H-p)$
are of type \rm{III}$_1$
\smallskip

\hspace{3.75 cm} whenever $p\neq 1_H\,$.
\end{itemize}
\end{theorem}
\smallskip

\noindent {\bf Remarks.} (1) is proved in the case of bounded
$\varphi$ in \cite{Bo5}, Theorem 2.4, and for the general case we
shall essentially repeat the same proof.

If $\varphi$ is bounded, then the equality $\displaystyle
\bigcap\limits_{s\geq 0}\gamma_s (M) = \{ x\in M\, ;\,\gamma_s (x)
= x\, ,\, s\geq 0\}$ in (3) follows from \cite{Lo1}, Corollary 2.2,
but for the proof of (3) in our setting we need a different method.

Finally, for bounded $\varphi\,$, the inclusion
$M^\varphi\subset\{ x\in M\, ;\,\gamma_s (x)=x\, ,\, s\geq 0\}$
and the type of $M\, (1_H -p)$ and $N\, (1_H -p)$ were established
in \cite{Wie1} if $M$ is a factor, and in \cite{Bo4}, \cite{Bo5}
in the case of a general $M\,$. However, there is a gap in the
proof in \cite{Bo4}, \cite{Bo5}: it is shown only that, for every
$e\in M^\varphi$ majorized by $1_H -p\,$, the spectrum of
$\,\Delta_\varphi\, |\, e H$ is $[0,+\infty)\,$, while the right
proof requires that the spectrum of the modular operator of
$\varphi_e\,$, which is $\,\Delta_\varphi\, |\, eJ_\varphi
eJ_\varphi H$ (see (\ref{red.mod})), be equal to $[0,+\infty)\,$.

We don't know if the above inclusion still holds without assuming
the strong operator density of
$\mathfrak M_\varphi\cap M^\varphi$ in $M^\varphi\,$.
We notice that, by (3) and (4) in the above theorem, if
$\displaystyle M^\varphi\subset\bigcap\limits_{s\geq 0}\gamma_s (M)$
would hold in general, then we would always have:
\smallskip

\noindent\hspace{2.65 cm}$\displaystyle M^\varphi\subset \big\{ x
\in M\, ;\,\gamma_s (x) =x\, ,\, s\geq 0\big\} =
\bigcap\limits_{s\geq 0}\gamma_s (M)\,$,
\smallskip

\noindent\hspace{2.1 cm}$M (1_H-p)\; ,\, N (1_H-p)$ are of type
III$_1$ in the case $p\neq 1_H\,$.
\smallskip

\bigskip
\noindent (c) {\sc The Analytic Extension Theorem}.
\bigskip

Let $\beta\in\mathbb R\, ,\,\beta\neq 0\,$. Set
\begin{equation*}
\mathbb S_\beta =\{ z\in\mathbb C\, ;\, 0<\beta^{-1}\,\Im z <1\}\, .
\end{equation*}
$H^\infty (\mathbb S_\beta )$ will denote the Banach algebra
of all bounded analytic complex functions on $\mathbb S_\beta\,$.

Most parts of the following Analytic Extension Theorem are known, but
we shall give a proof for the convenience of the reader:
\smallskip

\begin{theorem}\label{anal.ext}
{\rm (Analytic Extension Theorem)}
Let $A$ and $B$ be invertible positive selfadjoint linear operators
on the Hilbert spaces $H$ and $K$ respectively, $0\neq \beta\in
\mathbb R\,$, and $T\in B(K,H)\,$. Then the next statements
$(1) - (5)$ are equivalent {\rm :}
\begin{itemize}
\item[(1)] $\mathbb R\ni s\longmapsto A^{is}\, T\, B^{-is}$
has a uniformly bounded $so$-continuous extension
\begin{equation}\label{ext}
\overline{\mathbb S_\beta}\ni z\longmapsto T(z)\in B(K,H)
\end{equation}
which is analytic in $\mathbb S_\beta\,$.
\item[(2)] $\mathbb R\ni s\longmapsto A^{is}\, T\, B^{-is}$
has a $wo$-continuous extension $(\ref{ext})$ which is analytic in
$\mathbb S_\beta\,$.
\item[(3)] There exists a Borel set $\Xi_o\subset\mathbb R$ of non-zero
Lebesgue measure such that, for every $\xi\in K$ and $\eta\in H\,$,
there exists an $f_{\xi ,\eta}\in H^\infty (\mathbb S_\beta )$
satisfying
\begin{equation}\label{ae}
\lim_{0<t/\beta \to 0} f_{\xi ,\eta}(s+i\, t)=(A^{is}\, T\, B^{-is}\,
\xi\, |\,\eta )
\end{equation}
for almost all $s\in \Xi_o\,$.
\item[(4)] $A^{-\beta}\, T\, B^\beta$ is defined and bounded on a
core of $B^\beta\,$.
\item[(5)] {\rm Dom}$\, (A^{-\beta}\, T\, B^\beta )=
\,${\rm Dom}$\, (B^\beta)$ and $A^{-\beta}\, T\, B^\beta$ is bounded.
\end{itemize}
Moreover, if the above conditions are satisfied, then, for every
$z\in\overline{\mathbb S_\beta}\,$,
\begin{equation}\label{dom}
\text{\rm Dom}\, (A^{iz}\, T\, B^{-iz})=\,\text{\rm Dom}\, (B^{-iz})\, ,
\end{equation}
\begin{equation}\label{op}
A^{iz}\, T\, B^{-iz}\subset T(z)\, ,
\end{equation}
\begin{equation}\label{transl}
A^{it}\, T(z)\, B^{-it} =T(z+t)\, ,\qquad t\in\mathbb R\, .
\end{equation}
\end{theorem}
\smallskip

\noindent {\bf Remark.} A somewhat novel feature is the non-null
Borel set $\Xi_o$ in (3). We shall apply this Theorem in the case
$K=H\,$, $T=1_H\,$, $A=\Delta_M\,$, $B=\Delta_N$ and
$\beta =-\frac 12$ in the proof of the Modular Extension Theorem
(Theorem \ref{char.herm}).
\smallskip

We notice that, for $A$ and $B$ as in Theorem \ref{anal.ext},
$\mathbb R\ni s\longmapsto \alpha_s =A^{is}\cdot B^{-is}$ is an
$so$-continuous one-parameter group of linear isometries on
$B(K,H)\,$. If, for some $T\in B(K,H)$ and $z\in\mathbb C\,$, the
orbit $\mathbb R\ni s\longmapsto A^{is}\, T\, B^{-is}$ has an
$so$-continuous extension
\begin{equation*}
\overline{\mathbb S_{\Im z}}\ni\zeta\longmapsto T(\zeta )\in
B(K,H)\, ,
\end{equation*}
which is analytic in $\mathbb S_{\Im z}\,$, then we say that $T$
belongs to the domain of $\alpha_z$ and define $\alpha_z (T)=
T(z)\,$. According to Theorem \ref{anal.ext}, for $T\in B(K,H)$
the conditions
\begin{itemize}
\item[--] $\; T\in\,\text{Dom}\, (\alpha_z )\,$,
\item[--] $\; A^{i z}\, T\, B^{-i z}$ is defined and bounded on
a core of $B^{-i z}\,$,
\item[--] $\;\text{Dom}\, (A^{i z}\, T\, B^{-i z}) =\,\text{Dom}\,
(B^{-i z})$ and $A^{i z}\, T\, B^{-i z}$ is bounded
\end{itemize}
are equivalent and, if they are satisfied, then $\alpha_z (T)=
\overline{A^{i z}\, T\, B^{-i z}}\,$.
For a more detailed study of the analytic extension operator
$\alpha_z\,$, especially in the most relevant case $z=-i\,$,
when it is called the \emph{analytic generator} of the group
$\alpha\,$, we refer to \cite{C-Z}, \cite{Z1} and \cite{Z2}.
\smallskip

Let us point out that, in general, the five equivalent statements
in Theorem \ref{anal.ext} are not equivalent with
\begin{itemize}
\item[$(4')$] $A^{-\beta}\, T\, B^\beta$ is densely defined and
bounded.
\end{itemize}
Indeed:
\smallskip

\begin{proposition}\label{counterex.}
There exist an invertible positive selfadjoint linear operator $A$
on a Hilbert space $H$ and a unitary $v\in B(H)\,$, such that
\smallskip

\noindent\hspace{3 cm} $A^{-1} v A$ is densely defined
and bounded, but

\noindent\hspace{3 cm} $\text{\rm Dom}\, (A^{-1} v A)$ is
not a core for $A\,$.
\end{proposition}
\smallskip

\noindent {\bf Remark.} Actually the proof of Proposition
\ref{counterex.}, which will be given in Section 3, works to
prove the next more general statement:
\smallskip

{\it If $A$ and $B$ are invertible positive selfadjoint linear
operators on a non-zero Hilbert space, such that
\begin{equation*}
A^{it}B^{is}A^{-it} =e^{-its} B^{is}\, ,\qquad t\, ,\, s\in
\mathbb R\, ,
\end{equation*}
hence $A^{it}\underbrace{\exp\, (-B^\pi )}_{=\, b} A^{-it} =
\big(\exp\, (-B^\pi )\big)^{e^{-\pi t}}\,$, $t\in\mathbb R\,$,
and further,
\begin{equation*}
A^{it}\, b^{is} A^{-it} =b^{ise^{-\pi t}}\, ,\qquad t\, ,\, s
\in\mathbb R\, ,
\end{equation*}
then, for every $s>0\,$,
\begin{equation*}
\text{\rm Dom}\, (A^{-1}\, b^{is} A) =\text{\rm Dom}\, (A)
\text{ and }A^{-1}\, b^{is} A\subset b^{-is}\, ,
\end{equation*}
while

\noindent\hspace{2.45 cm} $A^{-1}\, b^{-is} A$ is densely defined
and $A^{-1}\, b^{-is} A\subset b^{is}\,$, but

\noindent\hspace{2.45 cm} $\text{\rm Dom}\, (A^{-1}\, b^{-is} A)$
is not a core of $A\,$.}
\medskip

Nevertheless, there are situations in which the above statement
$(4')$ is equivalent with the statements $(1) - (5)$ in Theorem
\ref{anal.ext}.
One of such situations occurs in \cite{T1}, Lemma 15.15 and
Theorem 15.3, namely:
\smallskip

{\it Let $M\subset B(H)$ be a von Neumann algebra, in standard
form with respect to a normal semifinite faithful weight $\varphi$
on $M$. Then, for any $a\in\mathfrak N_\varphi\,$, the following
are equivalent {\rm :}
\begin{itemize}
\item[(a)] $\,\Delta_\varphi^{-1/2} a\Delta_\varphi^{1/2}$ is
densely defined and $\|\Delta_\varphi^{-1/2}a\,\Delta_\varphi^{1/2}
\|\leq 1\,$,
\item[(b)] $\,\|\Delta_\varphi^{1/2} a^*\Delta_\varphi^{-1/2}\|
\leq 1\,$,
\item[(c)] $\,\text{\rm Dom}\, (\Delta_\varphi^{-1/2} a\,
\Delta_\varphi^{1/2}) =\text{\rm Dom}\, (\Delta_\varphi^{1/2})$
and $\|\Delta_\varphi^{-1/2} a\,\Delta_\varphi^{1/2}\|\leq 1\,$.
\end{itemize}}
\smallskip

Indeed, if (a) holds and $\eta\in\text{Dom}\, (\Delta_\varphi^{1/2}
a^*\Delta_\varphi^{-1/2})\,$, then
\begin{equation*}
|(\Delta_\varphi^{1/2} a^*\Delta_\varphi^{-1/2}\eta\, |\,\xi )|=
|(\eta\, |\,\Delta_\varphi^{-1/2} a\,\Delta_\varphi^{1/2}\xi )|\leq
\|\eta\|\,\|\xi\|\, ,\quad\; \xi\in
\text{Dom}\, (\Delta_\varphi^{-1/2} a\Delta_\varphi^{1/2})\, .
\end{equation*}
Since $\text{Dom}\, (\Delta_\varphi^{-1/2} a\,\Delta_\varphi^{1/2})$
is dense in $H\,$, we get $\|\Delta_\varphi^{1/2} a^*
\Delta_\varphi^{-1/2}\eta\|\leq\|\eta\|\,$.

Next, $\text{Dom}\, (\Delta_\varphi^{1/2}a^*\Delta_\varphi^{-1/2})$
always contains the core $\{ J_\varphi x_\varphi\, ;\, x\in
\mathfrak A_\varphi \}$ of $\,\Delta_\varphi^{-1/2}$. Indeed, if
$x\in\mathfrak A_\varphi\,$, then $xa\in\mathfrak A_\varphi$ and
so $a^*\Delta_\varphi^{-1/2}
J_\varphi x_\varphi =a^* S_\varphi x_\varphi =a^* (x^*)_\varphi
=\big( (x a)^*\big)_\varphi$ belongs to $\text{Dom}\, S_\varphi
=\text{Dom}\,\Delta_\varphi^{1/2}\,$.

Consequently, according to Theorem \ref{anal.ext}, (b) implies
that $a^*\in\text{Dom}\, (\sigma^\varphi_{-i/2})$ and
$\|\sigma^\varphi_{-i/2}(a^*)\|\leq 1\,$. But then $a\in
\text{Dom}\, (\sigma^\varphi_{i/2})$ and $\sigma^\varphi_{i/2}
(a)^* =\sigma^\varphi_{-i/2}(a^*)\,$, hence
$\|\sigma^\varphi_{i/2} (a)\| =\|\sigma^\varphi_{-i/2}(a^*)\|
\leq 1\,$. Using again Theorem \ref{anal.ext}, we obtain that
(c) holds.

Finally, the implication $(c)\,\Rightarrow\, (a)$ is trivial.

\bigskip
\noindent (d) {\sc Lebesgue continuity, Tomita algebras}.
\bigskip

Let $M\subset B(H)$ be a von Neumann algebra, in standard form with
respect to a normal semifinite faithful weight $\varphi$ on $M\,$.

The next lemma shows that $1_H$ can be approximated by particularly
regular elements of $\mathfrak M_\varphi$ with respect to the
$so$-topology:
\smallskip

\begin{lemma}\label{anal.a.u}
There is an increasing net $\{ a_\iota\}_\iota$ in
$\,\mathfrak M_\varphi\cap M^+$ such that, for any $\iota\,$,
the orbit $\mathbb R\ni s\longmapsto \sigma^\varphi_s (a_\iota )
\in M$ has an entire extension $\mathbb C\ni z\longmapsto
\sigma^\varphi_z (a_\iota )\in M$ and
\begin{itemize}
\item[--] $\;\sigma^\varphi_z (a_\iota )\in\mathfrak M_\varphi\,
,\,\sigma^\varphi_z (a_\iota )^* =\sigma^\varphi_{\overline{z}}
(a_\iota )\, ,\,
\| \sigma^\varphi_z (a_\iota )\|\leq e^{(\Im z)^2}$ for
all $\,\iota$ and $\, z\in\mathbb C\,$,
\item[--] $\; so-\lim\limits_\iota\, \sigma^\varphi_z (a_\iota )
=1_H\,$ for all $\, z\in\mathbb C\,$.
\end{itemize}
\end{lemma}
\smallskip

Nets $\{ a_\iota\}_\iota$ as in Lemma \ref{anal.a.u} (called in
\cite{Z3}, \S$\, 1\,$, regularizing nets for $\varphi$) will be
used to prove the following description of $\mathfrak N_\varphi\,$:
\smallskip

\begin{lemma}\label{char.N}
{\rm (1)} For $x\in M$ and $c\geq 0\,$,
\begin{equation*}
x\in\mathfrak N_\varphi\,\text{ and }\,\| x_\varphi\|\leq c\;
\Longleftrightarrow\,\| x\, J_\varphi\, y_\varphi\|\leq c\,\| y\|\,
\text{ for all }\, y\in\mathfrak M_\varphi\, .
\end{equation*}

{\rm (2)} For $x\in M$ and $\xi\in H\,$,
\begin{equation*} \;\;
x\in\mathfrak N_\varphi\,\text{ and }\, x_\varphi =\xi\,
\phantom{xx}
\Longleftrightarrow\, x\, J_\varphi\, y_\varphi =J_\varphi\, y\,
J_\varphi\,\xi\,\text{ for all }\, y\in\mathfrak M_\varphi\, .
\end{equation*}
\end{lemma}
\smallskip

Using the above lemma, we get immediately
\begin{equation}\label{integr}
f\in L^1(\mathbb R )\, ,\,
x\in\mathfrak N_\varphi\,\Longrightarrow \left\{\!
\begin{array}{l}\displaystyle
wo-\int_{\mathbb R} f(t)\,\sigma^\varphi_t (x)\, dt\in
\mathfrak N_\varphi\,\text{ and} \\ \displaystyle
\bigg( wo-\int_{\mathbb R}^{\substack{{}\\{}}} f(t)\,
\sigma^\varphi_t (x)\, dt \bigg)_{\!\! \varphi} =\int_{\mathbb R} f(t)\,
\Delta_\varphi^{i t}\, x_\varphi\, dt \\
\hspace{3.95 cm} =\hat{f}(\log \Delta_\varphi )\, x_\varphi\, ,
\end{array} \right.
\end{equation}
where $\hat{f}$ is the inverse Fourier transform of $f\,$:
\begin{equation}\label{inv.Fourier}
\hat{f}(\lambda )=\int_{\mathbb R} f(t)\, e^{i\,\lambda\, t}\, dt\, ,
\qquad \lambda\in\mathbb R\, .
\end{equation}
Indeed, by (\ref{j}) and (\ref{delta}) we have for every $y\in
\mathfrak N_\varphi$
\begin{align*}
\bigg( wo-\int_{\mathbb R} f(t)\,\sigma^\varphi_t (x)\, dt\bigg)
J_\varphi\, y_\varphi &=\int_{\mathbb R} f(t)\,\big( \sigma^\varphi_t
(x)\, J_\varphi\, y_\varphi\big)\, dt \\
&=\int_{\mathbb R} f(t)\,\big( J_\varphi\, y\, J_\varphi\,
\sigma^\varphi_t (x)_\varphi\big)\, dt \\
&=J_\varphi\, y\, J_\varphi\, \int_{\mathbb R} f(t)\,
\Delta_\varphi^{i t}\, x_\varphi\, dt\, ,
\end{align*}
so we can apply Lemma \ref{char.N}.(2) to $\,\displaystyle
wo-\int_{\mathbb R} f(t)\,\sigma^\varphi_t (x)\, dt$ and
$\displaystyle \,\int_{\mathbb R} f(t)\, \Delta_\varphi^{i t}\,
x_\varphi\, dt\,$.

If $\varphi$ is bounded, then the linear mapping $M\ni x\longmapsto
x_\varphi =x\,\xi_\varphi\in H$ is bounded, but its inverse is in
general not bounded. For unbounded $\varphi$ even $x\longmapsto
x_\varphi$ is not bounded. Nevertheless, both $\mathfrak N_\varphi
\ni x\longmapsto x_\varphi\in H$ and its inverse have a dominated
continuity property with respect to the $wo$-topology on $M$ and
the weak topology on $H\,$, called in \cite{Z1}, \S $\, 2\,$,
\emph{Lebesgue continuity}.
For the proof of Theorem \ref{hsmi} we need the following variant
of \cite{Z1}, Section 4.6, Propositions 1 and 2, concerning the
Lebesgue continuity of $x\longmapsto x_\varphi$ and
$x_\varphi\longmapsto x\,$:
\smallskip

\begin{proposition}\label{Leb}
Let $M\subset B(H)$ be a von Neumann algebra, in standard form with
respect to a normal semifinite faithful weight $\varphi$ on $M\,$,
and $\{ x_\iota\}_\iota\subset\mathfrak N_\varphi$ a net.

{\rm (1)} If $\, wo-\lim\limits_\iota x_\iota =x\in M$ and
$\,\sup\limits_\iota \| (x_\iota)_\varphi\| <\infty\,$, then
$x\in\mathfrak N_\varphi$ and $\, (x_\iota)_\varphi
\overset{\iota}{\longrightarrow} x_\varphi$ in the weak topology of
$H\,$.

{\rm (2)} If $\, (x_\iota)_\varphi \overset{\iota}{\longrightarrow}
\xi\in H$ in the weak topology of $H$ and $\,\sup\limits_\iota
\| x_\iota\| <\infty\,$, then there exists
$x\in\mathfrak N_\varphi$ such that $\, wo-\lim\limits_\iota
x_\iota =x$ and $x_\varphi =\xi\,$.
\end{proposition}
\smallskip

Let $\mathfrak T_\varphi$ denote the set of all
$x\in\mathfrak A_\varphi$ such that $\mathbb R\ni s\longmapsto
\sigma^\varphi_s (x)\in M$ has an entire extension $\mathbb C\ni
z\longmapsto\sigma^\varphi_z (x)\in M$ satisfying
$\sigma^\varphi_z (x)\in\mathfrak A_\varphi$ for all
$z\in\mathbb C\,$. Since
\begin{align*}
&x\, ,\, y\in\mathfrak T_\varphi\,\Longrightarrow\,
x\, y\in\mathfrak T_\varphi\,\text{ and }\,\sigma^\varphi_z (x\, y)
=\sigma^\varphi_z (x)\,\sigma^\varphi_z (y)\, ,\, z\in\mathbb C\, ,\\
&x\in\mathfrak T_\varphi\phantom{xx}\,\,
\Longrightarrow\,x^*\in\mathfrak T_\varphi\,
\text{ and }\,\sigma^\varphi_z (x^*)=\sigma^\varphi_{\overline z}
(x)^*\, ,\, z\in\mathbb C\, ,
\end{align*}
$\mathfrak T_\varphi$ is a $*$-subalgebra of $\mathfrak A_\varphi\,$,
called the (maximal) \emph{Tomita algebra} of $\varphi\,$.

In the next variant of \cite{S-Z}, 10.21, Corollary 1, certain
standard properties of the Tomita algebra $\mathfrak T_\varphi$
are formulated.
\smallskip

\begin{proposition}\label{Tom}
Let $M\subset B(H)$ be a von Neumann algebra, in standard form with
respect to a normal semifinite faithful weight $\varphi$ on $M\,$.
Then
\begin{equation}\label{intertwin}
x\in\mathfrak T_\varphi\, ,\, z\in\mathbb C\,\Longrightarrow\,
x_\varphi\in\,\text{\rm Dom}\, (\Delta_\varphi^{i z})\text{ and }
\,\sigma^\varphi_z (x)_\varphi =\Delta_\varphi^{i z}\, x_\varphi
\end{equation}
and, for every $y\in\mathfrak A_\varphi\,$, there exists a sequence
$\{ y_n\}_{n\geq 1}$ in $\mathfrak T_\varphi$ such that
\begin{itemize}
\item[--] $y_n\overset{so}{\longrightarrow} y$ and
$y_n^{\, *}\overset{so}{\longrightarrow} y_{}^*\,$,
\item[--] $(y_n)_\varphi\longrightarrow y_\varphi$ and
$(y_n^{\, *})_\varphi\longrightarrow (y_{}^*)_\varphi$ in the
norm-topology of $H\,$.
\item[--] $\| \sigma^\varphi_z (y_n )\|\leq e^{n (\Im z)^2}\| y\|$
for all $n\geq 1$ and $z\in\mathbb C\,$,
\item[--] $\| \Delta_\varphi^{i z}\, (y_n )_\varphi\|\leq
e^{n (\Im z)^2} \| y_\varphi\|\,$, $\| \Delta_\varphi^{i z}\,
(y_n^{\, *})_\varphi\|\leq e^{n (\Im z)^2} \| (y_{}^*)_\varphi\|$
for all $n\geq 1$ and $z\in\mathbb C\,$.
\end{itemize}
\end{proposition}

We notice that the set $\mathfrak S_\varphi$ of all
$x\in\mathfrak T_\varphi\,$, for which
\begin{equation*}
\|\sigma^\varphi_z (x)\|\leq e^{c(x)\,\Im z}\| x\|\, ,\;\; 
\| \Delta_\varphi^{iz}\, x_\varphi\|\leq e^{c(x)\,\Im z}\| x_\varphi\|
\, ,\qquad z\in\mathbb C
\end{equation*}
with $c(x)\geq 0$ a constant depending only on $x\,$, is a
$*$-subalgebra of $\mathfrak T_\varphi$ and for every
$y\in\mathfrak A_\varphi$ there exists a sequence $\{ y_n\}_{n\geq 1}$
in $\mathfrak S_\varphi$ such that
\begin{itemize}
\item[--] $y_n\overset{so}{\longrightarrow} y$ and
$y_n^{\, *}\overset{so}{\longrightarrow} y_{}^*\,$,
\item[--] $(y_n)_\varphi\longrightarrow y_\varphi$ and
$(y_n^{\, *})_\varphi\longrightarrow (y_{}^*)_\varphi$ in the
norm-topology of $H$
\end{itemize}
(see \cite{S-Z}, 10.22).

\bigskip
\noindent (e) {\sc Hermitian maps}.
\bigskip

Let $H\, ,\, K$ be Hilbert spaces and $M\subset B(H)\, ,\, N\subset B(K)$
von Neumann algebras, in standard form with respect to the normal
semifinite faithful weights $\varphi$ on $M$ and $\psi$ on $N\,$. An
essential role will be played by the fixed point real linear subspaces
of $K$ and $H$ under $S_\psi$ and $S_\varphi\,$, respectively:
\begin{equation*}
K^{S_\psi} =\{\xi\in\,\text{Dom}\, (S_\psi )\, ;\, S_\psi\,\xi =\xi\}
\, ,\quad H^{S_\varphi} =\{\eta\in\,\text{Dom}\, (S_\varphi )\, ;\,
S_\varphi\,\eta =\eta\}\, .
\end{equation*}
They have been used by various authors earlier: see, for example,
\cite{Ri-VD} and \cite{Dav}.

Let us formulate the basic properties, for example, of $K^{S_\psi}\,$:
\smallskip

\begin{lemma}\label{fixed}
{\rm (1)} $\; K^{S_\psi} =\overline{\{ x_\psi\, ;\, x^*=x\in
\mathfrak N_\psi\}}\,$.

{\rm (2)} $\; \xi\in K$ belongs to $K^{S_\psi}$ if and only if
$\, (\xi\, |\, J_\psi\, x_\psi )\in\mathbb R\,$ for all $x^*=x\in
\mathfrak N_\psi\,$.

{\rm (3)} $\; \xi\in K$ belongs to $K^{S_\psi}$ if and only if
$\, (\xi\, |\, J_\psi\, x_\psi ) =(J_\psi\, (x^*)_\psi\, |\,\xi )\,$
for all \\ \phantom{xxxxxx} $x\in\mathfrak A_\psi\,$.

{\rm (4)} $\;${\rm Dom}$\, S_\psi =K^{S_\psi} +i\, K^{S_\psi}\,$.
\end{lemma}
\smallskip

\begin{definition}\label{def.herm}
{\rm (1)} $\, T\in B(K,H)$ is said to be Hermitian with respect to the
weight pair $(\psi ,\varphi )$ if
\begin{equation*}
T\, K^{S_\psi}\subset H^{S_\varphi}\, .
\end{equation*}

{\rm (2)} $\; T\in B(K,H)$ is said to implement $\psi$ in $\varphi$ if
\begin{equation*}
x\in \mathfrak N_\psi\;\Longrightarrow\; T\, x\, T^*\in
\mathfrak N_\varphi\, ,\; (T\, x\, T^*)_\varphi =T\, x_\psi\, .
\end{equation*}
\end{definition}
\smallskip

Statement (3) in the next lemma explains why we call the fulfilment
of the implication in Definition \ref{def.herm}.(2) ``implementation
of $\psi$ in $\varphi$ by $T$''.
\smallskip

\begin{lemma}\label{herm}
{\rm (1)} $\; T\in B(K,H)$ is Hermitian with respect to $(\psi ,\varphi )$
whenever it \\ \phantom{xxxxxx} implements $\psi$ in $\varphi\,$.

{\rm (2)} $\;$If $T\in B(K,H)$ implements $\psi$ in $\varphi\,$, then
$T\, N\, T^*\subset M\,$.

{\rm (3)} $\;$If an isometric $T\in B(K,H)$ implements $\psi$ in
$\varphi\,$, then $N\ni x\mapsto T\, x\, T^*\in M$ \\
\phantom{xxxxxx} is an injective $*$-homomorphism and
\begin{equation*}
\psi (a) =\varphi (T\, a\, T^*)\, ,\qquad 0\leq a\in\mathfrak M_\psi\, .
\end{equation*}

{\rm (4)} $\;$For bounded $\psi$ and $\varphi$ and the corresponding
cyclic and separating vectors \\ \phantom{xxxxxx} $\xi_\psi =(1_K)_\psi$
and $\eta_\varphi =(1_H)_\varphi\,$, an injective $T\in B(K,H)$
implements $\psi$ in $\varphi$ \\ \phantom{xxxxxx} if and only if
\begin{equation*}
T\, N\, T^*\subset M\,\text{ and }\, T^*\,\eta_\varphi =\xi_\psi\, .
\end{equation*}
\end{lemma}
\smallskip

The following result provides important criteria for Hermiticity:
\smallskip

\begin{theorem}\label{char.herm}
{\rm (Modular Extension Theorem)}
Let $M\subset B(H)\, ,\, N\subset B(K)$ be von Neumann algebras, in
standard form with respect to the normal semifinite faithful weights
$\varphi$ on $M$ and $\psi$ on $N\,$. Then for $T\in B(K,H)$ the
following conditions $(1) - (8)$ are equivalent $:$
\begin{itemize}
\item[(1)] $T$ is Hermitian with respect to $(\psi ,\varphi )\, ;$
\item[(2)] $T\, x_\psi\in H^{S_\varphi}\,$ for all $\, x^*=x\in
\mathfrak N_\psi\, ;$
\item[(3)] $(T\, x_\psi\, |\, J_\varphi\, y_\varphi )\in\mathbb R\,$
for all $\, x^*=x\in\mathfrak N_\psi\,$ and $\, y^*=y\in
\mathfrak N_\varphi\, ;$
\item[(4)] For every $x\in\mathfrak A_\psi\,$ and $\, y\in
\mathfrak A_\varphi\,$ we have
\begin{equation*}
(T\, x_\psi\, |\, J_\varphi\, y_\varphi )=(J_\varphi\, (y^*)_\varphi
\, |\, T\, (x^*)_\psi )\, ;
\end{equation*}
\item[(5)] $T\, S_\psi\subset S_\varphi\, T\, ;$
\item[(6)] $\Delta_\varphi^{1/2}\, T\,\Delta_\psi^{-1/2}$ is
defined on $\,${\rm Dom}$\;\Delta_\psi^{-1/2}$ and coincides there
with $\, J_\varphi\, T\, J_\psi\, ;$
\item[(7)] $J_\psi\, T^* J_\varphi\,$ is Hermitian with respect to
$(\varphi ,\psi )\, ;$
\item[(8)] {\rm (Modular Extension Condition)} $\mathbb R\ni s
\longmapsto\Delta_\varphi^{is}\, T\,\Delta_\psi^{-is}\in B(K,H)$
extends to a bounded $\, so$-continuous map
\begin{equation*}
\overline{\mathbb S_{-1/2}}\ni z\longmapsto T(z)\in B(K,H)\, ,
\end{equation*}
analytic in $\mathbb S_{-1/2}$ and satisfying
\begin{equation}\label{-i/2}
T\Big( -\frac i2 \Big) =J_\varphi\, T\, J_\psi\, .
\end{equation}
\end{itemize}
Moreover, if the above equivalent conditions are satisfied, then,
with the notation from the {\rm Modular Extension Condition (8)},
we have
\begin{equation}\label{norm}
\| T(z)\|\leq\| T\|\, ,\qquad z\in \overline{\mathbb S_{-1/2}}\, ,
\frac {}{}
\end{equation}
\begin{equation}\label{mod.transl}
T(z+t)=\Delta_\varphi^{it}\, T(z)\,\Delta_\psi^{-it}\, ,\qquad
z\in \overline{\mathbb S_{-1/2}}\, ,\, t\in\mathbb R\, ,
\end{equation}
\begin{equation}\label{transl-i/2}
T\Big( s-\frac i2\Big) =J_\varphi\, T(s)\, J_\psi\, ,\qquad
s\in\mathbb R
\end{equation}
and $T(s)$ is Hermitian with respect to $(\psi ,\varphi )$ for
all $s\in\mathbb R\,$.
\end{theorem}

\bigskip
\noindent (f) {\sc Generalization of the structure theorem of
Borchers}.
\bigskip

Let $M\, ,\, N\, ,\,\varphi\, ,\,\psi$ be as in the preceding
subsection, and $T\in B(K,H)$ Hermitian with respect to
$(\psi ,\varphi )\,$. Then, by Theorem \ref{char.herm}, the
orbit $\mathbb R\ni s\longmapsto\Delta_\varphi^{is}\, T\,
\Delta_\psi^{-is}$ of $T$ has a bounded $so$-continuous extension
$T(\,\cdot\, )$ to $\overline{\mathbb S_{-1/2}}\,$, analytic in
$\mathbb S_{-1/2}\,$, which satisfies the boundary conditions
\begin{center}
$T(s)$ is Hermitian with respect to $(\psi ,\varphi )$
for all $s\in \mathbb R\,$,
\end{center}
\begin{center}
$\displaystyle J_\varphi\, T\Big(s-\frac i2\Big)\, J_\psi =T(s)$
is Hermitian with respect to $(\psi ,\varphi )$
for all $s\in \mathbb R\,$.
\end{center}
The next extension of a structure theorem of H. J. Borchers
(\cite{Bo3}, Theorem B, see also \cite{Bo2}, Theorem 11.9 and
\cite{Wie1}, Theorem 2) shows, in particular, that also the converse
statement holds, that is any bounded $so$-continuous map
$\overline{\mathbb S_{-1/2}}\longrightarrow B(K,H)\,$, which is
analytic in $\mathbb S_{-1/2}$ and satisfies the above boundary
conditions, arises from a Hermitian $T\in B(K,H)$ as above.
\smallskip

\begin{theorem}\label{borchers}
{\rm (Generalized Structure Theorem)}
Let $M\subset B(H)$ and $N\subset B(K)$ be von Neumann algebras, in
standard form with respect to the normal semifinite faithful weights
$\varphi$ on $M$ and $\psi$ on $N\,$. Further let $0\neq\beta\in
\mathbb R\,$, $\Xi_o$ and $\Xi_1$ be Lebesgue null sets in
$\mathbb R\,$, and
\begin{equation*}
\overline{\mathbb S_\beta}\;\backslash\,\big(\Xi_o\cup (\Xi_1 +i\beta)
\big)\ni z\longmapsto T(z)\in B(K,H)
\end{equation*}
be a bounded map which is analytic in $\mathbb S_\beta$ and satisfies
the boundary conditions
\begin{itemize}
\item[(i)] $T(s)$ is Hermitian with respect to $(\psi ,\varphi )$
for all $s\in\mathbb R\,\backslash\,\Xi_o$ and
\begin{equation}\label{cont+}
T(s)=wo - \lim_{0<t/\beta \to 0} T(s+it)\, ,\qquad
s\in\mathbb R\,\backslash\,\Xi_o\, ,
\end{equation}
\item[(ii)] $J_\varphi\, T(s+i\beta )\, J_\psi$ is Hermitian with
respect to $(\psi ,\varphi )$ for all $s\in\mathbb R\,\backslash
\,\Xi_1$ and
\begin{equation}\label{cont-}
T(s+i\beta )=wo - \lim_{1>t/\beta \to 1} T(s+it)\, ,\qquad
s\in\mathbb R\,\backslash\,\Xi_1\, .
\end{equation}
\end{itemize}
Then, for some $\, T\in B(K,H)$ which is Hermitian with
respect to $(\psi ,\varphi )\,$,
\begin{equation}\label{formula}
T(s)=\Delta_\varphi^{\textstyle -i\frac s{2\beta}}\, T\,
\Delta_\psi^{\textstyle i\frac s{2\beta}}\, ,\qquad s\in\mathbb R\,
\backslash\,\Xi_o\, .
\end{equation}
Hence the given map $z\mapsto T(z)$ extends to an $so$-continuous map
on the whole $\overline{\mathbb S_\beta}$ and, with
the same notation $T(\,\cdot\, )$ for the extension, it satisfies
\begin{equation}\label{borch.transl}
T(z+2\beta t)=\Delta_\varphi^{-it}\, T(z)\,\Delta_\psi^{it}\, ,
\qquad z\in\overline{\mathbb S_\beta}\, ,\, t\in\mathbb R\, ,
\end{equation}
\begin{equation}\label{transl+ibeta}
T(s+i\beta )= J_\varphi\, T(s)\, J_\psi\, ,\qquad s\in\mathbb R\, .
\;\phantom{xx}
\end{equation}
\end{theorem}
\smallskip

\noindent {\bf Remark.} Our theorem owes much to Borchers' work,
its proof being based on the main idea of the proof of Theorem B in
\cite{Bo3}. Nevertheless, our approach has several features of
generality:

(a) $\, z\longmapsto T(z)$ is not assumed to be $so$-continuous on the
whole $\overline{\mathbb S_\beta}\,$, but only the existence of radial
limits are assumed almost everywhere on the boundary. In our
application to the proof of Theorem \ref{hsmi} we shall use
Theorem \ref{borchers} with $\Xi_o=\{ 0\}$ and $\Xi_1=\emptyset\,$.

(b) We are considering the case of arbitrary normal semifinite 
faithful weights $\varphi$ and $\psi\,$, without assuming their
boundedness.

(c) On the boundary we assume only the Hermiticity of $T(s)$ and
$J_\varphi\, T(s+i\beta )\, J_\psi$ rather than the implementation
of $\psi$ in $\varphi$ by these operators. The advantage of our
assumption consists in its linearity, which allows ``mollification'',
while Borchers' assumption is of quadratic nature, more difficult
to handle.

(d) Our proof is made more elementary, avoiding most arguments of
the two-dimensional complex analysis and using instead of the
Malgrange-Zerner Theorem only the elementary Osgood Lemma (the
Hartogs Theorem for continuous functions) along with the Morera
Theorem (one-dimensional edge-of-the-wedge theorem).

\bigskip
\noindent (g) {\sc Complements to the implementation theorem
of Borchers}.
\bigskip

Based on the ideas from \cite{Arv}, an invariant subspace theory
was developed in \cite{Z2} for the ``bounded analytic'' elements
associated to an $so$-continuous one-parameter group
$(\alpha_t )_{t\in\mathbb R}$
of $*$-automorphisms of a von Neumann algebra $M\subset B(H)\,$.
This theory allows, starting with an already existent one-parameter
group of unitaries on $H$ which implements $\alpha\,$, to construct
canonically a new implementing group of unitaries on $H\,$, which
has a minimality property and inherits certain properties of the
$*$-automorphism group $\alpha$ (see \cite{Arv}, Proposition in
Section 3, where the idea is formulated in the realm of a particular
situation, and \cite{Z2}, Theorem 5.3, Corollary 5.4, Lemma 5.11
for the general theory).

The above method yields a proof for
the one-parameter version of the celebrated implementation theorem
of Borchers \cite{Bo1}, claiming the innerness of $\alpha$ whenever
it is implemented by a one-parameter group of unitaries $\big(
U(s)\big)_{s\in\mathbb R}$
having positive generator (see \cite{Arv}, Theorem 3.1 and \cite{Z2},
Corollary 5.7). Moreover, as we shall see in the next theorem, the
obtained canonical inner implementing group of unitaries inherits
certain commutation properties of the $*$-automorphism group
$\alpha\,$.

We recall that, if $M$ is a von Neumann algebra and
$(\alpha_s )_{s\in\mathbb R}$ is an $so$-continuous one-parameter
group of $*$-automorphisms of $M\,$, then the \emph{spectral subspace}
of $\alpha$ corresponding to a closed set $F\subset\mathbb R$
is defined by
\begin{equation*}
M^\alpha (F) =\Big\{\, x\in M\, ;\, wo-\int_{\mathbb R} f(s)\,
\alpha_s(x)\, ds =0\,\text{ if }f\in L^1(\mathbb R )\, ,\, F\cap\,
\text{supp} (\hat{f}) =\emptyset\,\Big\}\; ,
\end{equation*}
where $\hat{f}$ denotes the inverse Fourier transform
(\ref{inv.Fourier}) of $f$ (see \cite{Arv}, Definition 2.1).
\smallskip

\begin{theorem}\label{pos.gen}
Let $M\subset B(H)$ be a von Neumann algebra and $P$ a selfadjoint
operator in $H\,$, such that $P$ is bounded below and
$\text{\rm Ad}\,\exp (isP)$ leaves $M$ invariant for all
$s\in\mathbb R\,$, defining thus an $so$-continuous one-parameter
group $(\alpha_s )_{s\in\mathbb R}$ of $*$-automorphisms of $M\,$.
Then there exists a unique injective $b\in M\, ,\, 0\leq b\leq 1_H\,$,
such that
\begin{itemize}
\item[(i)] $\alpha_s (x)=b^{-i s} x\, b^{i s}\, ,\quad
s\in\mathbb R\, ,\, x\in M\,$,
\item[(ii)] for any injective $d\in M\, ,\, 0\leq d\leq 1_H\,$,
such that the implementation relation $\alpha_s (x)=d^{-i s} x\,
d^{i s}\, ,\, s\in\mathbb R\, ,\, x\in M$ holds, we have
\begin{equation*} \phantom{xx}
\chi_{\substack{{}\\{(0,e^\lambda ]}}}(b) \leq
\chi_{\substack{{}\\{(0,e^\lambda ]}}}(d)\, ,\qquad
\lambda\in\mathbb R\, ,
\end{equation*}
where $\chi_{\substack{{}\\{(0,e^\lambda ]}}}$ stands for the
characteristic function of $(0,e^\lambda ]\,$.
\end{itemize}
Moreover,
\begin{itemize}
\item[(iii)] for every $\lambda\in\mathbb R\,$,
$\chi_{\substack{{}\\{(0,e^\lambda ]}}}(b)$ is the orthogonal
projection onto
\begin{equation*}
\bigcap_{\mu >\lambda}\;\text{the closed linear span of }
M^\alpha\big( [-\mu ,+\infty )\big)\, H\, ,
\end{equation*}
\item[(iv)] for any $*$-automorphism $\sigma$ of $M$ and
$\lambda_\sigma >0\,$, such that $\sigma\circ\alpha_s =
\alpha_{\substack{{}\\ {\lambda_\sigma s}}}\circ\sigma$ for all
$s\in\mathbb R\,$, we have $\sigma (b) =b^{\lambda_\sigma}\,$.
\end{itemize}
\end{theorem}
\smallskip

The above theorem will be used in the proof of Theorem \ref{type}.

\bigskip
\noindent (h) {\sc Summary of the remaining part of the paper}.
\bigskip

The remainder of this paper presents proofs for the above results:
\begin{itemize}
\item[--] Theorem \ref{anal.ext} and Proposition \ref{counterex.}
in Section 3,
\item[--] Lemma \ref{anal.a.u}, Lemma \ref{char.N},
Proposition \ref{Leb} and Proposition \ref{Tom} in Section 4,
\item[--] Lemma \ref{fixed}, Lemma \ref{herm} and Theorem
\ref{char.herm} in Section 5,
\item[--] Theorem \ref{borchers} in Section 6,
\item[--] Theorem \ref{hsmi} in Section 7 and, finally,
\item[--] Theorem \ref{pos.gen}
and Theorem \ref{type} in Section 8.
\end{itemize}

\bigskip
\section{The Analytic Extension Theorem}
\bigskip

The aim of this section is to prove Theorem \ref{anal.ext} and
Proposition \ref{counterex.}.
\medskip

{\it Proof of }$\,${\bf Theorem \ref{anal.ext}}.
\medskip

The equivalence of conditions (2), (4) and (5), as well as the
three additional statements (\ref{dom}), (\ref{op}), (\ref{transl})
were proved in \cite{C-Z}, Theorem 6.2.

For the proof of the remaining part, we introduce the following
notation: let $K_c(B)$ and $H_c(A)$ be the set of all vectors
$\xi\in K$ and $\eta\in H\,$, respectively, with compact spectral
support for $\log\, B$ and $\log\, A\,$, respectively. For such $\xi$
and $\eta\,$, $\mathbb C\ni z\longmapsto B^{iz}\,\xi$ and
$\mathbb C\ni z\longmapsto A^{iz}\,\eta$ are analytic functions of
exponential type with respect to Im$\, z$ and they are uniformly
bounded in $\overline{\mathbb S_\beta}\,$. Furthermore,
$K_c(B)\subset K$ and $H_c(A)\subset H$ are dense linear subspaces
and they are cores of $B^{iz}$ and $A^{iz}$ for every
$z\in\mathbb C\,$, respectively.
\medskip

{\it Proof of }(2)$\,\Rightarrow\,$(1). The uniform boundedness of
$\overline{\mathbb S_\beta}\ni z\longmapsto T(z)$ and its
$so$-continuity are to be proved. The latter is automatic on
$\mathbb S_\beta\,$, where $T(\,\cdot\, )$ is analytic.

Let $\xi\in K_c(B)$ and $\eta\in H_c(A)\,$. Then
\begin{equation*}
(T(z)\,\xi\, |\,\eta ) =(T\, B^{-iz}\,\xi\, |\, A^{-i\,\overline z}\,
\eta )\, ,\qquad z\in \overline{\mathbb S_\beta}\, ,
\end{equation*}
because the analytic function $\mathbb C\ni z\longmapsto
(T\, B^{-iz}\,\xi\, |\, A^{-i\,\overline z}\,\eta )$ and the
continuous function $\overline{\mathbb S_\beta}\ni z\longmapsto
(T(z)\,\xi\, |\,\eta )\,$, which is analytic in the interior,
coincide on $\mathbb R\,$. By (\ref{dom}) and by the density of
$H_c(A)$ in $H\,$, it follows that
\begin{equation}\label{expl}
T(z)\,\xi =A^{i\, z}\,T\, B^{-iz}\,\xi\, ,\qquad
z\in \overline{\mathbb S_\beta}\; ,\;\xi\in K_c(B)\, .
\end{equation}

Since $\| A^{iz}\, T\, B^{-iz}\| =\| T\|$ for $z\in\mathbb R$ and
$\| A^{iz}\, T\, B^{-iz}\| =\| A^{-\beta}\, T\, A^\beta\|$ for
$z\in\mathbb R +i\,\beta$ by (5), we have by the Three Line Theorem
\begin{equation*}
|(T(z)\,\xi\, |\,\eta )|\leq\max\,\big(\,\| T\|\, ,\,\| A^{-\beta}\,
T\, A^\beta\|\,\big)\,\|\xi\|\,\|\eta\|\, ,\qquad z\in
\overline{\mathbb S_\beta}\, ,
\end{equation*}
obtaining thus the uniform boundedness of $\overline{\mathbb S_\beta}\ni
z\longmapsto T(z)\,$.

Due to this uniform boundedness, it suffices to prove the convergences
\begin{equation}\label{str.op.1}
\lim_{\overline{\mathbb S_\beta}\,\ni z\to s} \| T(z)\,\xi -
A^{is}\, T\, B^{-is}\,\xi\| =0\, ,\qquad s\in\mathbb R\, ,
\end{equation}
\begin{equation}\label{str.op.2}
\lim_{\overline{\mathbb S_\beta}\,\ni z\to s+i\,\beta} \| T(z)\,\xi -
A^{is-\beta}\, T\, B^{-is+\beta}\,\xi\| =0\, ,\qquad s\in\mathbb R
\end{equation}
for $\xi\in K_c(B)\,$. We give a proof explicitly only for
$\beta >0\,$, the treatment of the case $\beta <0$ being completely
similar.

Let $E$ denote the spectral projection of $\log\, A$ corresponding
to $(-\infty ,0]\,$. Owing to (\ref{expl}) we can split $T(z)\,\xi$
as follows:
\begin{equation*}
T(z)=\big( A^{iz}\, (1_H -E)\big)\, T\, B^{-iz}\,\xi +
\big( A^{iz+\beta}\, E\big)\, (A^{-\beta}\, T\, B^\beta )\,
B^{-iz-\beta}\,\xi\, ,\qquad z\in \overline{\mathbb S_\beta}\, .
\end{equation*}
We note that
\begin{center}
$A^{iz}\, (1_H -E)$ is defined on $H$ and $\| A^{iz}\, (1_H -E)\|
\leq 1$ for all $z\in\mathbb C\, ,\;$Im$\, z\geq 0\,$,
\end{center}
\hspace{0.2 cm} $A^{iz+\beta}\, E$ is defined on $H$ and
$\| A^{iz+\beta}\, E\|\leq 1$ for all $z\in\mathbb C\, ,\;$Im$\, z
\leq\beta\,$.

\noindent Now the norm-continuity of $\overline{\mathbb S_\beta}\ni
z\mapsto B^{-iz}\,\xi$ and $\overline{\mathbb S_\beta}\ni z\mapsto
B^{-iz-\beta}\,\xi\,$, the $so$-continuity of
$\overline{\mathbb S_\beta}\ni z\mapsto A^{iz}\, (1_H -E)$ and
$\overline{\mathbb S_\beta}\ni z\mapsto A^{iz+\beta}\, E\,$, and
the boundedness of $A^{-\beta}\, T\, B^\beta$ on $K_c(B)$ yield the
convergences (\ref{str.op.1}) and (\ref{str.op.2}).
\medskip

{\it Proof of }(1)$\,\Rightarrow\,$(3): obvious, with $\Xi_o =
\mathbb R\,$.
\medskip

{\it Proof of }(3)$\,\Rightarrow\,$(4). We proceed in three steps.
\smallskip

{\it Step 1}. First we quote some results from the theory of the
Hardy spaces on the disc. Let $H^\infty (\mathbb D )$ be the Banach
algebra of all bounded analytic complex functions on the unit disc
\begin{equation*}
{\mathbb D} =\{ z\in \mathbb C\, ;\, |z|<1\}\, .
\end{equation*}
Any $g\in H^\infty (\mathbb D )$ has a non-tangential limit
\begin{equation*}
\widetilde g (\zeta ) =\lim_{\mathbb D\,\ni z\to\zeta} g(z)
\end{equation*}
for almost all $\zeta$ in the boundary $\partial\,\mathbb D$ of
$\mathbb D$ (the unit circle) due to Fatou's Theorem (see, for
example, the second Corollary on page 38 of \cite{Ho} or the theorems
on pages 5 and 14 of \cite{Ko}). Furthermore, the map $H^\infty
(\mathbb D )\ni g\longmapsto\widetilde g\in L^\infty
(\partial\,\mathbb D)$ obtained this way is an isometric algebra
homomorphism. On the other hand, the range $\{\widetilde g\, ;\,
g\in H^\infty (\mathbb D )\}$ of the above homomorphism is equal to
\begin{equation*}
\{ \psi\in L^\infty (\partial\,\mathbb D)\; ;\, \int_0^{2\pi}
\psi (e^{is})\, e^{iks}\, ds =0\,\text{ for all }\, k=1\, ,\, 2
\, ,\,\ldots\,\}
\end{equation*}
and hence it is weak${}^*$closed (see, for example, \cite{Con},
\S$\,$20.1). We notice also that, according to a uniqueness theorem
of the Riesz brothers (see, for example, the second Corollary on
page 52 of \cite{Ho} or the Theorem on page 76 of \cite{Ko}),
if for some $g\in H^\infty (\mathbb D )$ the boundary function
$\widetilde g$ vanishes almost everywhere on a Borel subset of
$\partial\,\mathbb D$ with non-zero arc length measure, then
$g=0\,$.

We consider the one point compactification of the right half
and the left half of $\overline {\mathbb S_\beta}$ and denote
each added point by $+\infty$ and $-\infty\,$, respectively.
We extend the function
\begin{equation*}
\overline {\mathbb D}\,\backslash\,\{ +1\, ,\, -1\}\ni
\zeta\longmapsto \Phi_\beta (\zeta ) =\frac \beta\pi \,\log
\left( i\,\frac{1+\zeta}{1-\zeta}\right)\in
\overline {\mathbb S_\beta}
\end{equation*}
to be $+\infty$ at $\zeta =+1$ and $-\infty$ at $\zeta =-1\,$.
Then the extended function
\begin{equation*}
\overline {\mathbb D}\ni
\zeta\longmapsto \Phi_\beta (\zeta ) =\frac \beta\pi \,\log
\left( i\,\frac{1+\zeta}{1-\zeta}\right)\in
\overline {\mathbb S_\beta}\,\cup\{ -\infty\, ,\, +\infty\}
\end{equation*}
is a homeomorphism, mapping $\mathbb D$ onto $\mathbb S_\beta$
conformally, and the boundary $\partial\,\mathbb D$ onto
$\partial\,\mathbb S_\beta\cup\{ -\infty\, ,\, +\infty\}$ absolutely
bicontinuously with respect to the arc length measures: if $\Xi$
is a Borel set in $\partial\,\mathbb D\,$, then $\Xi$ has arc
length measure $0$ if and only if $\Phi_\beta (\Xi )$ has arc
length measure $0\,$. Moreover, $\Phi_\beta$ maps paths in
$\mathbb D$ tending to a $\zeta\in\partial\,\mathbb D\,\backslash\,
\{ 1,-1\}$ from within a sector of opening $<\pi$ having vertex
at $\zeta\,$, and symmetric about the inner normal to
$\partial\,\mathbb D$ in $\zeta\,$, to paths tending to
$\Phi_\beta (\zeta )\in\partial\,\mathbb S_\beta$ in a similar
non-tangential way. Therefore, if $f\in H^\infty
(\mathbb S_\beta )\,$, the non-tangential limit
\begin{equation*}
\widetilde f (\zeta ) =\lim_{\mathbb D\,\ni z\to\zeta} f(z)
\end{equation*}
exists for almost all $\zeta$ in $\partial\,\mathbb S_\beta$
by Fatou's Theorem applied to $f\circ\Phi_\beta\,$. Similarly
we can transcribe the above quoted results concerning
$H^\infty (\mathbb D)$ in the setting of
$H^\infty (\mathbb S_\beta)\,$:
$H^\infty (\mathbb S_\beta )\ni g\longmapsto\widetilde g\in
L^\infty (\partial\,\mathbb S_\beta)$ is an isometric algebra
homomorphism with weak${}^*$closed range in
$L^\infty (\partial\,\mathbb S_\beta)$ and
$f\in H^\infty (\mathbb S_\beta)$ is equal to zero whenever
$\widetilde f$ vanishes almost everywhere on a Borel subset of
$\partial\,\mathbb S_\beta$ with non-zero arc length measure.
\smallskip

{\it Step 2}. We consider the map
\begin{equation*}
F : K\times H\ni (\xi\, ,\,\eta )\longmapsto \widetilde f_{\xi ,\eta}
\in \{\widetilde f\, ;\, f\in H^\infty (\mathbb S_\beta )\}\subset
L^\infty (\partial\,\mathbb S_\beta )\, ,
\end{equation*}
where, as noticed in Step 1, $\{\widetilde f\, ;\, f\in H^\infty
(\mathbb S_\beta )\}$ is a weak${}^*$closed subalgebra of $L^\infty
(\partial\,\mathbb S_\beta )\,$. The function $f_{\xi ,\eta}$ is
uniquely determined by (\ref{ae}) due to the uniqueness result
quoted in Step 1. Since the right hand side of (\ref{ae}) is
sesquilinear in $\xi$ and $\eta\,$, the mapping $F$ is also
sesquilinear. We shall prove in this step that $F$ is bounded.

We first prove that the graph of $F$ is closed. Suppose
that  $\xi_n\longrightarrow\xi_o\,$, $\eta_n\longrightarrow\eta_o$
and $\widetilde f_{\xi_n ,\eta_n}\longrightarrow \widetilde f_o$
with respect to the norm of $L^\infty (\partial\,\mathbb S_\beta
)\,$, hence also $f_{\xi_n ,\eta_n}\longrightarrow f_o$ uniformly.
By the continuity of the right hand side of (\ref{ae}) in $\xi$
and $\eta\,$, $\widetilde f_o$ has to satisfy (\ref{ae}) for $\xi =
\xi_o$ and $\eta =\eta_o$ almost everywhere on the set $\Xi_o\,$.
Therefore $\widetilde f_o =\widetilde f_{\xi_o ,\eta_o}\,$, again
by the uniqueness theorem of the Riesz brothers, proving that
the graph of $F$ is closed.

Let us consider $F(\xi\, ,\,\cdot\, )$ for a fixed $\xi\,$. Denote
by $\overline H$ the conjugate of the Hilbert space $H$ and by
$\overline \eta$ the canonical image of $\eta\in H$ in
$\overline H\,$. By the above proved closedness of the graph of
$F\,$, the graph of the linear map ${\overline H}\ni{\overline \eta}
\longmapsto F(\xi\, ,\,\eta )\in
L^\infty (\partial\,\mathbb S_\beta )$ is closed and hence, by the
Closed Graph Theorem,
\begin{equation*}
\| F(\xi\, ,\,\eta )\|
\leq c_\xi\,\|\eta\|\, ,\qquad \eta\in H
\end{equation*}
for some constant $c_\xi\geq 0$ depending on $\xi\,$. Thus
$F(\xi\, ,\,\cdot\, ) : \overline H\longrightarrow
L^\infty (\partial\,\mathbb S_\beta )$ is bounded.

Now we prove that the graph of the linear map
\begin{equation*}
K\ni\xi\longmapsto F(\xi\, ,\,\cdot\, )\in
B\big(\, {\overline H}\, ,\,
L^\infty (\partial\,\mathbb S_\beta )\big)
\end{equation*}
is closed. Suppose that $\xi_n\longrightarrow\xi_o$ and
$F(\xi_n\, ,\,\cdot\, )\longrightarrow T_o$ with respect to the norm
of $B\big(\, {\overline H}\, ,\,
L^\infty (\partial\,\mathbb S_\beta )\big)\,$. Then, for every
$\eta\in H\,$, $F(\xi_n\, ,\,\eta )\longrightarrow T_o\,\eta$ and
by the closedness of the graph of $F$ it follows that $T_o\,\eta
= F(\xi_o\, ,\,\eta )\,$. Thus $T_o=F(\xi_o\, ,\,\cdot\, )\,$.

By the Closed Graph Theorem
\begin{equation*}
\| F(\xi\, ,\,\cdot\, )\|\leq c\,\|\xi\|\, ,\qquad \xi\in K
\end{equation*}
for some constant $c\geq 0\,$, so
\begin{equation*}
\| F(\xi\, ,\,\eta )\| =
\underset {\zeta\in\,\partial\,\mathbb S_\beta}{\text{ess sup}}\;
|\widetilde f_{\xi ,\eta}(\zeta )| =\sup_{z\in\, \mathbb S_\beta}
|f_{\xi ,\eta}(z)|\leq c\,\|\xi\|\,\|\eta\|\, ,\qquad \xi\in K\, ,\,
\eta\in H\, .
\end{equation*}
\smallskip

{\it Step 3}. We now take $\xi\in K_c(B)\, ,\,\eta\in H_c(A)\,$.
Then
\begin{equation*}
\mathbb C\ni z\longmapsto g_{\xi ,\eta}(z)=(T\, B^{-iz}\,\xi\, |\,
A^{-i\,\overline z}\,\eta )
\end{equation*}
is an entire function satisfying the boundary condition (\ref{ae}),
so that $g_{\xi ,\eta}=f_{\xi ,\eta}$ by the  uniqueness theorem of
the Riesz brothers. Therefore
\begin{equation*}
|g_{\xi ,\eta}(z)|\leq c\,\|\xi\|\,\|\eta\|\, ,\qquad
z\in\mathbb S_\beta
\end{equation*}
and hence the same estimate holds for all $z\in
\overline {\mathbb S_\beta}$ by continuity. This implies
\begin{equation*}
T\, B^{-iz}\,\xi\in\,\text{Dom}\, (A^{iz})\, ,
\end{equation*}
because $H_c(A)$ is a core of $A^{iz}\,$, and
\begin{equation*}
\| A^{iz}\, T\, B^{-iz}\,\xi\|\leq c\,\|\xi\|\, ,\qquad z\in
\overline {\mathbb S_\beta}\; .
\end{equation*}
Since $K_c(B)$ is a core of $B^{-iz}\,$, setting $z=i\,\beta$
we obtain (4).

\hfill $\square\quad$

\medskip

{\it Proof of }$\,${\bf Proposition \ref{counterex.}}.
\medskip

Let us denote:
\smallskip

\noindent\hspace{1.05 cm} by $\lambda_t$ the translation operator
$\xi\longmapsto\xi (\,\cdot\, -t)$ on $L^2(\mathbb R )\, ,\, t\in
\mathbb R\,$,
\smallskip

\noindent\hspace{1.05 cm} by $\alpha_t$ the $*$-automorphism
Ad$\, (\lambda_t )$ of $B(H)\, ,\, t\in\mathbb R\,$, and
\smallskip

\noindent\hspace{1.05 cm} by $b$ the multiplication operator with
$e^{-e^{\pi\,\cdot\,}}$on $L^2(\mathbb R )\,$.
\smallskip

\noindent Clearly, $(\alpha_t )_{t\in\mathbb R}$ is an
$so$-continuous one-parameter group of $*$-automorphisms of
$B(H)\,$, $0\leq b\leq 1_{L^2(\mathbb R )}$ and $b$ is injective.
Since
\begin{equation*} \phantom{x}
\alpha_t^{} (b) =\lambda_t^{}\, b\, \lambda_t^* =b^{\, e^{-\pi t}}
\, ,\qquad t\in\mathbb R\, ,
\end{equation*}
we have
\begin{equation}\label{orbit}
\alpha_t^{} (b^{is}) =\lambda_t^{}\, b^{is} \lambda_t^* =
b^{\, i s e^{-\pi t}}\, ,\qquad t\, ,\, s\in\mathbb R\, .
\end{equation}
By the Stone Representation Theorem there exists an invertible
positive selfadjoint linear operator $A$ on $L^2(\mathbb R )$
such that $\lambda_t =A^{it}\, ,\, t\in \mathbb R\,$. Then
\begin{equation}\label{Stone}
\alpha_t =\text{Ad}\, (A^{it})\, ,\qquad t\in\mathbb R\, .
\end{equation}

Let $s>0$ be arbitrary. Since $0\leq b\leq 1_{L^2(\mathbb R )}$ and
$\Re\, (ise^{-\pi z}) =s\, e^{-\pi \Re z} \sin\, (\pi \Im z )\,$,
the complex power $b^{i s e^{-\pi z}}\in B\big( L^2(\mathbb R )
\big)$ is defined and $\|\, b^{i s e^{-\pi z}}\|\leq 1$ for
every $z$ in the closed strip $\overline{\mathbb S_1}\,$.
Using (\ref{orbit}), it is easily seen that
\begin{equation*}
F_1 : \overline{\mathbb S_1}\ni z\longmapsto b^{i s e^{-\pi z}}
\in B\big( L^2(\mathbb R )\big)
\end{equation*}
is an $so$-continuous extension of $\,\mathbb R\ni t\longmapsto
\alpha_t (b^{is})\,$, which is analytic in $\mathbb S_1$ and
whose value at $i$ is $b^{-is}$. Taking into account (\ref{Stone}),
Theorem \ref{anal.ext} yields that
\begin{equation*}
\text{Dom}\, (A^{-1} b^{is} A) =\text{Dom}\, (A)\,\text{ and }\,
A^{-1} b^{is} A\subset b^{-is}\, ,
\end{equation*}
that is $A^{-1} b^{is} A =b^{-is}\, |\,\text{Dom}\, (A)\,$.
Therefore
\begin{equation}\label{yes}
\text{Dom}\, (A^{-1} b^{-is} A) =b^{-is}\big(\,\text{Dom}\, (A)
\big)\text{ is dense in } H\text{ and }\, A^{-1} b^{-is} A\subset
b^{is}\, .
\end{equation}
But
\begin{equation}\label{not}
\text{Dom}\, (A^{-1} b^{-is} A)\text{ is not a core of }A\, .
\end{equation}

Indeed, assuming that $\text{Dom}\, (A^{-1} b^{-is} A)$ is a core
of $A\,$, (\ref{yes}) and Theorem \ref{anal.ext} imply that
$\mathbb R\ni t\longmapsto \alpha_t (b^{-is})$ has a uniformly
bounded $so$-continuous extension
\begin{equation*}
F_2 : \overline{\mathbb S_1}\longrightarrow B\big( L^2(\mathbb R )
\big)\, ,
\end{equation*}
which is analytic in $\mathbb S_1$ and whose value at $i$ is
$b^{is}\,$. Then
\begin{equation*}
\overline{\mathbb S_{-1}}\ni z\longmapsto F_2 (\overline{z})^*\in
B\big( L^2(\mathbb R )\big)
\end{equation*}
is an $so$-continuous extension of $\,\mathbb R\ni t\longmapsto
\alpha_t (b^{is})\,$, which is analytic in $\mathbb S_{-1}$ and
whose value at $-i$ is $b^{-is}$. Consequently, $\,\mathbb R\ni
t\longmapsto\alpha_t (b^{is})$ has a uniformly bounded
$so$-continuous extension
\begin{equation*}
F : \overline{\mathbb S_1\cup \mathbb S_{-1}}\ni z\longmapsto
\left\{ \begin{array}{ll} F_1 (z) &\text{if }\,\Im z\geq 0 \\
F_2 (\overline{z})^* &\text{if }\,\Im z\leq 0 \end{array}\right.
\, ,
\end{equation*}
which is analytic in the interior and which takes the same value
$b^{-is}$ at $i$ and $-i\,$. Then, by (\ref{transl}), $F$ is periodic
of period $2\, i\,$, so it extends to a uniformly bounded entire
mapping, which must be constant by the Liouville Theorem. Thus
the orbit $\mathbb R\ni t\longmapsto\alpha_t (b^{is})$ is constant,
that is $b ^{is}$ commutes with every $\lambda_t\,$. Since $b ^{is}$
is the multiplication operator with $\,\mathbb R\ni r\longmapsto
e^{-ise^{\pi r}}$ on $L^2(\mathbb R )\,$, this means that the above
function is constant, what is plainly not true.

By (\ref{yes}) and by (\ref{not}) we conclude that, choosing
$v=b ^{-is}$ with $s>0\,$, $A^{-1} v A$ is densely defined and
bounded, but $\text{\rm Dom}\, (A^{-1} v A)$ is not a core for $A\,$.

\hfill $\square\quad$

\bigskip
\section{Lebesgue Continuity, Tomita Algebras}
\bigskip

In this section we prove Lemmas \ref{anal.a.u} and \ref{char.N},
as well as Propositions \ref{Leb} and \ref{Tom}. Throughout this
section $M\subset B(H)$ will stand for a von Neumann algebra, in
standard form with respect to a normal semifinite faithful weight
$\varphi$ on $M\,$.
\medskip

{\it Proof of }$\,${\bf Lemma \ref{anal.a.u}}.
\medskip

Since $\mathfrak M_\varphi$ is a hereditary $*$-subalgebra of $M\,$,
there is an increasing approximate unit $\{ b_\iota\}_\iota$ for
$\mathfrak M_\varphi$ (for example, the upward directed set
$\{ b\in\mathfrak M_\varphi\cap M^+\, ;\,\| b\| <1\}\,$, labeled by
itself). Then, by the $so$-density of $\mathfrak M_\varphi$ in
$M\,$, we have $\, so-\lim\limits_\iota\, b_\iota =1_H\,$. Setting
\begin{equation*}
a_\iota =\frac 1{\sqrt \pi}\; wo-\int_{-\infty}^\infty e^{-t^2}
\sigma^\varphi_t (b_\iota )\, dt\, ,
\end{equation*}
$\{ a_\iota\}_\iota$ is an increasing net in $M^+$ such that every
orbit
\begin{equation*}
\mathbb R\ni s\longmapsto \sigma^\varphi_s (a_\iota ) =
\frac 1{\sqrt \pi}\; wo-\int_{-\infty}^\infty e^{-(t-s)^2}
\sigma^\varphi_t (b_\iota )\, dt\in M
\end{equation*}
has an entire extension
\begin{equation*}
\mathbb C\ni z\longmapsto \sigma^\varphi_z (a_\iota ) =
\frac 1{\sqrt \pi}\; wo-\int_{-\infty}^\infty e^{-(t-z)^2}
\sigma^\varphi_t (b_\iota )\, dt\in M\, .
\end{equation*}
Clearly, $\,\sigma^\varphi_z (a_\iota )^* =
\sigma^\varphi_{\overline{z}} (a_\iota )$ for all $\iota$ and
$z\in\mathbb C\,$.
Since, for every $z\in\mathbb C\,$, the function
\begin{equation*}
\mathbb R\ni t\longmapsto
e^{-(t-z)^2}=e^{-(t-\Re z)^2 +(\Im z)^2} e^{2\, i\, (t-\Re z)\,\Im z}
\end{equation*}
is of the form $f_1-f_2+i\, (f_3-f_4)$ with $0\leq f_j\in
L^1(\mathbb R )\, ,\, 1\leq j\leq 4\,$, using (\ref{inv}) we deduce
easily that
\begin{equation*}
\sigma^\varphi_z (a_\iota )\in\mathfrak M_\varphi\,\text{ and }\,
\| \sigma^\varphi_z (a_\iota )\|\leq e^{(\Im z)^2}\,\text{ for all }
\,\iota\,\text{ and }\, z\in\mathbb C\, .
\end{equation*}
On the other hand, $\, so-\lim\limits_\iota\, b_\iota =1_H$
yields
\begin{equation*}
so-\lim\limits_\iota\, \sigma^\varphi_z (a_\iota ) =1_H\,\text{ for
all }\, z\in\mathbb C\, .
\end{equation*}

\hfill $\square\quad$

\medskip

{\it Proof of }$\,${\bf Lemma \ref{char.N}}.
\medskip

First we prove that
\begin{equation}\label{-1/2}
y\in\mathfrak A_\varphi\cap\text{Dom}\, (\sigma^\varphi_{-\frac i2})
\, ,\,\sigma^\varphi_{-\frac i2} (y)\in\mathfrak A_\varphi\;
\Longrightarrow\; \sigma^\varphi_{-\frac i2} (y)_\varphi =
J_\varphi\, (y^*)_\varphi\, .
\end{equation}
For let $x\in\mathfrak A_\varphi$ be arbitrary. Then (\ref{j}) yields
\begin{equation}\label{-1/2first}
\sigma^\varphi_{-\frac i2} (y)\, J_\varphi\, x_\varphi =
J_\varphi\, x\, J_\varphi\,\sigma^\varphi_{-\frac i2} (y)_\varphi\, .
\end{equation}
On the other hand, since $J_\varphi\, x_\varphi =J_\varphi\, S_\varphi
(x^*)_\varphi =\Delta_\varphi^{\frac 12}\, (x^*)_\varphi\in\,\text{Dom}
\, (\Delta_\varphi^{-\frac 12})\,$,
using Theorem \ref{anal.ext} we obtain
\begin{equation}\label{-1/2second}
\sigma^\varphi_{-\frac i2} (y)\, J_\varphi\, x_\varphi =
\Delta_\varphi^{\frac 12}\, y\,\Delta_\varphi^{-\frac 12}\, J_\varphi
\, x_\varphi =\Delta_\varphi^{\frac 12}\, (y\, x^*)_\varphi =
J_\varphi\, S_\varphi (y\, x^*)_\varphi =J_\varphi\, x\,
(y^*)_\varphi
\end{equation}
Now, (\ref{-1/2first}) and (\ref{-1/2second}) imply
\begin{equation*}
J_\varphi\, x\, J_\varphi\,\sigma^\varphi_{-\frac i2} (y)_\varphi =
J_\varphi\, x\, (y^*)_\varphi\,
,\qquad x\in\mathfrak A_\varphi
\end{equation*}
and by the $so$-density of $\,\mathfrak A_\varphi$ in $M\ni 1_H$
we conclude that $\,\sigma^\varphi_{-\frac i2} (y)_\varphi =
J_\varphi\, (y^*)_\varphi\,$.
\smallskip

(1) If $x\in\mathfrak N_\varphi\,$, then by (\ref{j})
\begin{equation*}
\| x\, J_\varphi\, y_\varphi\| =\| J_\varphi\, y\, J_\varphi\,
x_\varphi\|\leq \| y\|\,\| x_\varphi\|\, ,\qquad y\in
\mathfrak N_\varphi\supset\mathfrak M_\varphi\, .
\end{equation*}

Conversely, assume that $x\in M$ and $c\geq 0$ are such that
\begin{equation}\label{weak}
\| x\, J_\varphi\, y_\varphi\|\leq c\,\| y\|\, ,\qquad y\in
\mathfrak M_\varphi\, .
\end{equation}
Let $\{ a_\iota\}_\iota$ be a net as in Lemma \ref{anal.a.u}.
Then we have for every $\iota$
\begin{equation*}
\| \sigma^\varphi_{-\frac i2} (a_\iota )^* x^*x\,
\sigma^\varphi_{-\frac i2} (a_\iota )\|\leq \| x\|^2
\| \sigma^\varphi_{-\frac i2} (a_\iota )\|^2\leq \| x\|^2
e^{1/2}
\end{equation*}
and, according to (\ref{-1/2}) and (\ref{weak}),
\begin{equation*}
\varphi\big( \sigma^\varphi_{-\frac i2} (a_\iota )^* x^*x\,
\sigma^\varphi_{-\frac i2} (a_\iota )\big) =\|\, x\,
\sigma^\varphi_{-\frac i2} (a_\iota )_\varphi\|^2 =\|\, x\,
J_\varphi\, (a_\iota )_\varphi\|^2\leq c^2\,\| a_\iota \|^2\leq
c^2\, .
\end{equation*}
Since $\sigma^\varphi_{-\frac i2} (a_\iota )^* x^*x\,
\sigma^\varphi_{-\frac i2} (a_\iota ) =
\sigma^\varphi_{\frac i2} (a_\iota )\, x^*x\,
\sigma^\varphi_{-\frac i2} (a_\iota )
\overset{\iota}{\longrightarrow} x^*x$ in the $so$-topology and
$\varphi$ is lower $wo$-semi-continuous on the bounded subsets
of $M^+\,$, it follows that
\begin{equation*}
\varphi (x^*x)\leq c^2\, ,\text{ that is }\, x\in
\mathfrak N_\varphi\,\text{ and }\,\| x_\varphi\|\leq c\, .
\end{equation*}

(2) Since the implication ``$\Longrightarrow$'' is an immediate
consequence of (\ref{j}), we have to prove only the converse
implication. For let $x\in M$ and $\xi\in H$ be such that
\begin{equation}\label{strong}
x\, J_\varphi\, y_\varphi =J_\varphi\, y\, J_\varphi\,\xi\, ,
\qquad y\in\mathfrak M_\varphi\, .
\end{equation}
Then (\ref{weak}) holds with $c=\|\xi\|\,$, so by the above part
of the proof we have $x\in\mathfrak N_\varphi\,$. But then
(\ref{j}) and (\ref{strong}) yield
\begin{equation*}
J_\varphi\, y\, J_\varphi\, x_\varphi =x\, J_\varphi\, y_\varphi
= J_\varphi\, y\, J_\varphi\,\xi\, ,
\end{equation*}
so by the $so$-density of $\mathfrak M_\varphi$ in $M\ni 1_H$
we conclude that $x_\varphi =\xi\,$.

\hfill $\square\quad$

\medskip

{\it Proof of }$\,${\bf Proposition \ref{Leb}}.
\medskip

(1) Let $\mathcal D$ be the linear span of $\{ J_\varphi\, a^*
J_\varphi\, b_\varphi\, ;\, a\, ,\, b\in\mathfrak N_\varphi\}\,$
, which is dense in $H\,$. Define the linear functional
$F : \mathcal D \longrightarrow \mathbb C$ by
\begin{equation*}
F(\eta )=\lim_\iota\, (\eta\, |\, (x_\iota )_\varphi )\, ,\qquad
\eta\in\mathcal D\, ,
\end{equation*}
where the limit exists due to the convergence
\begin{equation}\label{F}
\begin{split}
(J_\varphi\, a^* J_\varphi\, b_\varphi\, |\, (x_\iota )_\varphi )
&=(b_\varphi\, |\, J_\varphi\, a\, J_\varphi\, (x_\iota )_\varphi )
\overset{(\ref{j})}{=} (b_\varphi\, |\, x_\iota\, J_\varphi\,
a_\varphi ) \\
&\overset{\iota}{\longrightarrow}
(b_\varphi\, |\, x\, J_\varphi\, a_\varphi )\, .
\end{split}
\end{equation}
Since $F$ is bounded by
\begin{equation*}
\| F\|\leq \sup_\iota\, \| (x_\iota )_\varphi\|<\infty\, ,
\end{equation*}
it extends to a continuous linear functional on $H\,$, and hence
there exists $\xi\in H$ satisfying $F(\eta ) =(\eta\, |\,\xi )$
for all $\eta\in\mathcal D\,$. In particular, by (\ref{F}),
\begin{equation*}
(b_\varphi\, |\, x\, J_\varphi\, a_\varphi ) =
F( J_\varphi\, a^* J_\varphi\, b_\varphi ) =
(J_\varphi\, a^* J_\varphi\, b_\varphi\, |\,\xi ) =
(b_\varphi\, |\, J_\varphi\, a\, J_\varphi\,\xi )\, ,\qquad
a\, ,\, b\in\mathfrak N_\varphi\, .
\end{equation*}
This implies that $x\, J_\varphi\, a_\varphi =J_\varphi\, a\,
J_\varphi\,\xi$ for all $a\in\mathfrak N_\varphi$ and by
Lemma \ref{char.N}.(2) we get
\begin{equation*}
x\in\mathfrak N_\varphi\,\text{ and }\, x_\varphi =\xi\, .
\end{equation*}
Furthermore,
\begin{equation*}
\lim_\iota\, (\eta\, |\, (x_\iota )_\varphi ) =F(\eta ) =(\eta\,
|\,\xi )=(\eta\, |\, x_\varphi )\, ,\qquad\eta\in\mathcal D\, ,
\end{equation*}
the density of $\mathcal D$ in $H$ and the boundedness of the net
$\{ (x_\iota )_\varphi\}_\iota$ yield that
\begin{equation*}
(x_\iota )_\varphi \overset{\iota}{\longrightarrow} x_\varphi\,
\text{ in the weak topology of }\, H\, .
\end{equation*}

(2) Let $x\in M$ be any $wo$-limit point of the bounded net
$\{ x_\iota\}_\iota\,$. Then, for every $a\in\mathfrak N_\varphi\,$,
$x\, J_\varphi\, a_\varphi$ is a weak limit point of the net
$\{ x_\iota\, J_\varphi\, a_\varphi\}_\iota \overset{(\ref{j})}{=}
\{ J_\varphi\, a\, J_\varphi\, (x_\iota)_\varphi\}_\iota\,$. Since
$(x_\iota)_\varphi\overset{\iota}{\longrightarrow}\xi\,$, we deduce
\begin{equation*}
x\, J_\varphi\, a_\varphi =J_\varphi\, a\, J_\varphi\,\xi\, ,
\qquad a\in\mathfrak N_\varphi
\end{equation*}
and, using Lemma \ref{char.N}.(2), we obtain
\begin{equation*}
x\in\mathfrak N_\varphi\,\text{ and }\, x_\varphi =\xi\, .
\end{equation*}
By the injectivity of the mapping $\,\mathfrak N_\varphi\ni y
\longmapsto y_\varphi\,$, the uniqueness of the $wo$-limit point
$x$ of $\{ x_\iota\}_\iota$ follows and we conclude that
$wo-\lim\limits_\iota\, x_\iota =x\,$.

\hfill $\square\quad$

\medskip

{\it Proof of }$\,${\bf Proposition \ref{Tom}}.
\medskip

First we show that every $y\in\mathfrak A_\varphi$ can be approximated
by a sequence $\{ y_n\}_{n\geq 1}$ in $\mathfrak T_\varphi$ as
required in the statement and such that (\ref{intertwin}) holds for
$x=y_n\, ,\, n\geq 1\,$.

Set
\begin{equation*}
y_n=\sqrt{\frac n\pi}\; wo-\int_{-\infty}^\infty e^{-n t^2}
\sigma^\varphi_t(y)\, dt\, ,\qquad n\geq 1\, .
\end{equation*}
Then every orbit
\begin{equation*}
\mathbb R\ni s\longmapsto \sigma^\varphi_s (y_n) =
\sqrt{\frac n\pi}\; wo-\int_{-\infty}^\infty e^{-n (t-s)^2}
\sigma^\varphi_t (y)\, dt\in M
\end{equation*}
has the entire extension
\begin{equation}\label{compl.y_n}
\mathbb C\ni z\longmapsto \sigma^\varphi_z (y_n) =
\sqrt{\frac n\pi}\; wo-\int_{-\infty}^\infty e^{-n (t-z)^2}
\sigma^\varphi_t (y)\, dt\in M
\end{equation}
and by (\ref{integr}) we have $\sigma^\varphi_z (y_n)\in
\mathfrak N_\varphi$ for every $z\in\mathbb C\,$. Similarly,
\begin{equation*}
\mathbb R\ni s\longmapsto \sigma^\varphi_s (y_n^{\, *}) =
\sqrt{\frac n\pi}\; wo-\int_{-\infty}^\infty e^{-n (t-s)^2}
\sigma^\varphi_t (y^*)\, dt\in M
\end{equation*}
has an entire extension $\,\mathbb C\ni z\longmapsto
\sigma^\varphi_z (y_n^{\, *})$ and $\sigma^\varphi_z (y_n^{\, *})
\in\mathfrak N_\varphi$ for every $z\in\mathbb C\,$. Since
$\sigma^\varphi_z (y_n)^*=\sigma^\varphi_{\overline{z}}
(y_n^{\, *})\,$, we have
\begin{equation*}
\sigma^\varphi_z (y_n)\in (\mathfrak N_\varphi )^* \cap\,
\mathfrak N_\varphi =\mathfrak A_\varphi\, ,\qquad n\geq 1\, ,\,
z\in\mathbb C\, ,
\end{equation*}
that is $y_n\in\mathfrak T_\varphi$ for all $n\geq 1\,$. By the
$so$-continuity of $\mathbb R\ni t\longmapsto \sigma^\varphi_t
(y)\in M$
and $\mathbb R\ni t\longmapsto \sigma^\varphi_t (y^*)\in M$ we get
$y_n\overset{so}{\longrightarrow}y$ and $y_n^{\, *}
\overset{so}{\longrightarrow}y_{}^{\, *}\,$, while using
\begin{equation}\label{exp}
\big|\, e^{-n (t-z)^2}\big| =\big|\, e^{-n (t-\Re z )^2}
e^{n (\Im z )^2} e^{2\, n\, i (t-\Re z )\,\Im z}\big| =
e^{-n (t-\Re z )^2} e^{n (\Im z )^2}
\end{equation}
it is easily seen that $\| \sigma^\varphi_z (y_n )\|\leq
e^{n (\Im z)^2}\| y\|$ for all $n\geq 1$ and $z\in\mathbb C\,$.

On the other hand, by (\ref{integr}) we have
\begin{equation*}
(y_n)_\varphi =\sqrt{\frac n\pi}\,\int_{-\infty}^\infty
e^{-n t^2}\Delta_\varphi^{i t}\, y_\varphi\, dt\, ,\quad
(y_n^{\, *})_\varphi =\sqrt{\frac n\pi}\,\int_{-\infty}^\infty
e^{-n t^2}\Delta_\varphi^{i t}\, (y^*)_\varphi\, dt
\end{equation*}
and by the norm-continuity of $\mathbb R\ni t\longmapsto
\Delta_\varphi^{i t}\, y_\varphi\in H$ and $\mathbb R\ni t
\longmapsto\Delta_\varphi^{i t}\, (y^*)_\varphi\in H$ we get
the convergences $(y_n)_\varphi\longrightarrow y_\varphi$ and
$(y_n^{\, *})_\varphi\longrightarrow (y^*)_\varphi$ in the
norm-topology. Furthermore, every orbit
\begin{equation*}
\mathbb R\ni s\longmapsto \Delta_\varphi^{i s}\, (y_n)_\varphi
=\sqrt{\frac n\pi}\,\int_{-\infty}^\infty e^{-n (t-s)^2}
\Delta_\varphi^{i t}\, y_\varphi\, dt\in H
\end{equation*}
has the entire extension
\begin{equation*}
\mathbb C\ni z\longmapsto
\sqrt{\frac n\pi}\, \int_{-\infty}^\infty e^{-n (t-z)^2}
\Delta_\varphi^{i t}\, y_\varphi\, dt\in H
\end{equation*}
and thus (see \cite{P-T}, Lemma 3.2 and \cite{C-Z}, Theorem 6.1)
$\displaystyle (y_n)_\varphi\in\bigcap_{z\in\mathbb C}
\text{Dom}\,\Delta_\varphi^{i z}$ and
\begin{equation}\label{compl.y_n_phi}
\Delta_\varphi^{i z}\, (y_n)_\varphi =\sqrt{\frac n\pi}\,
\int_{-\infty}^\infty e^{-n (t-z)^2} \Delta_\varphi^{i t}\,
y_\varphi\, dt\, ,\qquad z\in\mathbb C\, .
\end{equation}
Similarly, $\displaystyle (y_n^{\, *})_\varphi\in
\bigcap_{z\in\mathbb C}\,\text{Dom} \Delta_\varphi^{i z}$ and
\begin{equation*}
\Delta_\varphi^{i z}\, (y_n^{\, *})_\varphi =\sqrt{\frac n\pi}\,
\int_{-\infty}^\infty e^{-n (t-z)^2} \Delta_\varphi^{i t}\,
(y^*)_\varphi\, dt\, ,\qquad z\in\mathbb C\, .
\end{equation*}
Moreover, using (\ref{exp}), we get for every $n\geq 1$ and $z\in
\mathbb C$
\begin{equation*}
\| \Delta_\varphi^{i z}\, (y_n)_\varphi\|\leq e^{n (\Im z)^2}
\| y_\varphi\|\, ,\quad \| \Delta_\varphi^{i z}\,
(y_n^{\, *})_\varphi\|\leq e^{n (\Im z)^2}\| (y^*)_\varphi\|
\, .
\end{equation*}

Finally, for every $n\geq 1\,$, (\ref{compl.y_n}), (\ref{integr})
and (\ref{compl.y_n_phi}) yield
\begin{equation}\label{intertwiny_n}
\sigma^\varphi_z (y_n)_\varphi =\sqrt{\frac n\pi}\,
\int_{-\infty}^\infty e^{-n (t-z)^2} \sigma^\varphi_t (y)_\varphi
\, dt =\Delta_\varphi^{i z}\, (y_n)_\varphi\, ,\qquad z\in
\mathbb C\, ,
\end {equation}
hence (\ref{intertwin}) holds for $x=y_n\,$.

It remains to prove that (\ref{intertwin}) holds in full generality.
First we show that
\begin{equation}\label{part.intertwin}
x\in\mathfrak T_\varphi\, ,\,    
z\in\mathbb C\; ,\, x_\varphi\in\text{Dom}\, (\Delta_\varphi^{i z})\;
\Longrightarrow\,\sigma^\varphi_z (x)_\varphi =\Delta_\varphi^{i z}\,
x_\varphi\, .
\end{equation}
By Lemma \ref{char.N} this is equivalent to the implication
\begin{equation*}
x\in\mathfrak T_\varphi\, ,\,    
z\in\mathbb C\; ,\, x_\varphi\in\text{Dom}\, (\Delta_\varphi^{i z})
\, ,\, y\in\mathfrak A_\varphi\;\Longrightarrow\,\sigma^\varphi_z (x)
\, J_\varphi\, y_\varphi =J_\varphi\, y\, J_\varphi\,
\Delta_\varphi^{i z}\, x_\varphi\, ,
\end{equation*}
what we now are going to prove. Choose a sequence
$\{ y_n\}_{n\geq 1}$ in $\mathfrak T_\varphi$ as in the above part
of the proof. For each $n\geq 1\,$, $(y_n)_\varphi\in\text{Dom}\,
(\Delta_\varphi^{-i\,\overline{z}})$ implies by (\ref{mod.comm}) that
\begin{equation*}
J_\varphi\, (y_n)_\varphi\in\text{Dom}\, (\Delta_\varphi^{-i z})\,
\text{ and }\,\Delta_\varphi^{-i z} J_\varphi\, (y_n)_\varphi =
J_\varphi\,\Delta_\varphi^{-i\,\overline{z}}\, (y_n)_\varphi\, ,
\end{equation*}
so, according to Theorem \ref{anal.ext},
\begin{equation}\label{first}
\begin{split}
&x\,\Delta_\varphi^{-i z} J_\varphi\, (y_n)_\varphi =
x\, J_\varphi\,\Delta_\varphi^{-i\,\overline{z}}\, (y_n)_\varphi
\in\text{Dom}\, (\Delta_\varphi^{i z})\,\text{ and} \\
&\sigma^\varphi_z (x)\, J_\varphi\, (y_n)_\varphi =
\Delta_\varphi^{i z}\, x\,\Delta_\varphi^{-i z} J_\varphi\,
(y_n)_\varphi =\Delta_\varphi^{i z}\, x\, J_\varphi\,
\Delta_\varphi^{-i\,\overline{z}}\, (y_n)_\varphi\, .
\end{split}
\end{equation}
Since, by (\ref{intertwiny_n}) and by (\ref{j}),
\begin{equation*}
x\, J_\varphi\,\Delta_\varphi^{-i\,\overline{z}}\, (y_n)_\varphi =
x\, J_\varphi\,\sigma^\varphi_{-\overline{z}} (y_n)_\varphi =
J_\varphi\, \sigma^\varphi_{-\overline{z}} (y_n)\, J_\varphi\,
x_\varphi\, ,
\end{equation*}
(\ref{first}) yields
\begin{equation*}
J_\varphi\, \sigma^\varphi_{-\overline{z}} (y_n)\, J_\varphi\,
x_\varphi\in\text{Dom}\, (\Delta_\varphi^{i z})\,\text{ and }\,
\sigma^\varphi_z (x)\, J_\varphi\, (y_n)_\varphi =
\Delta_\varphi^{i z}\, J_\varphi\,\sigma^\varphi_{-\overline{z}}
(y_n)\, J_\varphi\, x_\varphi\, .
\end{equation*}
Using again (\ref{mod.comm}), we obtain
\begin{equation}\label{second}
\sigma^\varphi_{-\overline{z}} (y_n)\, J_\varphi\, x_\varphi
\in\text{Dom}\, (\Delta_\varphi^{i\,\overline{z}})\,\text{ and }\,
\sigma^\varphi_z (x)\, J_\varphi\, (y_n)_\varphi =J_\varphi\,
\Delta_\varphi^{i\,\overline{z}}\,\sigma^\varphi_{-\overline{z}}
(y_n)\, J_\varphi\, x_\varphi\, .
\end{equation}
Taking into account that $x_\varphi\in\text{Dom}\,
(\Delta_\varphi^{i z})$ and, by Theorem \ref{anal.ext} and by
(\ref{mod.comm}),
\begin{equation*}
J_\varphi\,\Delta_\varphi^{i\,\overline{z}}\,
\sigma^\varphi_{-\overline{z}} (y_n)\, J_\varphi\supset J_\varphi\,
y_n\,\Delta_\varphi^{i\,\overline{z}}\, J_\varphi\supset J_\varphi\,
y_n\, J_\varphi\,\Delta_\varphi^{i z}\, ,
\end{equation*}
(\ref{second}) implies the equality $\sigma^\varphi_z (x)\,
J_\varphi\, (y_n)_\varphi =J_\varphi\, y_n\, J_\varphi\,
\Delta_\varphi^{i z}\, x_\varphi\,$. Passing now to the limit for
$n\to\infty\,$, we conclude that
\begin{equation*}
\sigma^\varphi_z (x)\, J_\varphi\, y_\varphi =J_\varphi\, y\,
J_\varphi\,\Delta_\varphi^{i z}\, x_\varphi\, .
\end{equation*}

Next we show by induction on $k$ that
\begin{equation}\label{k}
x\in\mathfrak T_\varphi\;\Longrightarrow\, x_\varphi\in\text{Dom}\,
\big(\Delta_\varphi^{\frac k2}\big)
\end{equation}
holds for every integer $k\geq 1\,$. Indeed,
\begin{equation*}
x\in\mathfrak T_\varphi\subset\mathfrak A_\varphi\;\Longrightarrow\,
x_\varphi\in\text{Dom}\, (S_\varphi ) =\text{Dom}\,\big(
\Delta_\varphi^{\frac 12}\big)
\end{equation*}
is clear and if (\ref{k}) holds for some $k\geq 1$ and $x\in
\mathfrak T_\varphi\,$, then we have by (\ref{part.intertwin})
\begin{equation*}
\Delta_\varphi^{\frac k2}\, x_\varphi =\sigma^\varphi_{-\frac k2 i}
(x)_\varphi\in\text{Dom}\, (S_\varphi ) =\text{Dom}\,\big(
\Delta_\varphi^{\frac 12}\big)\, ,\text{ that is } x_\varphi\in
\text{Dom}\,\big(\Delta_\varphi^{\frac{k+1}2}\big)\, .
\end{equation*}
On the other hand, for every $x\in\mathfrak T_\varphi$ and $k\geq 1$
we have $\sigma^\varphi_{\frac k2 i}(x)\in\mathfrak T_\varphi\,$, so
(\ref{k}) yields $\sigma^\varphi_{\frac k2 i}(x)_\varphi\in
\text{Dom}\,\big(\Delta_\varphi^{\frac k2}\big)$ and, using
(\ref{part.intertwin}), we deduce
\begin{equation*}
x_\varphi =\sigma^\varphi_{-\frac k2 i}\big(
\sigma^\varphi_{\frac k2 i}(x)\big)_\varphi =\Delta_\varphi^{\frac k2}
\,\sigma^\varphi_{\frac k2 i}(x)_\varphi\in\text{Dom}\,\big(
\Delta_\varphi^{-\frac k2}\big)\, .
\end{equation*}
Therefore (\ref{k}) holds also for every integer $k\leq -1\,$ that is
$\displaystyle \mathfrak T_\varphi\subset\bigcap_{z\in\mathbb C}
\text{Dom}\,\Delta_\varphi^{\, z}\,$.
This last inclusion together with (\ref{part.intertwin}) imply
(\ref{intertwin}).

\hfill $\square\quad$

\bigskip
\section{Hermitian Maps}
\bigskip

In this section we analyse the notion introduced in Definition
\ref{def.herm} by proving Lemma \ref{fixed}, Lemma \ref{herm} and
Theorem \ref{char.herm}.
\medskip

{\it Proof of }$\,${\bf Lemma \ref{fixed}}.
\medskip

Let $\xi\in\,$Dom$\, S_\psi$ be arbitrary. Then there is a 
sequence $(x_n)_{n\geq 1}$ in $\mathfrak A_\psi$ such that
\begin{equation*}
(x_n^{})_\psi\longrightarrow\xi\, ,\qquad (x_n^*)_\psi\longrightarrow
S_\psi\,\xi\, .
\end{equation*}
Then, denoting
\begin{equation*}
\xi_+ =\frac 12\, (\xi +S_\psi\,\xi )\, ,\qquad \xi_- =\frac 1{2\, i}\,
(\xi -S_\psi\,\xi )\, ,\qquad a_n =\frac 12\, (x_n^{} +x_n^*)\, ,
\end{equation*}
we have
\begin{equation}\label{split}
\xi =\xi_+ +i\,\xi_-\, ,\qquad S_\psi\,\xi_\pm =\xi_\pm\, ,\,
\text{ i.e. }\,\xi_\pm\in K^{S_\psi}\, ,
\end{equation}
\begin{equation}\label{appr}
a_n^*=a_n^{}\in\mathfrak N_\psi\, ,\qquad (a_n^{})_\psi\longrightarrow
\xi_+\, .
\end{equation}
\vspace{0.1 pt}

Since $\xi_+=\xi$ if $\xi\in K^{S_\psi}\,$, (\ref{appr}) implies that
$K^{S_\psi}\subset\overline{\{ x_\psi\, ;\, x^*=x\in
\mathfrak N_\psi\}}\,$. The converse inclusion being trivial, the
equality $K^{S_\psi} =\overline{\{ x_\psi\, ;\, x^*=x\in
\mathfrak N_\psi\}}$ follows. On the other hand, (\ref{split})
implies that Dom$\, S_\psi =K^{S_\psi} +i\, K^{S_\psi}\,$. This
proves (1) and (4) in Lemma \ref{fixed}.

For (2) and (3) we first notice that, for every $\xi\in K^{S_\psi}$
and $x\in\mathfrak A_\psi\,$,
\begin{equation*}
(\xi\, |\, J_\psi\, x_\psi )=(\xi\, |\, J_\psi\, S_\psi\, (x^*)_\psi )
=(\xi\, |\,\Delta_\psi^{1/2} (x^*)_\psi )=(\Delta_\psi^{1/2} \xi\, |\,
(x^*)_\psi )=(J_\psi\, (x^*)_\psi\, |\,\xi )\, .
\end{equation*}
In particular, $(\xi\, |\, J_\psi\, x_\psi )\in\mathbb R$ whenever
$x=x^*\,$.

Conversely, let us assume that $\xi\in K$ is such that $(\xi\, |\, J_\psi\,
x_\psi )\in\mathbb R$ if $x^*=x\in\mathfrak N_\psi\,$. For every
$x\in\mathfrak A_\psi$ we have
\begin{equation*}
a=\frac 12\, (x+x^*)\, ,\, b=\frac 1{2\, i}\, (x-x^*)\in
\mathfrak A_\psi\,\text{ are selfadjoint and }\; x=a+i\, b\, ,
\end{equation*}
hence, by our assumption on $\xi\,$,
\begin{align*}
(\xi\, |\, J_\psi\,x_\psi ) &=\frac 12\,\big( (\xi\, |\, J_\psi\,
a_\psi ) +i\, (\xi\, |\, J_\psi\,b_\psi )\big) =\frac 12\,\big(
(J_\psi\, a_\psi\, |\,\xi ) +i\, (J_\psi\, b_\psi\, |\,\xi )
\big) = \\
&=\big( J_\psi\, \frac 12\, (a-i\, b)_\psi\, |\,\xi\big) =
(J_\psi\, (x^*)_\psi\, |\,\xi )\, .
\end{align*}
It follows that $(x_\psi\, |\, J_\psi\,\xi ) =(\,\Delta_\psi^{1/2}\,
x_\psi\, |\,\xi )$ for all $x\in\mathfrak A_\psi\,$, hence,
$\{ x_\psi\, ;\, x\in\mathfrak A_\psi\}$ being a core of
$\Delta_\psi^{1/2}\,$, $\xi$ belongs to the domain of
$(\Delta_\psi^{1/2})^*=\Delta_\psi^{1/2}$ and $\Delta_\psi^{1/2}\,\xi
=J_\psi\,\xi\,$. In other words, $\xi\in\,$Dom$\, S_\psi$ and
$S_\psi\,\xi =J_\psi\,\Delta_\psi^{1/2}\,\xi =\xi\,$, i.e.
$\xi\in K^{S_\psi}\,$.

\hfill $\square\quad$
\medskip

{\it Proof of }$\,${\bf Lemma \ref{herm}}.
\medskip

If $T\in B(K,H)$ implements $\psi$ in $\varphi$ and $x^*=x\in
\mathfrak N_\psi\,$, then
\begin{equation*}
T\, x_\psi =(T\, x\, T^*)_\varphi\,\text{ with }\, (T\, x\, T^*)^* =
T\, x\, T^*\in\mathfrak N_\varphi
\end{equation*}
and hence Lemma \ref{fixed}.(1) implies $T\, K^{S_\psi}\subset
H^{S_\varphi}\,$, proving (1). In this case the inclusion
$T\,\mathfrak N_\psi\, T^*\subset\mathfrak N_\varphi$ and the
$wo$-density of $\mathfrak N_\psi$ in $N$ imply $T\, N\, T^*
\subset M\,$, proving (2).

If $T$ is isometric in addition, then
\begin{equation*}
x=T^*(T\, x\, T^*)\, T\, ,\qquad x\in N
\end{equation*}
shows the injectivity of the map $N\ni x\longmapsto T\, x\, T^*\in
M\,$, which is clearly also a $*$-homomorphism. Furthermore, $0\leq
a\in\mathfrak M_\psi$ implies $a^{1/2}\in\mathfrak N_\psi$ and
\begin{equation*}
\psi (a) =\| (a^{1/2})_\psi\|^2 =\| T\, (a^{1/2})_\psi\|^2 =
\| (T\, a^{1/2}\, T^*)_\varphi\|^2 =\varphi (T\, a\, T^*)\, .
\end{equation*}
Therefore (3) holds.

Let us finally assume that $\psi$ and $\varphi$ are bounded and
$\xi_\psi =(1_K)_\psi\, ,\,\eta_\varphi =(1_H)_\varphi\,$. If
$T\in B(K,H)$ is injective and implements $\psi$ in $\varphi\,$,
then $T\, N\, T^*\subset M$ by the above proved (2) and
\begin{equation*}
T\, T^*\eta_\varphi =(T\, 1_K\, T^*)_\varphi =T\, (1_K)_\psi =
T\,\xi_\psi \;\;\Longrightarrow\;\; T^*\eta_\varphi =\xi_\psi
\end{equation*}
by the injectivity of $T\,$. Conversely, if $T\in B(K,H)$ is
injective and satisfies $T\, N\, T^*\subset M$ and $T^*\eta_\varphi
=\xi_\psi\,$, then for $x\in N$
\begin{equation*}
(T\, x\, T^*)_\varphi =T\, x\, T^*\eta_\varphi =T\, x\,\xi_\psi
=T\, x_\psi\, .
\end{equation*}
Hence we have (4).

\hfill $\square\quad$
\medskip

{\it Proof of }$\,${\bf Theorem \ref{char.herm}}.
\medskip

(1), (2) and (3) in Lemma \ref{fixed} imply the equivalences
(1)$\,\Leftrightarrow\,$(2), (2)$\,\Leftrightarrow\,$(3) and
(2)$\,\Leftrightarrow\,$(4), respectively.

Let us assume that (1) holds.
By (4) in Lemma \ref{fixed}, every $\xi\in\,$Dom$\, S_\psi$ is of the
form $\xi =\xi_1 +i\,\xi_2$ with $\xi_1\, ,\,\xi_2\in K^{S_\psi}\,$.
Hence we get
\begin{align*}
&T(\xi )=T(\xi_1 )+i\, T(\xi_2 )\in H^{S_\varphi} +i\, H^{S_\varphi}
\subset\,\text{Dom}\, S_\varphi\; , \\
&S_\varphi\big( T(\xi )\big) =T(\xi_1 )-i\, T(\xi_2 )=T(\xi_1 -i\,
\xi_2 ) =T\big( S_\psi (\xi )\big)\; ,
\end{align*}
proving (5). Conversely, if (5) holds, then we have for every
$\xi\in K^{S_\psi}\subset\,\text{Dom}\, S_\psi$
\begin{equation*}
T\, (\xi )\in\,\text{Dom}\, S_\varphi\,\text{ and }\,S_\varphi
\big( T(\xi )\big) =T\big( S_\psi (\xi )\big) =T(\xi )\, ,
\end{equation*}
so $T(\xi )\in H^{S_\varphi}\,$. Therefore (1)$\,\Leftrightarrow
\,$(5).

Since $J_\psi$ is involutive and, by (\ref{mod.comm-}), $S_\psi =
\Delta_\psi^{-1/2}\, J_\psi$ and $S_\varphi =\Delta_\varphi^{-1/2}\,
J_\varphi\,$, (5) is equivalent to
\begin{equation*}
T\,\Delta_\psi^{-1/2}\subset \Delta_\varphi^{-1/2}\, J_\varphi\, T\,
J_\psi\, .
\end{equation*}
This equation, in turn, is equivalent to the validity of
\begin{equation*}
\Delta_\varphi^{1/2}\, T\,\Delta_\psi^{-1/2}\xi =J_\varphi\, T\,
J_\psi\,\xi\, ,\qquad \xi\in\,\text{Dom}\,\Delta_\psi^{-1/2}
\end{equation*}
and thus (5)$\,\Leftrightarrow\,$(6).

We have already seen that (1)$\,\Leftrightarrow\,$(3). Applying this
equivalence to $J_\psi\, T^* J_\varphi\,$, it follows that
$J_\psi\, T^* J_\varphi\,$ is Hermitian with respect to
$(\varphi ,\psi )$ if and only if
\begin{equation*}
(J_\psi\, T^* J_\varphi\, y_\varphi\, |\, J_\psi x_\psi )=
(x_\psi\, |\, T^* J_\varphi\, y_\varphi )=(T\, x_\psi\, |\,
J_\varphi\, y_\varphi )
\end{equation*}
is real for all $x^*=x\in\mathfrak N_\psi$ and $y^*=y\in
\mathfrak N_\varphi\,$. But this means exactly (3), so (3)$\,
\Leftrightarrow\,$(7).

By the equivalence of statements (1) and (5) in Theorem
\ref{anal.ext} with $A=\Delta_\varphi\,$, $B=\Delta_\psi\,$,
$\beta =-\frac 12\,$, and taking into account that they imply
(\ref{dom}) and (\ref{op}), we obtain the equivalence (6)$\,
\Leftrightarrow\,$(8).

Now let us assume that the equivalent conditions $(1) - (8)$ are
satisfied. Then, by (\ref{transl}) in Theorem \ref{anal.ext}, we
get (\ref{mod.transl}). Further, using (\ref{mod.comm-}), we
obtain (\ref{transl-i/2}) immediately from (\ref{mod.transl})
and (\ref{-i/2}). Since the map $\, T(\,\cdot\, )$ is bounded and
\begin{equation*}
\| T(s)\| =\|\Delta_\varphi^{is}\, T\,\Delta_\psi^{-is}\| =\| T\|\, ,
\; \Big\| T\Big( s-\frac i2 \Big)\Big\|
\overset{(\ref{transl-i/2})}{=}
\| J_\varphi\, T(s)\, J_\psi\| =\| T\|\, ,\qquad s\in\mathbb R\, ,
\end{equation*}
we get also (\ref{norm}) by the Three Line Theorem. Finally,
since $K^{S_\psi}$ and $H^{S_\varphi}$ are invariant under
$\Delta_\psi^{-is}$ and $\Delta_\varphi^{is}\,$, respectively,
for every $s\in\mathbb R\,$, due to (\ref{delta}) and Lemma
\ref{fixed}.(1), we obtain the Hermiticity of $T(s)=
\Delta_\varphi^{is}\, T\,\Delta_\psi^{-is}$ from (1).

\hfill $\square\quad$

\bigskip
\section{Generalization of the Structure Theorem of Borchers}
\bigskip

We prove Theorem \ref{borchers} in two steps: first we prove it
for the case where $\Xi_o$ and $\Xi_1$ are empty, and then we reduce
the proof of the general case to the above special case.

\bigskip

{\it Step 1}. {\sc Proof in the case of $\Xi_o =\Xi_1 =\emptyset$
and $wo$-continuous $\, T(\,\cdot\, )\,$.}
\medskip

By our assumptions in this step,
\begin{equation*}
\overline{\mathbb S_\beta}\ni z\longmapsto T(z)\in B(K,H)
\end{equation*}
is a bounded $wo$-continuous map which is analytic in $\mathbb S_\beta$
and satisfies the boundary conditions
\begin{itemize}
\item[(i)] $T(s)$ is Hermitian with respect to $(\psi ,\varphi )$
for all $s\in\mathbb R\,$,
\item[(ii)] $J_\varphi\, T(s+i\beta )\, J_\psi$ is Hermitian with
respect to $(\psi ,\varphi )$ for all $s\in\mathbb R\,$.
\end{itemize}

Let $x\in\mathfrak T_\psi$ and $y\in\mathfrak T_\varphi$ be arbitrary
(for the Tomita algebras $\mathfrak T_\psi$ and $\mathfrak T_\varphi$
see the comments before Proposition \ref{Tom}) such that
\begin{equation}\label{growth}
\|\Delta_\varphi^{i z}\, x_\varphi\|\leq e^{c\, (\Im z)^2} \| x_\varphi\|
\, ,\;\|\Delta_\varphi^{i z}\, y_\varphi\|\leq e^{c\, (\Im z)^2}
\| y_\varphi\|\, ,\qquad z\in\mathbb C
\end{equation}
for some constant $c\geq 0\,$. Consider the functions
\begin{equation*}
f_1 : \mathbb C\times \overline{\mathbb S_\beta}\ni (z_1,z_2)
\longmapsto\big( T(z_2)\,\Delta_\psi^{-i z_1}\, x_\psi\, |\,
J_\varphi\,\Delta_\varphi^{-i z_1}\, y_\varphi\big)\, ,\phantom{xxx}
\end{equation*}
\begin{equation*}
f_2 : \mathbb C\times \overline{\mathbb S_{-\beta}}\ni (z_1,z_2)
\longmapsto\big( \Delta_\varphi^{-i z_1 +\frac 12}\, y_\varphi\, |\,
T(\overline{z_2} )\, J_\psi\,\Delta_\psi^{-i z_1 +\frac 12}\, x_\psi
\big)\, .
\end{equation*}
They are continuous and, according to (\ref{growth}), bounded on
any set of the form
\begin{equation*}
\{ z_1\in\mathbb C\, ;\, |\Im z_1 |\leq\delta \}\times
\overline{\mathbb S_\beta}\,\text{ and }\{ z_1\in\mathbb C\, ;\,
|\Im z_1 |\leq\delta \}\times\overline{\mathbb S_{-\beta}}\, ,
\text{ respectively },\quad \delta >0\, .
\end{equation*}
Furthermore, the partial functions
\begin{equation*}
\mathbb C\ni z_1\longmapsto f_1(z_1,z_2)\, ,\quad z_2\in
\overline{\mathbb S_\beta}\, ,\qquad\;\phantom{xx}
\mathbb S_\beta\ni z_2\longmapsto f_1(z_1,z_2)\, ,\quad z_1\in
\mathbb C\, ,
\end{equation*}
\begin{equation*}
\mathbb C\ni z_1\longmapsto f_2(z_1,z_2)\, ,\quad z_2\in
\overline{\mathbb S_{-\beta}}\, ,\qquad\;
\mathbb S_{-\beta}\ni z_2\longmapsto f_2(z_1,z_2)\, ,\quad z_1\in
\mathbb C\phantom{x}
\end{equation*}
are analytic.

Now, by (\ref{intertwin}) and Theorem \ref{char.herm}.(4), (i) implies,
for every $z_1\in\mathbb C$ and $s\in\mathbb R\,$,
\begin{align*}
f_1(z_1,s) &= \big( T(s)\,\sigma^\psi_{-z_1}(x)_\psi\, |\,
J_\varphi\,\sigma^\varphi_{-z_1}(y)_\varphi\big) =\\
&=\big( J_\varphi\,\big( \sigma^\varphi_{-z_1}(y)^*\big)_\varphi
\, |\, T(s)\,\big( \sigma^\psi_{-z_1}(x)^*\big)_\psi \big) =\\
&=\big( J_\varphi\, S_\varphi\,\sigma^\varphi_{-z_1}(y)_\varphi
\, |\, T(s)\, S_\psi\,\sigma^\psi_{-z_1}(x)_\psi \big) =\\
&=\big( \Delta_\varphi^{-i z_1 +\frac 12}\, y_\varphi\, |\,
T(s)\, J_\psi\,\Delta_\psi^{-i z_1 +\frac 12}\, x_\psi \big)
=f_2(z_1,s)\, .
\end{align*}
Therefore
\begin{equation*}
f : \mathbb C\times \big\{ z_2\in\mathbb C\, ;\, |\Im z_2 |\leq
|\beta |\big\} \ni (z_1,z_2)\,\longmapsto \begin{cases}
\; f_1(z_1,z_2) &\text{if $z_2\in\overline{\mathbb S_\beta}\,$,} \\
\; f_2(z_1,z_2) &\text{if $z_2\in\overline{\mathbb S_{-\beta}}$} \\
\end{cases}
\end{equation*}
is a well defined continuous function, bounded on every set of the
form
\begin{equation*}
\big\{ z_1\in\mathbb C\, ;\, |\Im z_1 |\leq\delta \big\}\times
\big\{ z_2\in\mathbb C\, ;\, |\Im z_2 |\leq |\beta |\big\}\, ,
\qquad \delta >0\, .
\end{equation*}
For each fixed $z_1\in\mathbb C\,$, the function $\mathbb S_\beta
\cup \mathbb S_{-\beta}\ni z_2\longmapsto f(z_1,z_2)$ is analytic.
Hence, by the Morera Theorem (the one-dimensional edge-of-the-wedge
theorem, see for example \cite{Ber-Gay}, 2.1.9.(2) or \cite{Car},
II.2.7), it can be analytically extended across $\mathbb R\,$,
that is the partial functions
\begin{equation*}
\big\{ z_2\in\mathbb C\, ;\, |\Im z_2 | <|\beta |\big\}\ni z_2\,
\longmapsto f(z_1,z_2)\, ,\qquad z_1\in\mathbb C
\end{equation*}
are analytic. Thus we can apply to $f$ the Osgood Lemma (the
Hartogs Theorem for continuous functions, see for example
\cite{G-R}, Theorem I.A.2) and deduce that it is analytic, as
function of two complex variables, on $\mathbb C\times
\big\{ z_2\in\mathbb C\, ;\, |\Im z_2 | <|\beta |\big\}\,$.

For every $z_1\in\mathbb C$ and $s\in\mathbb R\,$, (ii) implies
by (\ref{intertwin}) and Theorem \ref{char.herm}.(4),
\begin{align*}
f\Big( z_1+\frac i2\, ,\, s+i\,\beta \Big) &=f_1\Big( z_1+
\frac i2\, ,\, s+i\,\beta \Big) =\\
&=\big( T(s+i\beta )\, J_\psi\, S_\psi \Delta_\psi^{-i z_1}\,
x_\psi\, |\, S_\varphi\, \Delta_\varphi^{-i z_1}\, y_\varphi
\big) = \\
&=\big( T(s+i\beta )\, J_\psi\, \big(\sigma^\psi_{-z_1}
(x)^*\big)_\psi\, |\, \big(\sigma^\varphi_{-z_1}(y)^*
\big)_\varphi \big) =\\
&=\big( J_\varphi\, \big(\sigma^\varphi_{-z_1}(y)^*\big)_\varphi
\, |\, J_\varphi\, T(s+i\beta )\, J_\psi\, \big(
\sigma^\psi_{-z_1} (x)^*\big)_\psi \big) =\\
&=\big( J_\varphi\, T(s+i\beta )\, J_\psi\, \sigma^\psi_{-z_1}
(x)_\psi\, |\, J_\varphi\, \sigma^\varphi_{-z_1}(y)_\varphi
\big) =\\
&=\big( \sigma^\varphi_{-z_1}(y)_\varphi\, |\, T(s+i\beta )\,
J_\psi\, \sigma^\psi_{-z_1} (x)_\psi \big) =\\
&=\big( \Delta_\varphi^{-i (z_1-\frac i2 )+\frac 12}\, y_\varphi
\, |\, T(\overline{s-i \beta} )\, J_\psi\,
\Delta_\psi^{-i (z_1-\frac i2 )+\frac 12}\, x_\psi \big) =\\
&=f_2\Big( z_1-\frac i2\, ,\, s-i\,\beta \Big) =
f\Big( z_1-\frac i2\, ,\, s-i\,\beta \Big)\, .
\end{align*}

Therefore, for each $s\in\mathbb R\,$, the bounded, continuous
function
\begin{equation*}
g_s : \Big\{ \zeta\in\mathbb C\, ;\, |\Im \zeta |\leq\frac 12\Big\}
\ni\zeta\longmapsto f\big(\zeta\, ,\, s+2\,\beta\,\zeta\big)\, ,
\end{equation*}
which is analytic in the interior, satisfies
\begin{equation*}
g_s \Big( t+\frac i2\Big) =g_s \Big( t-\frac i2\Big)\, ,\qquad
t\in\mathbb R\, .
\end{equation*}
By the Morera Theorem, $g_s$ extends to a periodic entire function
with period $i\,$, still denoted by $g_s\,$, which is bounded. By
the Liouville Theorem it follows that $g_s$ is constant, hence we
get successively
\begin{equation*}
f_1(0\, ,s)=f(0\, ,s)=g_s(0)=g_s\Big(\,\frac {-s}{2\,\beta}\,\Big)
=f\Big(\,\frac {-s}{2\,\beta}\, ,\, 0\Big) =
f_1\Big(\,\frac {-s}{2\,\beta}\, ,\, 0\Big)\, ,
\end{equation*}
\begin{align*}
\big( T(s)\, x_\psi\, |\, J_\varphi\, y_\varphi\big) =\;\,& \big( T(0)
\,\Delta_\psi^{i \frac s{2\,\beta}}\, x_\psi\, |\, J_\varphi\,
\Delta_\varphi^{i \frac s{2\,\beta}}\, y_\varphi \big) =\\
\overset{(\ref{mod.comm-})}{=}&
\big( \Delta_\varphi^{-i \frac s{2\,\beta}}\, T(0)\,
\Delta_\psi^{i \frac s{2\,\beta}}\, x_\psi\, |\, J_\varphi\,
y_\varphi\big)\, .
\end{align*}

By the density property of $\mathfrak T_\varphi$ stated in
Proposition \ref{Tom}, the above equalities imply that
\begin{equation*}
T(s)=\Delta_\varphi^{-i \frac s{2\,\beta}}\, T(0)\,
\Delta_\psi^{i \frac s{2\,\beta}}\, ,\qquad s\in\mathbb R\, ,
\end{equation*}
hence (\ref{formula}) holds with $T=T(0)\,$. From (\ref{mod.transl})
and (\ref{transl-i/2}) in Theorem \ref{char.herm}, we obtain
(\ref{borch.transl}) and (\ref{transl+ibeta}).
\bigskip

{\it Step 2}. {\sc Proof in the general case.}
\medskip

Let us consider, for any integer $n\geq 1\,$, the entire function
\begin{equation*}
\mathbb C\ni z\longmapsto f_n(z)=\sqrt{\frac n{\pi}}\, e^{-n z^2}
\end{equation*}
and the mollification of $T(\,\cdot\, )$
\begin{equation}\label{moll}
\overline{\mathbb S_\beta}\ni z\longmapsto T_n(z)=
wo-\int_{-\infty}^\infty f_n(t)\, T(t+z)\, dt\in B(K,H)\, .
\end{equation}
We notice that the mapping $\mathbb R\ni t\longmapsto T(t+z)\in
B(K,H)$ is norm-continuous for $z\in\mathbb S_\beta$ and, due to
the continuity conditions (\ref{cont+}) and (\ref{cont-}), 
$wo$-measurable with respect to the Lebesgue measure for
$z\in\partial\,\mathbb S_\beta\,$. Since $T(\,\cdot\, )$ is bounded
and $f_n(t)\, dt$ is a probability measure, the integrals in
(\ref{moll}) exist and
\begin{equation}\label{bdd}
\sup\,\{\,\| T_n(z)\|\, ;\, z\in\overline{\mathbb S_\beta}\, ,\, n
\geq 1\}\leq \sup\,\{\,\| T(z)\|\, ;\, z\in\mathbb S_\beta\}
<\infty\, .
\end{equation}
Further, (\ref{cont+}) and (\ref{cont-}) yield by the Lebesgue
Dominated Convergence Theorem
\begin{equation}\label{moll.cont+}
T_n(s)=wo - \lim_{0<t/\beta \to 0} T_n(s+it)\, ,\qquad
s\in\mathbb R\, ,
\end{equation}
\begin{equation}\label{moll.cont-}
T_n(s+i\beta )=wo - \lim_{1>t/\beta \to 1} T_n(s+it)\, ,\qquad
s\in\mathbb R\, .
\end{equation}

We compare the operator valued function $T_n(\,\cdot\, )$ with
\begin{equation}\label{mod.moll}
\begin{split}
\mathbb C\ni z\longmapsto T_{\zeta ,n}(z)
&=wo-\int_{\mathbb R +\zeta} f_n(w-z)\, T(w)\, dw\\
&=wo-\int_{-\infty}^\infty f_n(t+\zeta -z)\, T(t+\zeta )\, dt
\in B(K,H)\, ,
\end{split}
\end{equation}
where $\zeta\in\mathbb S_\beta\,$. Due to
\begin{equation*}
f_n(t+\zeta -z)=\sqrt{\frac n{\pi}}\, e^{-n t^2}\,
e^{-2 n t (\zeta -z) -n (\zeta -z)^2}\, ,
\end{equation*}
the integral in (\ref{mod.moll}) is convergent and defines an entire
mapping $\, T_{\zeta ,n}(\,\cdot\, )\,$.

For $\zeta_1\, ,\,\zeta_2\in\mathbb S_\beta$ and $z\in\mathbb C$
we have by the Cauchy Integral Theorem
\begin{equation*}
T_{\zeta_1 ,n}(z) =
wo-\int_{\mathbb R +\zeta_1} f_n(w-z)\, T(w)\, dw =
wo-\int_{\mathbb R +\zeta_2} f_n(w-z)\, T(w)\, dw =
T_{\zeta_2 ,n}(z)\, ,
\end{equation*}
so $T_{\zeta ,n}(\,\cdot\, )$ does not depend on $\zeta\in
\mathbb S_\beta\,$. Therefore, for any $\zeta\in\mathbb S_\beta\,$,
\begin{equation}\label{openstr}
T_{\zeta ,n}(z)=T_{z ,n}(z)=T_n(z)\, ,\qquad z\in\mathbb S_\beta\, .
\end{equation}

Let $\zeta\in\mathbb S_\beta$ be arbitrary. Since
$\, T_{\zeta ,n}(\,\cdot\, )$ is an entire mapping, by
(\ref{moll.cont+}), (\ref{moll.cont-}) and (\ref{openstr}) we get
for every $s\in\mathbb R$
\begin{equation*}
T_n(s)=wo - \lim_{0<t/\beta \to 0} T_n(s+it)=
wo - \lim_{0<t/\beta \to 0} T_{\zeta ,n}(s+it)=T_{\zeta ,n}(s)\, ,
\end{equation*}
\begin{equation*}
T_n(s+i\beta )=wo - \lim_{1>t/\beta \to 1} T_n(s+it)=
wo - \lim_{1>t/\beta \to 1} T_{\zeta ,n}(s+it)=
T_{\zeta ,n}(s+i\beta )\, .
\end{equation*}

Consequently, the mapping $\,\overline{\mathbb S_\beta}\ni z
\longmapsto T_n(z)$ defined in (\ref{moll}) is a restriction of the
entire mapping $\, T_{\zeta ,n}(\,\cdot\, )\,$. In particular, it
is $so$-continuous (the role of $\, T_{\zeta ,n}(\,\cdot\, )$ is
just to prove this statement) and its restriction to $\mathbb S_\beta$
is analytic. We recall that its boundedness was already noticed in
(\ref{bdd}). 

Since $\mathbb R\ni t\longmapsto f_n(t)$ is a real function,
(1)$\,\Leftrightarrow\,$(3) in Theorem \ref{char.herm} implies that
the Hermiticity of $T(s)$ and $J_\varphi\, T(s+i\beta )\, J_\psi$
for $s\in\mathbb R\,\backslash\,\Xi_o$ respectively $s\in\mathbb R
\,\backslash\,\Xi_1$ is inherited by $T_n(s)$ and $J_\varphi\,
T_n(s+i\beta )\, J_\psi$ for all $s\in\mathbb R\,$.

Thus $T_n(\,\cdot\, )$ fulfils all the assumptions made in Step 1.
Consequently there exists $T_n\in B(K,H)$ satisfying
\begin{equation*}
T_n(s)=\Delta_\varphi^{-i \frac s{2\,\beta}}\, T_n\,
\Delta_\psi^{i \frac s{2\,\beta}}\, ,\qquad s\in\mathbb R\, .
\end{equation*}
It follows that
\begin{equation*}
T_n(s+t)=\Delta_\varphi^{-i \frac s{2\,\beta}}\, T_n(t)\,
\Delta_\psi^{i \frac s{2\,\beta}}\, ,\qquad t\, ,\; s\in\mathbb R\, ,
\end{equation*}
which yields by analytic extension
\begin{equation}\label{eq.n}
T_n(z+s)=\Delta_\varphi^{-i \frac s{2\,\beta}}\, T_n(z)\,
\Delta_\psi^{i \frac s{2\,\beta}}\, ,\qquad z\in
\overline{\mathbb S_\beta}\, ,\; s\in\mathbb R\, .
\end{equation}
On the other hand, from (\ref{moll}) we obtain
\begin{equation*}
\text{norm}-\lim_{n\to\infty} T_n(z)=T(z)\, ,\qquad z\in
\mathbb S_\beta\, ,
\end{equation*}
due to the boundedness and norm-continuity of $T(\,\cdot\, )$ in
$\mathbb S_\beta\,$. Thus we get by (\ref{eq.n})
\begin{equation*}
T(z+s)=\Delta_\varphi^{-i \frac s{2\,\beta}}\, T(z)\,
\Delta_\psi^{i \frac s{2\,\beta}}\, ,\qquad z\in\mathbb S_\beta\,
,\; s\in\mathbb R\, .
\end{equation*}

Now choose some $s_o\in\mathbb R\,\backslash\,\Xi_o$ and denote
\begin{equation*}
T=\Delta_\psi^{i \frac {s_o}{2\,\beta}}\, T(s_o)\,
\Delta_\varphi^{-i \frac {s_o}{2\,\beta}}\, .
\end{equation*}
Then $T$ is Hermitian with respect to $(\psi ,\varphi )$ and we have
for every $s\in\mathbb R\,\backslash\,\Xi_o$
\begin{align*}
T(s) &=wo - \lim_{0<t/\beta \to 0} T(s+it) \\
&=wo - \lim_{0<t/\beta \to 0}
\Delta_\varphi^{-i \frac {s-s_o}{2\,\beta}}\, T(s_o+it)\,
\Delta_\varphi^{i \frac {s-s_o}{2\,\beta}} \\
&=\Delta_\varphi^{-i \frac {s-s_o}{2\,\beta}} \Big(
wo - \lim_{0<t/\beta \to 0} T(s_o+it) \Big)\,
\Delta_\varphi^{i \frac {s-s_o}{2\,\beta}} \\
&=\Delta_\varphi^{-i \frac s{2\,\beta}} \Big(
\Delta_\psi^{i \frac {s_o}{2\,\beta}}\, T(s_o)\,
\Delta_\varphi^{-i \frac {s_o}{2\,\beta}} \Big)
\Delta_\psi^{i \frac s{2\,\beta}} \\
&=\Delta_\varphi^{-i \frac s{2\,\beta}}\, T\,
\Delta_\psi^{i \frac s{2\,\beta}}\, ,
\end{align*}
that is (\ref{formula}) holds. Using (1)$\,\Rightarrow\,$(8) in
Theorem \ref{char.herm}, (\ref{mod.transl}) and (\ref{transl-i/2}),
as well as the uniqueness theorem of the Riesz brothers (see, in
Section 3, Step 1 of the proof of (3)$\,\Rightarrow\,$(4) in Theorem
\ref{anal.ext}), we obtain that $T(\,\cdot\, )$ extends to an
$so$-continuous map on $\mathbb S_\beta\,$, for which
(\ref{borch.transl}) and (\ref{transl+ibeta}) hold.

\hfill $\square\quad$

\bigskip
\section{Proof of Theorem \ref{hsmi}}
\bigskip

We recall the setting:
\begin{itemize}
\item[--] $N\subset M\subset B(H)$ are von Neumann algebras,
\item[--] $\varphi$ is a normal semifinite faithful weight on $M$
such that its restriction $\psi$ to $N$ is semifinite,
\item[--] we assume that $M$ is in standard form with respect to
$\varphi$ and $\{ y_\varphi\, ;\, y\in\mathfrak N_\psi\}$ is dense
in $H\,$, hence $N$ is in standard form with respect to $\psi$ and
$y_\psi =y_\varphi$ for all $y\in\mathfrak N_\psi\,$.
\end{itemize}
We shall use the notations $\Delta_M =\Delta_\varphi\, ,\, J_M=
J_\varphi$ and $\,\Delta_N =\Delta_\psi\, ,\, J_N=J_\psi\,$.

The proof of Theorem \ref{hsmi} will be performed in nine steps.
\bigskip

{\it Step 1}. {\sc Application of the Modular Extension Theorem.}
\medskip

Since $S_\psi\subset S_\varphi\,$, hence $H^{S_\psi}\subset
H^{S_\varphi}\,$, the identity map $I$ on $H$ is Hermitian with
respect to $(\psi ,\varphi )\,$. By (1)$\,\Leftrightarrow\,$(8) in
Theorem \ref{char.herm}, (\ref{mod.transl}) and (\ref{transl-i/2}),
we obtain an $so$-continuous map
\begin{equation*}
\overline{\mathbb S_{-1/2}}\ni z\longmapsto I(z)\in B(H)\, ,
\end{equation*}
which is analytic in $\mathbb S_{-1/2}$ and satisfies the conditions
\begin{equation*}
I(s)=\Delta_M^{is}\,\Delta_N^{-is}\,\text{ and }\, I\Big(s-\frac i2
\Big) =J_M\, I(s)\, J_N\, ,\qquad\quad s\in\mathbb R\, ,
\end{equation*}
\begin{equation*} \phantom{xxxxx}\,
\| I(z)\|\leq 1\, ,\qquad z\in \overline{\mathbb S_{-1/2}}\, .
\end{equation*}
Therefore the mapping
\begin{equation*}
\overline{\mathbb S_{-1/2}}\ni \zeta\longmapsto W(\zeta )=
I(-\overline{\zeta})^*\in B(H)
\end{equation*}
is $so$-continuous, analytic in $\mathbb S_{-1/2}\,$, and such that
\begin{equation}\label{Wboundary}
W(s)=\Delta_N^{-is}\,\Delta_M^{is}\,\text{ and }\, W\Big(s-\frac i2
\Big) =J_N\, W(s)\, J_M\, ,\qquad\quad s\in\mathbb R\, ,
\end{equation}
\begin{equation}\label{Wbounded} \phantom{xxxxxx}\,
\| W(\zeta )\|\leq 1\, ,\qquad\zeta\in\overline{\mathbb S_{-1/2}}\, .
\end{equation}
\smallskip

{\it Step 2}. {\sc Hermiticity on the boundary.}
\medskip

First we show that
\begin{equation*}
\Big( W\Big(s-\frac i2\Big)\, x_\psi\,\Big|\, J_N\, y_\psi \Big)
\in\mathbb R\, ,\qquad x^*=x\in\mathfrak N_\psi\, ,\; y^*=y\in
\mathfrak N_\psi\, ,
\end{equation*}
what is equivalent, according to (1)$\,\Leftrightarrow\,$(3) in
Theorem \ref{char.herm}, to
\begin{equation}\label{Wherm.imag}
W\Big(s-\frac i2\Big)\,\text{ is Hermitian with respect to }\,
(\psi ,\psi )\, ,\qquad s\in\mathbb R\, .
\end{equation}
Indeed,
\begin{align*}
\Big( W\Big(s-\frac i2\Big)\, x_\psi\,\Big|\, J_N\, y_\psi \Big)
\overset{(\ref{Wboundary})}{=} &\big( J_N\,\Delta_N^{-is}
\Delta_M^{is}\, J_M\, x_\psi\,\big|\, J_N\, y_\psi \big)
=\big( \Delta_N^{is}\, y_\psi\,\big|\,\Delta_M^{is}\, J_M\,
x_\psi \big) \\
\overset{(\ref{mod.comm-})}{=} &\big( \Delta_N^{is}\, y_\psi\,
\big|\, J_M\, \Delta_M^{is}\, x_\varphi \big)
\overset{(\ref{delta})}{=} \big( \sigma^\psi_s(y)_\psi\,\big|\,
J_M\, \sigma^\varphi_s(x)_\varphi \big) \\
\overset{\phantom{xxx}}{=}\, & \big( \sigma^\psi_s(y)_\varphi\,
\big|\, J_M\, \sigma^\varphi_s(x)_\varphi \big)
\end{align*}
is real because of Lemma \ref{fixed}.(2).

Next we show, using the negative half-sided modular inclusion
assumption (\ref{hsmicond}), that
\begin{equation}\label{Wherm.-} \phantom{xxxxx}
\, W(s)\,\text{ is Hermitian with respect to }\, (\psi ,\psi )\, ,
\qquad s\leq 0\, ,
\end{equation}
\begin{equation}\label{Wherm.+}
J_N\, W(s)\, J_N\,\text{ is Hermitian with respect to }\,
(\psi ,\psi )\, ,\qquad s\geq 0\, .
\end{equation}

For (\ref{Wherm.-}), let $s\leq 0$ and $x^*=x\in\mathfrak N_\psi$ be
arbitrary. By (\ref{hsmicond}) we have $\sigma^\varphi_s(x)\in N\,$,
so $\sigma^\varphi_s(x)\in \mathfrak N_\psi$ and
$\sigma^\varphi_s(x)_\varphi =\sigma^\varphi_s(x)_\psi\,$. Thus
\begin{equation*}
W(s)\, x_\psi \overset{(\ref{Wboundary})}{=} \Delta_N^{-is}
\Delta_M^{is}\, x_\varphi \overset{(\ref{delta})}{=}
\Delta_N^{-is} \sigma^\varphi_s(x)_\varphi =
\Delta_N^{-is} \sigma^\varphi_s(x)_\psi \overset{(\ref{delta})}{=}
\sigma^\psi_{-s}\big(\sigma^\varphi_s(x)\big)_\psi
\end{equation*}
and the Hermiticity of $W(s)$ with respect to $(\psi ,\psi )$
follows by using (2)$\,\Rightarrow\,$(1) in
Theorem \ref{char.herm}.

Now, for (\ref{Wherm.+}), let $s\geq 0$ and $x^*=x\in\mathfrak N_\psi$
be arbitrary. Then, due to (\ref{hsmicond}), we have
\begin{equation*}
\sigma^\varphi_{-s}\big(\sigma^\psi_s(x)\big)\in \mathfrak N_\psi\,
\text{ and }\;\sigma^\varphi_{-s}\big(\sigma^\psi_s(x)\big)_\psi
\overset{(\ref{delta})}{=} \Delta_M^{-is}\Delta_N^{is}\, x_\varphi
\overset{(\ref{Wboundary})}{=} W(s)^* x_\psi\, ,
\end{equation*}
so $W(s)^* x_\psi\in H^{S_\psi}\,$. Therefore, owing to
(2)$\,\Rightarrow\,$(7) in Theorem \ref{char.herm},
$J_N\, W(s)\, J_N$ is Hermitian with respect to $(\psi ,\psi )\,$.

To summarize,
\begin{itemize}
\item[--] $W(\zeta )$ is Hermitian with respect to $(\psi ,\psi )$
for $\displaystyle \zeta\in (-\infty ,0\, ]\cup\Big(\mathbb R -
\frac i2\Big)$ and
\item[--] $J_N\, W(\zeta )\, J_N$ is Hermitian with respect to
$(\psi ,\psi )$ for $\zeta\in [0,\infty )\,$.
\end{itemize}
\bigskip

{\it Step 3}. {\sc Change of variable.}
\medskip

We shall use the analytic logarithm branches
\begin{equation*}
\log_+ : \Big\{ r\, e^{i \theta}\, ;\, r>0\, ,\, -\frac{\pi}2
<\theta <\frac{3 \pi}2 \Big\}\ni r\, e^{i \theta}\longmapsto
\log r +i\,\theta\, ,
\end{equation*}
\begin{equation*}
\log_{-} : \Big\{ r\, e^{i \theta}\, ;\, r>0\, ,\, -\frac{3 \pi}2
<\theta <\frac{\pi}2 \Big\}\ni r\, e^{i \theta}\longmapsto
\log r +i\,\theta\, .
\end{equation*}

For any $\beta\in\mathbb R\,$, $\beta\neq 0\,$, we consider (like
in Section 3, in the proof of (3)$\,\Rightarrow\,$(4) in Theorem
\ref{anal.ext}) the one point compactification of the right half
and the left half of $\overline{\mathbb S_\beta}$ and denote each
added point by $+\infty$ and $-\infty\,$, respectively.

In order to apply Theorem \ref{borchers} to $W\,$, we have
to map $\,\mathbb S_{-1/2}$ conformally onto some
$\mathbb S_\beta\,$, $\beta\neq 0\,$, such that
$\displaystyle (-\infty ,0\, )\cup\{ -\infty\}\cup\Big(\mathbb R
-\frac i2\Big)$ correspond to $\mathbb R\,$, and $(0,\infty )$ to
$\mathbb R +i\,\beta\,$. This is done, for $\beta =\pi\,$, by
\begin{equation*}
\overline{\mathbb S_{-1/2}}\,\backslash\,\{ 0\}\ni\zeta
\longmapsto\Psi (\zeta )=\log_+\big( 1-e^{2 \pi \zeta}\big)\in
\overline{\mathbb S_\pi}\,\backslash\,\{ 0\}\, ,
\end{equation*}
which extends to a homeomorphism $\Psi : \overline{\mathbb S_{-1/2}}
\cup\{-\infty ,+\infty\}\longrightarrow\overline{\mathbb S_\pi}
\cup\{-\infty ,+\infty\}$ satisfying $\Psi (0)=-\infty\,$,
$\Psi (-\infty )=0\,$ and $\Psi (+\infty )=+\infty\,$.

The inverse homeomorphism $\Psi^{-1} : \overline{\mathbb S_\pi}
\cup\{-\infty ,+\infty\}\longrightarrow\overline{\mathbb S_{-1/2}}
\cup\{-\infty ,+\infty\}$ is given by
\begin{equation*}
\overline{\mathbb S_\pi}\,\backslash\,\{ 0\}\ni z\longmapsto
\Psi^{-1} (z)=\frac 1{2 \pi} \log_-\big( 1-e^z\big)\in
\overline{\mathbb S_{-1/2}}\,\backslash\,\{ 0\}\, ,
\end{equation*}
\begin{equation*}
\Psi^{-1} (0)=-\infty\, ,\,\Psi^{-1} (-\infty )=0\, ,\,\Psi^{-1}
(+\infty )=+\infty\, ,
\end{equation*}
so it maps $\,\mathbb S_\pi$ conformally onto $\mathbb S_{-1/2}$
and
\begin{equation*}
(-\infty ,0)\,\text{ onto }\, (-\infty ,0)\qquad \Big(\,
0{\scriptscriptstyle -}\mapsto -\infty\, ,\;\;
-\infty\mapsto 0{\scriptscriptstyle -}\Big)\, ,
\end{equation*}
\begin{equation*}
(0,+\infty )\,\text{ onto }\,\mathbb R -\frac i2\qquad \Big(\,
0{\scriptscriptstyle +}\mapsto -\infty -\frac i2\, ,\;\;
+\infty\mapsto +\infty +\frac i2\Big)\, ,
\end{equation*}
\begin{equation*}
\mathbb R +i \pi\,\text{ onto }\, (0,+\infty )\qquad \Big(\,
-\infty +i \pi\mapsto 0{\scriptscriptstyle +}\, ,\;\;
i \pi\mapsto \frac{\log 2}{2\,\pi}\, ,\;\;
+\infty +i \pi\mapsto +\infty\Big)\, .
\end{equation*}
Thus we can consider the $so$-continuous mapping
\begin{equation}\label{map}
\overline{\mathbb S_\pi}\,\backslash\,\{ 0\}\ni z\longmapsto
V(z)=W\big( \Psi^{-1}(z)\big)\in B(H)\, ,
\end{equation}
which is analytic in $\mathbb S_\pi$ and, according to
(\ref{Wbounded}), (\ref{Wherm.imag}), (\ref{Wherm.-}) and
(\ref{Wherm.+}), satisfies
\begin{equation}\label{Vbounded} \hspace{5cm}
\| V(z)\|\leq 1\, ,\qquad z\in\overline{\mathbb S_\pi}\,
\backslash\,\{ 0\}\, ,
\end{equation}
\begin{equation}\label{Vherm.} \hspace{2.05cm}
V(s)\,\text{ is Hermitian with respect to }\, (\psi ,\psi )\, ,
\qquad s\in\mathbb R\,\backslash\,\{ 0\}\, ,
\end{equation}
\begin{equation}\label{Vherm.imag} \hspace{0.97cm}
J_N\, V(s+i\pi )\, J_N\,\text{ is Hermitian with respect to }\,
(\psi ,\psi )\, ,\qquad s\in\mathbb R\, .
\end{equation}
\smallskip

{\it Step 4}. {\sc Application of the Generalized Structure
Theorem.}
\medskip

Since the $so$-continuous mapping considered in (\ref{map}) is
analytic in $\,\mathbb S_\pi$ and (\ref{Vbounded}), (\ref{Vherm.}),
(\ref{Vherm.imag}) hold, it satisfies the
assumptions for $T(\,\cdot\, )$ in Theorem \ref{borchers} with
\begin{center}
$M$ replaced by $N\,$, $\varphi$ replaced by $\psi\,$,
$\beta =\pi\,$, $\Xi_o=\{ 0\}$ and $\Xi_1=\emptyset\,$.
\end{center}
By Theorem \ref{borchers} it follows that, for some $V\in B(H)$
which is Hermitian with respect to $(\psi ,\psi )\,$,
\begin{equation}\label{Vformula}
V(s)=\Delta_N^{\textstyle -i\frac s{2\pi}}\, V\,
\Delta_N^{\textstyle i\frac s{2\pi}}\, ,\qquad s\in\mathbb R\,
\backslash\,\{ 0\}
\end{equation}
and the mapping (\ref{map}) has an $so$-continuous extension
\begin{equation*}
\overline{\mathbb S_\pi}\ni z\longmapsto V(z)\in B(H)
\end{equation*}
with $V(0)=V\,$, satisfying
\begin{equation}\label{Vtransl}
V(z+2\pi t)=\Delta_N^{-it}\, V(z)\,\Delta_N^{it}\, ,
\qquad z\in\overline{\mathbb S_\pi}\, ,\; t\in\mathbb R\, ,
\end{equation}
\begin{equation}\label{transl+ipi}
V(s+i\pi )= J_N\, V(s)\, J_N\, ,\qquad s\in\mathbb R\, .
\;\phantom{xxx}
\end{equation}
In particular,
\begin{equation*}
V(0) = so-\lim_{0\neq s\to 0}V(s) =
so-\lim_{t\to -\infty} W(t)=
so-\lim_{t\to -\infty} \Delta_N^{-it}\,\Delta_M^{it}
\end{equation*}
is a wave operator.

We notice an additional continuity property of $V(\,\cdot\, )\,$:
since $\Psi^{-1}(-\infty )=0\,$, the limit
\begin{equation*}
so-\lim_{\overline{\mathbb S_\pi}\ni z\to -\infty} V(z) =
so-\lim_{\overline{\mathbb S_{-1/2}}\ni \zeta\to 0} W(\zeta ) =
W(0)=1_H
\end{equation*}
exists. Thus the mapping (\ref{map}) has actually an $so$-continuous
extension
\begin{equation}\label{fullmap}
\overline{\mathbb S_\pi}\cup\{ -\infty\}\ni z\longmapsto V(z)
\in B(H)
\end{equation}
with $\displaystyle V(0)=so-\lim_{t\to -\infty} \Delta_N^{-it}\,
\Delta_M^{it}$ and $V(-\infty )=1_H\,$.
\bigskip

{\it Step 5}. {\sc Further change of variable.}
\medskip

We recall that
\begin{equation*}
\overline{\mathbb S_\pi}\ni z\longmapsto \Theta (z)=e^z\in
\{\zeta\in\mathbb C\; ;\, \Im\zeta \geq 0\}\,\backslash\,\{ 0\}
\end{equation*}
is a homeomorphism mapping $\,\mathbb S_\pi$ conformally onto
$\{\zeta\in\mathbb C\, ;\, \Im\zeta >0\}\,$, which extends to a
homeomorphism $\Theta : \overline{\mathbb S_\pi}\cup\{ -\infty
,+\infty\}\longrightarrow\{\zeta\in\mathbb C\, ;\,\Im\zeta\geq 0\}
\cup\{\infty\}$ satisfying $\Theta (-\infty )=0$ and $\Theta
(+\infty )=\infty\,$.

The inverse homeomorphism $\Theta^{-1} : \{\zeta\in\mathbb C\,
;\,\Im\zeta\geq 0\}\cup\{\infty\}\longrightarrow
\overline{\mathbb S_\pi}\cup\{ -\infty ,+\infty\}\,$,
which maps $\{\zeta\in\mathbb C\, ;\, \Im\zeta >0\}$ conformally
onto $\mathbb S_\pi\,$, is given by
\begin{equation*}
\{\zeta\in\mathbb C\, ;\,\Im\zeta\geq 0\}\,\backslash\,\{ 0\}
\ni \zeta\longmapsto \Theta^{-1} (\zeta )=\log_+ \zeta\in
\overline{\mathbb S_\pi}\, ,
\end{equation*}
\begin{equation*}
\Theta^{-1}(0)=-\infty\, ,\quad \Theta^{-1}(\infty )=+\infty\, .
\end{equation*}
Since the mapping (\ref{fullmap}) is $so$-continuous, also the
mapping
\begin{equation*}
\{\zeta\in\mathbb C\, ;\,\Im\zeta\geq 0\}\ni \zeta\longmapsto
U(\zeta )=V\big(\Theta^{-1}(\zeta )\big)\in B(H)
\end{equation*}
is $so$-continuous. Moreover, since $\,\mathbb S_\pi\ni z
\longmapsto V(z)$ is analytic, the restriction of the above
mapping to $\{\zeta\in\mathbb C\, ;\, \Im\zeta >0\}$ is analytic.

For $\zeta =0$ and $\zeta =1$ we have
\begin{equation*}
U(0)=V(-\infty )=1_H\, ,\qquad U(1)=V(0)=so-\lim_{t\to -\infty}
\Delta_N^{-it}\,\Delta_M^{it}\, .
\end{equation*}
In particular, $U(0)$ and $U(1)$ are unitaries and
$\Delta_N^{is}\, U(1)=U(1)\,\Delta_M^{is}$ for all $s\in
\mathbb R\,$.
On the other hand, for $\zeta\in\mathbb C\, ,\,\Im \zeta\geq 0\,$,
$\zeta\neq 0\, , 1\,$,
\begin{equation*}
U(\zeta )=
V\big(\Theta^{-1}(\zeta )\big) =
V\big(\log_+ \zeta \big) =W\big(\Psi^{-1}(\log_+ \zeta )\big) =
W\left(\, \frac 1{2 \pi} \log_- (1-\zeta )\right)
\end{equation*}
holds. In particular, according to (\ref{Wboundary}), also the
operators
\begin{equation*}
\{ U(s)\, ;\, s\in\mathbb R\, ,\, s\neq 0\, , 1\} =
\{ V(z)\, ;\, z\in\partial\,\mathbb S_\pi\, ,\, z\neq 0\} =
\{ W(\zeta )\, ;\,\zeta\in\partial\,\mathbb S_{-i/2}\, ,\,\zeta
\neq 0\}
\end{equation*}
are unitaries.

We summarize:
\begin{align}
&U(0)=1_H\,\text{ and }\, U(s)\,\text{ is unitary for every }\,
s\in\mathbb R\, , \notag \\
&U(1)=so-\lim_{t\to -\infty}\Delta_N^{-it}\,\Delta_M^{it}\,
\text{ and }\;\Delta_N^{is}\, U(1)=U(1)\,\Delta_M^{is}\, ,\qquad
s\in\mathbb R\,  , \label{wave} \\
&U(\zeta ) =V\big( \log_+ \zeta \big)\, ,\qquad 0\neq \zeta \in
\mathbb C\, ,\, \Im \zeta\geq 0\, , \label{UbyV} \\
&U(\zeta ) =W\left(\, \frac 1{2 \pi} \log_- (1-\zeta )\right)\, ,
\qquad 1\neq \zeta\in\mathbb C\, ,\, \Im \zeta\geq 0\, .
\label{UbyW}
\end{align}

Using (\ref{UbyV}), we obtain from (\ref{Vbounded}),
(\ref{Vtransl}) and (\ref{transl+ipi})
\begin{align}
&\| U(\zeta )\|\leq 1\, ,\qquad \zeta\in\mathbb C\, ,\, \Im \zeta
\geq 0\, , \label{Ubounded} \\
&U\big( e^{2 \pi t} \zeta\big) =\Delta_N^{-it}\, U(\zeta )\,
\Delta_N^{it}\, ,\qquad t\in\mathbb R\, ,\;\zeta\in\mathbb C\, ,\,
\Im \zeta\geq 0\, , \label{Utransl} \\
&U(-s)=J_N\, U(s)\, J_N\, ,\qquad s\in\mathbb R\, .
\label{transl+i}
\end{align}
Indeed, (\ref{Vbounded}) implies that $\| U(\zeta )\|\leq 1$ for
$0\neq \zeta\in\mathbb C\, ,\, \Im \zeta\geq 0\,$, while the norm
of $U(0)=1_H$ is $\,\leq1\,$. Similarly, the equality in
(\ref{Utransl}) is an immediate consequence of (\ref{Vtransl}) for
$0\neq \zeta\in \mathbb C\, ,\, \Im \zeta\geq 0\,$, while it is
trivial for $\zeta =0\,$. Finally, for any $s>0\,$, (\ref{transl+ipi})
implies
\begin{equation*}
U(-s)= V( \log s +i\,\pi )= J_N\, V( \log s)\, J_N =J_N\, U(s)\,
J_N\, ,
\end{equation*}
hence also $\, J_N\, U(-s)\, J_N =J_N^{\, 2}\, U(s)\, J_N^{\, 2}
=U(s)\,$. Therefore the equality in (\ref{transl+i}) holds for
every $s\in\mathbb R$ (it is trivial for $s=0)\,$.

Furthermore, by (\ref{UbyW}) and (\ref{Wboundary}),
\begin{equation}\label{2}
U(2)=W\left(\, \frac 1{2 \pi} \log_- (-1)\right) =
W\left( -\frac i2 \right) =J_N\, W(0)\, J_M =J_N\, J_M\, .
\end{equation}
On the other hand, (\ref{UbyW}) is equivalent to the equality
\begin{equation*}
W(z)=U\big( 1-e^{2 \pi z}\big)\, ,\qquad z\in
\overline{\mathbb S_{-i/2}}\, ,
\end{equation*}
which yields
\begin{equation}\label{WbyU}
\Delta_N^{-it}\,\Delta_M^{it} =W(t) =U\big( 1-e^{2 \pi t}\big)\, ,
\qquad t\in\mathbb R\, .
\end{equation}
\smallskip

{\it Step 6}. {\sc Group property of $U(\,\cdot\, )\,$.}
\medskip

We now prove the group property
\begin{equation}\label{group}
U(s_1)\, U(s_2)=U(s_1+s_2)\, ,\qquad s_1\, ,\, s_2\in\mathbb R\, .
\end{equation}

Let $s>0$ and $t\in\mathbb R$ be arbitrary.
By (\ref{UbyV}), (\ref{Vformula}). (\ref{wave}) and (\ref{WbyU}),
we obtain
\begin{align*}
U(s) &=V( \log s)=\Delta_N^{\textstyle -i\frac {\log s}{2\pi}}
V(0)\,\Delta_N^{\textstyle i\frac {\log s}{2\pi}} =
\Delta_N^{\textstyle -i\frac {\log s}{2\pi}} U(1)\,
\Delta_N^{\textstyle i\frac {\log s}{2\pi}} = \\
&=U(1)\,\Delta_M^{\textstyle -i\frac {\log s}{2\pi}}
\Delta_N^{\textstyle i\frac {\log s}{2\pi}} =U(1)\; W\!\left(
\frac {\log s}{2\pi}\right)^{\! *} =U(1)\, U(1-s)^*\, ,
\end{align*}
hence $\, U(s)\, U(1-s) =U(1)\,$. By sandwiching this equation
by $\text{Ad}\; \Delta_N^{- i t}$ and taking into account
(\ref{Utransl}), we get
\begin{equation}\label{1-s}
U\big( e^{2 \pi t} s\big)\, U\big( e^{2 \pi t} (1-s)\big) =
U\big( e^{2 \pi t} \big)\, .
\end{equation}

Next let $r_1\, ,\, r_2\in\mathbb R$ be such that $r_1>0$ and
$r_1+r_2>0\,$. Then $\displaystyle s=\frac{r_1}{r_1+r_2}>0$ and,
with $\displaystyle t=\frac 1{2 \pi} \log\, (r_1+r_2)\in
\mathbb R\,$, we have $e^{2 \pi t} s=r_1\, ,\, e^{2 \pi t} (1-s)
=r_2\,$. Thus (\ref{1-s}) yields
\begin{equation}\label{semigroup}
U(r_1)\, U(r_2)=U(r_1+r_2)\, .
\end{equation}

Finally let $s_1\, ,\, s_2\in\mathbb R$ be arbitrary and choose
$s\in\mathbb R$ such that $s>0\,$, $s+s_1>0$ and $s+s_1+s_2>0\,$.
Then, using (\ref{semigroup}) with $(r_1 ,r_2)$ equal to
$(s\, ,s_1)\,$, $(s+s_1 ,s_2)$ and $(s\, ,s_1+s_2)\,$,
respectively, we obtain
\begin{equation*}
U(s)\, U(s_1)\, U(s_2) =U(s+s_1)\, U(s_2) =U(s+s_1+s_2) =
U(s)\, U(s_1+s_2)\, .
\end{equation*}
Since $U(s)$ is unitary, the above equality implies that
$U(s_1)\, U(s_2)=U(s_1+s_2)\,$.

Therefore (\ref{group}) is proved. In particular,
\begin{equation}\label{inverse}
U(s)\, U(-s)=U(0)=1_H\, ,\text{ that is }\, U(s)^*=U(-s)\, ,\qquad
s\in\mathbb R\, .
\end{equation}

Thus $\mathbb R\ni s\longmapsto U(s)\in B(H)$ is an $so$-continuous
one-parameter group of unitaries, which allows an $so$-continuous
extension $\{\zeta\in\mathbb C\, ;\,\Im \zeta\geq 0\}\ni\zeta
\longmapsto U(\zeta )\,$,
analytic in $\{\zeta\in\mathbb C\, ;\,\Im \zeta >0\}$ and
satisfying (\ref{Ubounded}). Consequently
\begin{equation}\label{P}
U(s) =\exp (i s P)\, ,\qquad s\in \mathbb R
\end{equation}
for some positive selfadjoint operator $P$ in $H\,$.
\bigskip

{\it Step 7}. {\sc Further properties of $U(\,\cdot\, )\,$.}
\medskip

Here we show that the above constructed group $\mathbb R\ni s
\longmapsto U(s)\in B(H)$ satisfies properties (1) - (7) in
Theorem \ref{hsmi}.

By (\ref{Utransl}), (\ref{WbyU}), (\ref{group}) and (\ref{inverse}),
we obtain for all $s\, ,\, t\in\mathbb R$
\begin{equation}\label{Utransl.bis}
\begin{split}
\Delta_M^{-it}\, U(s)\,\Delta_M^{it} &=
\big( \Delta_M^{-it} \Delta_N^{it}\big)\,\Delta_N^{-it}\, U(s)\,
\Delta_N^{it}\big( \Delta_N^{-it} \Delta_M^{it}\big) = \\
&=U\big( 1-e^{2 \pi t}\big)^{\! *}\, U\big( e^{2 \pi t} s\big)\,
U\big( 1-e^{2 \pi t}\big) =U\big( e^{2 \pi t} s\big)\, .
\end{split}
\end{equation}
This equality and (\ref{Utransl}) show that property (1) in
Theorem \ref{hsmi} is satisfied.

Similarly, (\ref{transl+i}), (\ref{2}), (\ref{group}) and
(\ref{inverse}) yield for every $s\in\mathbb R$
\begin{equation}\label{Mtransl+i}
\begin{split}
J_M\, U(s)\, J_M &=\big( J_M\, J_N\big) J_N\, U(s)\, J_N\big(
J_N\, J_M\big) \\
&=U(2)^*\, U(-s)\, U(2) =U(-s)\, . \hspace{2.4 cm}
\end{split}
\end{equation}
Now property (2) in Theorem \ref{hsmi} is (\ref{transl+i}) together
with (\ref{Mtransl+i}).

The validity of property (3) in Theorem \ref{hsmi} follows from
(\ref{WbyU}) and (\ref{wave}).

The first equality in property (4) in Theorem \ref{hsmi} is
(\ref{2}), while the second one follows from (\ref{inverse}),
(\ref{Mtransl+i}), (\ref{group}) and (\ref{2}) :
\begin{equation}\label{1}
U(1)\, J_M\, U(1)^* =U(1)\, J_M\, U(-1) =U(1)^2\, J_M =U(2)\, J_M
=J_N\, .
\end{equation}

Next we prove property (5). Since
\begin{equation*}
\Delta_M^{-i t}\,\Delta_N^{i t}\, N\,\Delta_N^{-i t}\,\Delta_M^{i t}
=\Delta_M^{-i t}\, N\,\Delta_M^{i t}\subset
\Delta_M^{-i t}\, M\,\Delta_M^{i t} =M\, ,\qquad t\in\mathbb R\, ,
\end{equation*}
(\ref{wave}) implies
\begin{equation}\label{overN}
U(1)^* N\, U(1)\subset M\, .
\end{equation}
On the other hand, by sandwiching this relation by $\text{Ad}\, J_M$
and using (\ref{1}), we obtain
\begin{equation*}
M' =J_M\, M\, J_M\supset J_M\, U(1)^* N\, U(1)\, J_M =
U(1)^* J_N\, N\, J_N\, U(1) =U(1)^* N'\, U(1)\, .
\end{equation*}
Passing to the commutants, this inclusion relation yields
\begin{equation*}
M\subset U(1)^* N\, U(1)\, ,
\end{equation*}
which together with (\ref{overN}) gives $U(1)^* N\, U(1)=M\,$,
that is
\begin{equation}\label{equalN}
N=U(1)\, M\, U(1)^*\, .
\end{equation}

For property (6) in Theorem \ref{hsmi} we notice that
(\ref{Wboundary}) and the negative half-sided modular inclusion
assumption (\ref{hsmicond}) imply that $W(t)\, N\, W(t)^*\subset N$
for all $t\leq 0\,$, what is by (\ref{WbyU}) equivalent to
\begin{equation*}
U(s)\, N\, U(s)^*\subset N\, ,\qquad 0\leq s\leq 1\, .
\end{equation*}
Using (\ref{equalN}), (\ref{group}) and (\ref{inverse}), we get
for every $0\leq s\leq 1$
\begin{equation*}
\begin{split}
U(s)\, M\, U(s)^* &=U(s)\, U(1)^* N\, U(1)\, U(s)^* =
U(1)^* U(s)\, N\, U(s)^* U(1) \\
&\subset U(1)^* N\, U(1)=M\, .
\end{split}
\end{equation*}
Using now induction on $n\,$, it follows that
\begin{equation*}
U(s)\, M\, U(s)^*\subset M\, ,\qquad 0\leq s\leq n
\end{equation*}
holds for every integer $n\geq 1\,$, that is $U(s)\, M\, U(s)^*
\subset M$ for all $s\geq 0\,$.

(7) is an immediate consequence of (6), (5) and (4).
\bigskip

{\it Step 8}. {\sc Invariance properties of $U(\,\cdot\, )$ with
respect to $\varphi\,$.}
\medskip

We show in the following that $\mathbb R\ni s\longmapsto U(s)\in B(H)$
satisfies property (8) in Theorem \ref{hsmi}, hence also property
(9), which is an immediate consequence of (8).

For any $y\in\mathfrak N_\psi$ and $t\in\mathbb R\,$, (\ref{delta})
yields
\begin{equation*}
\Delta_N^{it}\, y\,\Delta_N^{-it}\in\mathfrak N_\psi\subset
\mathfrak N_\varphi\,\text{ and }\,\big(\Delta_N^{it}\, y\,
\Delta_N^{-it}\big)_{\!\varphi} =\big(\Delta_N^{it}\, y\,
\Delta_N^{-it}\big)_{\!\psi} =\Delta_N^{it}\, y_\psi\, .
\end{equation*}
Setting $s=1-e^{2 \pi t}$ and using (\ref{WbyU}), we obtain
\begin{equation*}
U(s)^* y\, U(s) =\Delta_M^{-it}\,\Delta_N^{it}\, y\,\Delta_N^{-it}\,
\Delta_M^{it}\in\Delta_M^{-it}\,\mathfrak N_\psi\,\Delta_M^{it}\subset
\Delta_M^{-it}\,\mathfrak N_\varphi\,\Delta_M^{it}\subset
\mathfrak N_\varphi\, ,
\end{equation*}
\begin{align*}
\big(\, U(s)^* y\, U(s)\big)_{\!\varphi} &=\Delta_M^{-it}\big(
\Delta_N^{it}\, y\,\Delta_N^{-it}\big)_{\!\varphi} =\Delta_M^{-it}
\big(\Delta_N^{it}\, y\,\Delta_N^{-it}\big)_{\!\psi} =\Delta_M^{-it}
\Delta_N^{it}\, y_\psi \phantom{xxxxxx} \\ &=U(s)^* y_\psi\, .
\end{align*}
Therefore we have for all $s\in\mathbb R\, ,\, s<1\,$,
\begin{equation}\label{impl1}
U(s)^* y\, U(s)\in\mathfrak N_\varphi\,\text{ and }\,
\big(\, U(s)^* y\, U(s)\big)_{\!\varphi} =U(s)^* y_\psi\, .
\end{equation}
Moreover, according to the Lebesgue continuity result Proposition
\ref{Leb}, (\ref{impl1}) holds also for $s=1\,$.

According to (\ref{equalN}), $\gamma_1 : M\ni x\longmapsto U(1)\, x
\, U(1)^*\in N$ is a $*$-isomorphism with inverse $\gamma_1^{\, -1} :
N\ni y\longmapsto U(1)^* y\, U(1)\in M\,$. The modular automorphism
groups of the normal semifinite faithful weights $\psi$ and
$\varphi\circ\gamma_1^{\, -1}$ are equal. Indeed, by (\ref{wave}) we
have for every $y\in N$ and $t\in\mathbb R\,$:
\begin{align*}
\sigma^{\varphi\circ\gamma_1^{\, -1}}_t (y)=
\gamma_1\Big( \sigma^\varphi_t \big( \gamma_1^{\, -1} (y)\big)\Big)
&=U(1)\,\Delta_M^{it}\, U(1)^* y\, U(1)\,\Delta_M^{-it}\, U(1)^* \\
&=\Delta_N^{it}\, y\,\Delta_N^{-it} =\sigma^\psi_t (y)\, .
\end{align*}
On the other hand, for every $y\in\mathfrak N_\psi\,$, using
(\ref{impl1}) with $s=1\,$, we obtain
\begin{equation*}
\psi (y^*y)=\| y_\psi\|^2 =\|\, U(1)^* y_\psi\|^2 =\|\gamma_1^{\, -1}
(y)_\varphi\|^2 =\varphi\circ\gamma_1^{\, -1}(y^*y)\, ,
\end{equation*}
so $\varphi\circ\gamma_1^{\, -1}$ and $\psi$ coincide on
$\mathfrak M_\psi\,$. Thus \cite{P-T}, Proposition 5.9 yields
$\psi =\varphi\circ\gamma_1^{\, -1}\,$, that is
\begin{equation}\label{1inv}
\varphi\circ\gamma_1 =\psi\circ\gamma_1 =\varphi\, .
\end{equation}
In particular, for every $x\in M$ and $n\geq 0\,$,
\begin{equation}\label{finiteweight}
x\in\mathfrak N_\varphi\,\Longleftrightarrow\,\gamma_n (x)\in
\mathfrak N_\varphi\, .
\end{equation}

Let $x\in\mathfrak N_\varphi$ be arbitrary. By
(\ref{finiteweight}) and (\ref{equalN}) we have $\gamma_1 (x)\in
\mathfrak N_\psi\,$, so (\ref{impl1}) holds with $y=\gamma_1 (x)$
and any $0\leq s\leq 1\,$. Using the group property of
$U(\,\cdot\, )\,$, we deduce that, for every $0\leq s\leq 1\,$,
\begin{equation}\label{impl2}
\begin{split}
&\gamma_s (x) = U(s)\, x\, U(s)^* =
U(1-s)^* \gamma_1 (x)\, U(1-s)\in\mathfrak N_\varphi\,\text{ and} \\
&\gamma_s (x)_\varphi =U(1-s)^* \gamma_1 (x)_\psi\, .
\end{split}
\end{equation}
In particular, for $s=0$ we get $x_\varphi =U(1)^* \gamma_1 (x)_\psi
\,$, that is $\gamma_1 (x)_\psi =U(1)\, x_\varphi\,$. Thus
(\ref{impl2}) yields
\begin{equation}\label{impl3}
\gamma_s (x)\in\mathfrak N_\varphi\,\text{ and }\,\gamma_s
(x)_\varphi =U(s)\, x_\varphi\, ,\qquad 0\leq s\leq 1\, .
\end{equation}
Iterating (\ref{impl3}), we obtain
\begin{equation}\label{impl4}
x\in\mathfrak N_\varphi\, ,\; s\geq 0\;\Longrightarrow\,
\gamma_s (x)\in\mathfrak N_\varphi\,\text{ and }\,\gamma_s
(x)_\varphi =U(s)\, x_\varphi\, .
\end{equation}

On the other hand, for $x\in M$ and $s\geq 0\,$, denoting by $n$
the integer part of $s\,$, that is the integer $n\geq 0$ with
$n\leq s<n+1\,$, we have
\begin{equation*}
\gamma_s (x)\in\mathfrak N_\varphi\;
\overset{(\ref{impl4})}{\Longrightarrow}\; \gamma_{n+1}
(x)=\gamma_{n+1-s}\big(\gamma_s (x)\big)\in\mathfrak N_\varphi\;
\overset{(\ref{finiteweight})}{\Longrightarrow}\; x\in
\mathfrak N_\varphi\, .
\end{equation*}
Consequently property (8) in Theorem \ref{hsmi} is satisfied.
\bigskip

{\it Step 9}. {\sc Description of the generator $P\,$.}
\medskip

First we verify that statement (10) in Theorem \ref{hsmi} holds
with $P$ defined by (\ref{P}).

We recall from Subsection (b) of Section 2 that, if we endow
$\mathbb R^2$ with the Lie group structure defined by the
composition law
\begin{equation*}
(s_1 ,t_1)\cdot (s_2 ,t_2) = (s_1 +e^{-2\pi t_1} s_2\, ,\,
t_1+t_2 )\, ,
\end{equation*}
then $\mathbb R^2\ni (s ,t)\longmapsto T_s\, L_t\in
\mathcal P_+^\uparrow (1)$ is a Lie group isomorphism. Hence,
by (\ref{Utransl.bis}),
\begin{equation*}
\pi : \mathcal P_+^\uparrow (1)\ni T_s\, L_t\longmapsto
U(s)\,\Delta_M^{i t}
\end{equation*}
is an $so$-continuous unitary representation on $H$ and, according
to (\ref{WbyU}), the group $\pi \big(\mathcal P_+^\uparrow (1)\big)$
contains $\{ \Delta_M^{it}\, ,\,\Delta_N^{is}\, ;\, t\, ,\, s\in
\mathbb R\}$ and is generated by this set.

Let us consider the elements $X_1\, ,\, X_2\, ,\, X_3$ of the Lie
algebra $\mathfrak p_+^\uparrow (1)\equiv\mathfrak g$ defined in
(\ref{gener}). By (\ref{gen.exp}) and by the definition of $\pi\,$,
taking into account how $\mathcal P_+^\uparrow (1)$ was identified in
Subsection (b) of Section 2 with $\mathcal G\,$, we obtain for every
$t\in\mathbb R\,$:
\begin{align*}
&\pi\big( \exp (t\, X_1)\big) =\pi (L_t) =\Delta_M^{i t}=
\exp\big( i t \log \Delta_M \big)\, , \\
&\pi\big( \exp (t\, X_2)\big) =\pi \big( T_{1-e^{-2\pi t}}\, L_t\big)
=U\big( 1-e^{-2\pi t}\big)\,\Delta_M^{i t}
\overset{(\ref{WbyU})}{=}\Delta_N^{i t}=
\exp\big( i t \log \Delta_N \big)\, , \\
&\pi\big( \exp (t\, X_3)\big) =\pi (T_t) =U(t)
\overset{(\ref{P})}{=}\exp (i t P)\, .
\end{align*}
Therefore, according to (\ref{closure}),
\begin{equation}\label{mapping}
\overline{\text{\rm d}\pi (X_1)}=i\,\log \Delta_M\, ,\quad
\overline{\text{\rm d}\pi (X_2)}=i\,\log \Delta_N\, ,\quad
\overline{\text{\rm d}\pi (X_3)}=i\, P\, .
\end{equation}

Since any two of $X_1\, ,\, X_2\, ,\, X_3$ is a basis for
$\mathfrak g_+^\uparrow (1)$ and $\mathcal P_+^\uparrow (1)$ is
connected and simply connected, the representation $\pi$ is uniquely
determined by any two of the relations (\ref{mapping}) (see e.g.
\cite{Ba-Ra}, Ch. 11, \S$\,$5 or \cite{Schm}, Proposition 10.5.2).

Now, by (\ref{commutation}), (\ref{sumLie}) and (\ref{mapping}),
we conclude that
\begin{align*}
i P &=\overline{\text{\rm d}\pi (X_3)} =
\frac 1{2\pi}\,\overline{\text{\rm d}\pi (X_2 -X_1)}=
\frac 1{2\pi}\,\big(\, \overline{\overline{\text{\rm d}\pi (X_2)}
-\overline{\text{\rm d}\pi (X_1)}}\, \big) \\
&= \frac i{2\pi}\,\big(\, \overline{\log \Delta_N -\log \Delta_M}
\,\big)\, ,
\end{align*}
hence $P$ is the closure of $\,\displaystyle \frac 1{2\pi}\,
\big(\,\log \Delta_N -\log \Delta_M \big)\,$.

\hfill $\square\quad$

\bigskip
\section{Complements to the implementation theorem of Borchers
and the proof of Theorem \ref{type}}
\bigskip

First we prove Theorem \ref{pos.gen},
which will then be used to prove Theorem \ref{type}.
\medskip

{\it Proof of }$\,${\bf Theorem \ref{pos.gen}}.
\medskip

{\it Step 1}. {\sc The existence and the uniqueness of $b\,$}
(it is essentially the proof of \cite{Arv}, Theorem 3.1 and
\cite{Z2}, Corollary 5.7).
\medskip

By the lower boundedness of $P$ we have $d_o=\exp (-P)\in B(H)\,$.
Moreover, $d_o$ is clearly positive and injective. Denoting
$\beta_s =\alpha_{-s}\,$, $(\beta_s )_{s\in\mathbb R}$ is
an $so$-continuous group of $*$-automorphisms of $M$ such that
\begin{equation*}
\beta_s (x) ={d_o}^{is} x\, {d_o}^{-is}\, ,\qquad s\in\mathbb R
\, ,\, x\in M\, .
\end{equation*}

We recall that, according to \cite{Z2}, Theorem 1.4, we have for
every $\lambda\in\mathbb R\,$:
\begin{equation}\label{spectr.subsp.}
\begin{split}
&M^\alpha\big( [-\lambda , +\infty )\big) =
M^\beta\big( (-\infty ,\lambda ]\big) = \\
=\;&\Big\{ x\in \!\bigcap_{\substack{{z\in\mathbb C}\\{\Im z\geq 0}}}\!
{\rm Dom}\, (\alpha_z )\; ;\, \|\alpha_z (x)\|\leq e^{\lambda \Im z}
\| x\|\text{ for all } z\in\mathbb C\, ,\Im z\geq 0 \Big\} = \\
=\;&\Big\{ x\in \!\bigcap_{\substack{{z\in\mathbb C}\\{\Im z\leq 0}}}\!
{\rm Dom}\, (\beta_z )\; ;\, \|\beta_z (x)\|\leq e^{-\lambda \Im z}
\| x\|\text{ for all } z\in\mathbb C\, ,\Im z\leq 0 \Big\}\, .
\end{split}
\end{equation}
Denoting
\begin{align*}
H_\lambda =\; &\text{the closed linear span of }M^\alpha\big(
[-\lambda ,+\infty )\big)\, H = \\
=\; &\text{the closed linear span of }M^\beta\big(
(-\infty ,\lambda ]\big)\, H\, ,
\end{align*}
we have clearly $H_\lambda \supset M^\beta\big((-\infty ,0]\big)
\, H\supset 1_H H = H\,$, hence $H_\lambda =H\,$, for all $\lambda
\geq 0\,$. In particular, $H$ is an
\emph{invariant subspace of support} $1_H$ \emph{relative to}
$\beta\,$, as defined in Section 5 of \cite{Z2}.
Moreover, since the spectral subspace of $d_o$ corresponding to
$(0,\| d_o\| ]$ is $H\,$, the second statement in \cite{Z2},
Theorem 5.3 implies that $H$ is \emph{simply invariant}, that is
$\displaystyle \bigcap_{\lambda\in\mathbb R} H_\lambda =\{ 0\}\,$.
Furthermore, $M' H_\lambda\subset H_\lambda$ implies that the
orthogonal projection $p_\lambda$ onto $H_\lambda$ belongs to $M\,$.

Using now the first statement in \cite{Z2}, Theorem 5.3, it follows
that there exists an injective $b\in B(H)\, ,\, 0\leq b\leq 1_H\,$,
such that
\begin{equation}\label{impl.}
\beta_s (x) =b^{is} x\, b^{-is}\, ,\qquad s\in\mathbb R\, ,\,
x\in M
\end{equation}
and, for every $\lambda\in\mathbb R\,$, the spectral projection
$\chi_{\substack{{}\\{(0,e^\lambda ]}}} (b)$ is the orthogonal
projection onto $\displaystyle H_{\lambda +0} =\bigcap_{\mu >
\lambda} H_\mu\,$, hence it is equal to $p_{\lambda +0} =
so-\lim\limits_{\lambda <\mu\to\lambda}\, p_\mu\in M\,$.
In particular, $b\in M\,$.

Property (i) in the statement of Theorem \ref{pos.gen} holds by
(\ref{impl.}). In order to verify property (ii), let $d$ be
an arbitrary injective operator in $M$ such that $0\leq d\leq 1_H$
and $\alpha_s (x)=d^{-i s} x\, d^{i s} ,\, s\in\mathbb R\, ,\,
x\in M\,$, that is $\beta_s (x)=d^{i s} x\, d^{-i s} ,\, s\in
\mathbb R\, ,\, x\in M\,$. Since the spectral subspace of the
unitary group $(d^{i s})_{s\in\mathbb R}$ corresponding to
$(-\infty ,\lambda ]$ is
$\chi_{\substack{{}\\{(0,e^\lambda ]}}} (d)\, H\,$, \cite{Z2},
Corollary 2.6 yields
\begin{equation*}
M^\beta\big( (-\infty ,\lambda ]\big)\, H =
M^\beta\big( (-\infty ,\lambda ]\big)\,
\underbrace{\chi_{\substack{{}\\{(0,e^0 ]}}} (d)}_{=\, 1_H}\, H
\subset\chi_{\substack{{}\\{(0,e^\lambda ]}}} (d)\, H\, ,\qquad
\lambda\in\mathbb R\, .
\end{equation*}
Consequently, $\chi_{\substack{{}\\{(0,e^\lambda ]}}} (b) =
p_{\lambda +0}\leq\chi_{\substack{{}\\{(0,e^\lambda ]}}} (d)$
for all $\lambda\in\mathbb R\,$.

The uniqueness of $b$ is an immediate consequence of (ii).
\medskip

{\it Step 2}. {\sc Proof of {\rm (iii)} and {\rm (iv)}.}
\medskip

(iii) is clear from the construction of $b$ in Step 1.

In order to verify (iv), let $\sigma$ be a $*$-automorphism of
$M$ such that, for some $\lambda_\sigma >0\,$,
\begin{equation*}
\sigma\circ\alpha_s =\alpha_{\lambda_\sigma s}\circ\sigma\, ,
\qquad s\in\mathbb R\, .
\end{equation*}
Then holds clearly
\begin{equation}\label{intertw.}
\sigma\circ\alpha_z =\alpha_{\lambda_\sigma z}\circ\sigma
\, ,\qquad z\in\mathbb C\, .
\end{equation}

There exists a faithful unital normal $*$-representation
$\pi : M\longrightarrow B(K)\,$, which is covariant with respect
to $\,\sigma\,$, that is $\pi\big( \sigma (x)\big) =
U \pi (x)\, U^* ,\, x\in M\,$, where $U$ is an appropriate
unitary on $K\,$: for example, we can choose
\smallskip

\noindent \hspace{0.25 in} $K=l^2(\mathbb Z ; H)\,$, the space
of all square-summable two-sided sequences in $H\,$,
\smallskip

\noindent \hspace{0.25 in} $\pi (x)\,
(\xi_k)_{\substack{{}\\{k\in\mathbb Z}}} =
\big( \sigma^k (x)\,\xi_k\big)_{\substack{{}\\{k\in\mathbb Z}}}\,$
for $\, x\in M\, ,\, (\xi_k)_{\substack{{}\\{k\in\mathbb Z}}}\in
l^2(\mathbb Z ; H)\,$,

\noindent \hspace{0.25 in} $U\,
(\xi_k)_{\substack{{}\\{k\in\mathbb Z}}} =
(\xi_{k+1})_{\substack{{}\\{k\in\mathbb Z}}}\,$ for
$\, (\xi_k)_{\substack{{}\\{k\in\mathbb Z}}}\in
l^2(\mathbb Z ; H)\,$.
\smallskip

\noindent Then $\big(\pi\circ\alpha_s\circ\pi^{-1}
\big)_{s\in\mathbb R}$ is an $so$-continuous one-parameter group
of $*$-automorphisms of the von Neumann algebra
$\pi (M)\subset B(K)\,$, $\pi (b)$ is an injective element of
$\pi (M)$ with $0\leq\pi (b)\leq 1_K$ and
\begin{equation*}
\big(\pi\circ\alpha_s\circ\pi^{-1}\big) \big(\pi (x)\big) =
\pi (b)^{-is}\pi (x)\,\pi (b)^{is}\, ,\qquad s\in\mathbb R\, ,\,
x\in M\, .
\end{equation*}
Moreover, by the definition of $b\,$, for any injective
$\pi (d)\in\pi (M)\, ,\, 0\leq\pi (d)\leq 1_K\,$, such that
\begin{equation*}
\big(\pi\circ\alpha_s\circ\pi^{-1}\big) \big(\pi (x)\big) =
\pi (d)^{-is}\pi (x)\,\pi (d)^{is}\, ,\qquad s\in\mathbb R\, ,\,
x\in M\, ,
\end{equation*}
we have
\begin{equation*}
\chi_{\substack{{}\\{(0,e^\lambda ]}}} \big(\pi (b)\big)
\leq\chi_{\substack{{}\\{(0,e^\lambda ]}}} \big(\pi (d)\big)\,
,\qquad \lambda\in\mathbb R\, .
\end{equation*}
Applying the above proved (iii) to $\pi (M)\,$, $\big(\pi\circ
\alpha_s\circ\pi^{-1}\big)_{s\in\mathbb R}\,$, $\pi (b)$
instead of $M\,$, $\alpha\,$, $b\,$, we
obtain that, for every $\lambda\in\mathbb R\,$,
$\pi\big( \chi_{\substack{{}\\{(0,e^\lambda ]}}} (b)\big) =
\chi_{\substack{{}\\{(0,e^\lambda ]}}} \big(\pi (b)\big)$ is
the orthogonal projection onto
\begin{equation*} \hspace{0.3 cm}
\bigcap_{\mu >\lambda}\;\text{the closed linear span of }
\pi (M)^{\pi\circ\alpha\circ\pi^{-1}}\big( [-\mu ,+\infty )\big)
\, K =
\end{equation*}
\begin{equation*}
= \bigcap_{\mu >\lambda}\;\text{the closed linear span of }
\pi \Big( M^\alpha \big( [-\mu ,+\infty )\big)\Big)\, K
\hspace{1.2 cm}
\end{equation*}

For every $\lambda\in\mathbb R\,$ and $x\in M\,$,
by (\ref{spectr.subsp.}) and by (\ref{intertw.}), the following
four conditions are equivalent:
\begin{equation*}
x\in M^\alpha \big( [-\lambda ,+\infty )\big)\, ,
\end{equation*}
\begin{equation*}
x\in \!\bigcap_{\substack{{z\in\mathbb C}\\{\Im z\geq 0}}}\!
{\rm Dom}\, (\alpha_z )\,\text{ and }\,
\big\|\sigma\big( \alpha_z (x) \big)\big\| = \|\alpha_z (x)\|
\leq e^{\lambda \Im z}\| x\|\,\text{ for all }\, z\in\mathbb C
\, ,\,\Im z\geq 0\, ,
\end{equation*}
\begin{equation*}
\sigma (x)\in \!\bigcap_{\substack{{z\in\mathbb C}\\{\Im z\geq 0}}}
\! {\rm Dom}\, (\alpha_z )\,\text{ and }\,
\big\|\alpha_{\lambda_\sigma z}\big( \sigma (x) \big)\big\|\leq
e^{\lambda \Im z}\| x\|\,\text{ for all }\, z\in\mathbb C\, ,\,
\Im z\geq 0\, ,
\end{equation*}
\begin{equation*}
\sigma (x)\in M^\alpha \big( [-\lambda_\sigma^{\, -1}\lambda ,
+\infty )\big)\, .
\end{equation*}
Therefore
\begin{equation}\label{subsp.}
\sigma \big( M^\alpha \big( [-\lambda ,+\infty )\big)\big) =
M^\alpha \big( [-\lambda_\sigma^{\, -1}\lambda ,+\infty )\big)\, ,
\qquad \lambda\in\mathbb R\, .
\end{equation}

Let next $\lambda\in\mathbb R$ be arbitrary. By the covariance
property of $\pi$ and by (\ref{subsp.}), we have for every
$\mu >\lambda\,$:
\begin{equation*}
U\,\pi \big( M^\alpha \big( [-\mu ,+\infty )\big)\big)\, K =\;
U\,\pi \big( M^\alpha \big( [-\mu ,+\infty )\big)\big)\, U^*K =
\end{equation*}
\begin{equation*} \hspace{3.45 cm}
= \pi\Big( \sigma \big( M^\alpha \big( [-\mu ,+\infty )\big)
\big)\Big) K =
\end{equation*}
\begin{equation*} \hspace{3.32 cm}
= \pi\Big( M^\alpha \big( [-\lambda_\sigma^{\, -1}\mu ,+\infty )
\big)\Big) K\, .
\end{equation*}
Consequently
$\, U\,\pi\big(\chi_{\substack{{}\\{(0,e^{\lambda} ]}}}
(b)\big) K = \pi\big(
\chi_{\substack{{}\\{(0,e^{\lambda_\sigma^{\, -1}\lambda} ]}}}
\big(\pi (b)\big)\big) K\,$, and so
\begin{equation*}
\pi\big( \sigma\big(\chi_{\substack{{}\\{(0,e^{\lambda} ]}}}
(b)\big)\big) =
U\,\pi\big(\chi_{\substack{{}\\{(0,e^{\lambda} ]}}}
(b)\big)\, U^* =\pi\big(
\chi_{\substack{{}\\{(0,e^{\lambda_\sigma^{\, -1}\lambda} ]}}}
(b)\big)\, ,
\end{equation*}
\begin{equation}\label{proj.}
\chi_{\substack{{}\\{(0,e^{\lambda} ]}}} \big(\sigma (b)\big) =
\sigma\big(\chi_{\substack{{}\\{(0,e^{\lambda} ]}}} (b)\big) =
\chi_{\substack{{}\\{(0,e^{\lambda_\sigma^{\, -1}\lambda} ]}}}(b)
=\chi_{\substack{{}\\{(0,e^{\lambda} ]}}} (b^{\lambda_\sigma})\, .
\end{equation}
\smallskip

Now, by (\ref{proj.}) we conclude that
$\sigma (b) =b^{\lambda_\sigma}\,$.

\hfill $\square\quad$
\medskip

{\it Proof of }$\,${\bf Theorem \ref{type}}.
\medskip

{\it Proof of }(1).
Let us consider $\gamma_s =\text{Ad}\, U(s)$ for all
$s\in\mathbb R$ (not only for $s\geq 0$) and let $z$ be an
arbitrary selfadjoint element of the center $Z(M)$ of $M\,$.
By (6) in Theorem \ref{hsmi} we have $\gamma_s (z)\in M$ for all
$s\geq 0\,$, hence $z$ and $\gamma_s (z)$ commute for any
$s\geq 0\,$. But then the elements of the set $\{ \gamma_s (z)
\, ;\, s\in\mathbb R \}$ are mutually commuting, so the von
Neumann algebra $C$ generated by this set is commutative.
Since
\smallskip

\noindent\hspace{0.6 cm} $\gamma_s =\text{Ad}\; U(s)\, |\, M$
leaves $\, C$ invariant for every $s\in\mathbb R\,$ and
\smallskip

\noindent\hspace{0.6 cm} $U(s)=\exp (isP)\,$, $s\in\mathbb R\,$,
for some positive selfadjoint operator $P$ in $H\,$,
\smallskip

\noindent according to the implementation theorem of Borchers
\cite{Bo1} (see also Theorem \ref{pos.gen}) there exists an
element $b\in C\, ,\, 0\leq b\leq 1_H\,$, such that
\begin{equation*}
\gamma_s (x) =b^{-is} x\, b^{is} =x\, ,\qquad x\in C\, ,\,
s\in\mathbb R\, .
\end{equation*}
Consequently, $\gamma_s (z) =z$ for all $s\in\mathbb R\,$.
\medskip

{\it Proof of }(2).
By (5) in Theorem \ref{hsmi} and by the above proved (1), $Z(N)=
Z\big(\gamma_1 (M)\big) =Z(M)\,$. Now it is easy to see that
the projection family
\begin{center}
$\{\, q\in Z(M)\, ;\, M q =N q\,\}$
\end{center}
is upward directed and its lowest upper bound is the greatest
projection $p\in Z(M)$ satisfying $M p =N p\,$.

The implication $\, e\leq p\;\Longrightarrow\, U(s)\, e=e\,$
holds for any projection $e\in M\,$, because
\begin{center}
$M p =N p\;\Longrightarrow\; \varphi_p =\psi_p
\;\overset{(\ref{centr.red.mod})}{\Longrightarrow}\;
\Delta_M^{\, it}\, p = \Delta_N^{\, it}\, p\, ,\; t\in\mathbb R
\;\Longrightarrow\; U(s)\, p =p\, ,\; s\in\mathbb R\, .$
\end{center}
\medskip

Now let $e\in M$ be an arbitrary projection such that
\begin{equation}\label{fix}
U(s)\, e=e\, ,\qquad s\in\mathbb R\, .
\end{equation}
For every $a\in\mathfrak A_\varphi$ and $x\, , y\in
\mathfrak N_\varphi\,$, using (8) in Theorem \ref{hsmi}, we
deduce that
\begin{center}
$( U(s)\, a\, e\, J_M\, x_\varphi\, |\, J_M\, y_\varphi )
\overset{(\ref{fix})}{=}
\big(\gamma_s (a)\, e\, J_M\, x_\varphi\,\big|\, J_M\, y_\varphi
\big) =\big( e\, J_M\, x_\varphi\,\big|\, \gamma_s (a^*)\, J_M\,
y_\varphi\big)$
\end{center}

\noindent\hspace{3.88 cm} $\overset{(\ref{j})}{=} \big( e\, J_M\,
x_\varphi\,\big|\, J_M\, y\, J_M\,\gamma_s (a^*)_\varphi\big)$
\smallskip

\noindent\hspace{4.035 cm} $=\big( J_M\, y^*\, J_M\, e\, J_M\,
x_\varphi\,\big|\, U(s)\, (a^*)_\varphi\big)$
\smallskip

\noindent\hspace{4.035 cm} $=\big( e\, J_M\, y^* x_\varphi\,\big|
\, U(s)\, (a^*)_\varphi\big)$
\smallskip

\noindent\hspace{3.88 cm} $\overset{(\ref{fix})}{=}\big( e\, J_M\,
y^* x_\varphi\,\big|\, (a^*)_\varphi\big)$
\smallskip

\noindent does not depend on $s\geq 0\,$, so
$\big( \big( 1_H -U(s)\big)\, a\, e\, J_M\, x_\varphi\,\big|
\, J_M\, y_\varphi\big) =0\,\text{ for all }\, s\geq 0\,$. By
$\overline{\{ J_M\, x_\varphi\, ;\, x\in\mathfrak N_\varphi\} }
=H$ and $\,\overline{\mathfrak A_\varphi }^{\, so}=M\,$, we get
\begin{equation*}
\big( 1_H-U(s)\big)\, M\, e\, H=\{ 0\}\, ,\qquad s\in\mathbb R\, .
\end{equation*}
Since the orthogonal projection onto the closed linear span of
$M\, e\, H$ is the central support $z(e)\in Z(M)$ of $e\,$, we
obtain that $\big( 1_H-U(s)\big)\, z(e)=0$ for all $s\geq 0\,$,
hence
\begin{equation*}
U(s)\, z(e) =z(e)\, ,\qquad s\in\mathbb R\, .
\end{equation*}
Consequently, $N z(e) =U(1)\, M\, U(1)^*\, z(e) =M z(e)$ and so
$e\leq z(e)\leq p\,$.

Finally, let $e\in M^\varphi$ be a projection such that
\begin{equation}\label{Ue} \hspace{1.39 cm}
U(s)\, e =e\, U(s)\, ,\qquad\; s\in\mathbb R\, ,
\end{equation}
\begin{equation}\label{UeJeJ}
U(s)\, e J_M e J_M =e J_M e J_M\, ,\qquad\; s\in\mathbb R\, .
\phantom{x}
\end{equation}

If $\pi =\pi_{\substack{{}\\ e J_\varphi e J_\varphi}} :
e M e\longrightarrow B(e J_\varphi e J_\varphi H)$ is the
$*$-representation defined in (\ref{stand.red}), then, for every
$s\geq 0\,$, (\ref{Ue}) and (\ref{UeJeJ}) yield that
$\gamma_s (eMe)\subset eMe$ and
\begin{equation*}
\pi\big( \gamma_s (a)\big) =U(s)\, a\, U(-s)\, |\, e J_M e J_M H
=a\, |\, e J_M e J_M H =\pi (a)\, ,\qquad a\in eMe\, .
\end{equation*}
Since $\pi$ is faithful, we obtain that $U(s)\, a\, U(-s) =
\gamma_s (a) =a$ for all $s\geq 0$ and all $a\in eMe\,$. In other
words, every $U(s)$ commutes with every operator in $eMe\,$.  

Consequently, for every $s\in\mathbb R\,$, the unitary
$\, U(s)\, |\, eH : eH\rightarrow eH$ belongs to the commutant
of the reduced von Neumann algebra $\big\{\, x\, |\, eH : eH
\rightarrow eH\, ;\, x\in eMe\,\big\}\,$, hence to the induced
von Neumann algebra $\big\{\, x'\, |\, eH : eH\rightarrow
eH\, ;\, x'\in M'\,\big\}\,$. Since the kernel of the induction
$*$-homomorphism $M'\ni x'\mapsto x'\, |\, eH$ is
$M' \big( 1_H -z(e)\big)\,$, where $z(e)\in Z(M)$ stands for the
central support of $e\,$, there exists a
one-parameter group $(u'_{\! s} )_{s\in\mathbb R}$ of unitaries
in $M'z(e)$ such that
\begin{equation*}
U(s)\, |\, eH =u'_{\! s}\, |\, eH\, ,\text{ that is }\,
U(s)\, e =u'_{\! s} e\, ,\qquad s\in\mathbb R\, .
\end{equation*}
Setting $u_s =J_M u'_{-s} J_M\,$, $(u_s )_{s\in\mathbb R}$ is a
one-parameter group of unitaries in $M z(e)$ such that
\begin{equation*}
U(s) J_M e J_M =J_M U(-s)\, e J_M =J_M u'_{-s} e J_M =u_s
J_M e J_M\, ,\qquad s\in\mathbb R\, .
\end{equation*}
Therefore we have, for every $s\geq 0$ and $a\in M z(e)\,$,
\begin{align*}
\gamma_s (a) J_M e J_M &=J_M e J_M U(s)\, a\, U(-s) J_M e J_M =
J_M e J_M u_s a u_{-s} J_M e J_M \\
&= (u_s a u_{-s}) J_M e J_M\, .
\end{align*}
Since the kernel of the induction $*$-homomorphism $M\ni x\mapsto
x\, |\, J_M e J_M H$ is equal to $M\big( 1_H -z(J_M e J_M)\big)
=M\big( 1_H -z(e)\big)\,$, we obtain that
\begin{equation*}
\gamma_s (a) =u_s a u_{-s}\, ,\qquad s\geq 0\, ,\; a\in Mz(e)\, .
\end{equation*}
In particular, $Nz(e)=\gamma_1\big( M z(e)\big) =u_1 Mz(e) u_{-1}
=Mz(e)\,$, and so $e\leq z(e)\leq p\,$.
\medskip

{\it Proof of }(3).
Let us assume that $\displaystyle M_{-\infty} =
\bigcap\limits_{s\geq 0} \gamma_s (M)$ contains $M^\varphi\,$.

$\text{Ad}\; U(s)$ leaves $M_{-\infty}$ invariant for every
$s\in\mathbb R\,$, defining thus an $so$-continuous one-parameter
group $(\alpha_s)_{s\in\mathbb R}$ of $*$-automorphisms of
$M_{-\infty}\,$. Using Theorem \ref{pos.gen}, we get an
injective $b\in M_{-\infty}\,$, $0\leq b\leq 1_H\,$, such that
\begin{equation}\label{bound.implem}
\alpha_s (x) =b^{-is} x\, b^{is}\, ,\qquad s\in\mathbb R\, ,\,
x\in M_{-\infty}
\end{equation}
and
\begin{equation}\label{transfer}
\left. \begin{array}{l} \sigma\text{ $*$-automorphism of }
M_{-\infty}\text{ and }\lambda_\sigma >0 \\
\sigma\circ\alpha_s =\alpha_{\substack{{}\\ {\lambda_\sigma s}}}
\circ\sigma\text{ for all }s\in\mathbb R \end{array} \right\}
\;\,\Longrightarrow\;\,\sigma (b) =b^{\lambda_\sigma}\, .
\end{equation}
\smallskip
\noindent Since $b\in M_{-\infty}\,$, (\ref{bound.implem}) yields
that
\begin{equation}\label{bfixed}
U(s)\, b\, U(s)^*=\alpha_s(b) =b\, ,\text{that is }b
\text{ commutes with }U(s)\, ,\quad s\in\mathbb R\, .
\end{equation}
In particular, $b\in U(1)\, M\, U(1)^* =N\,$. Furthermore, by
(1) in Theorem \ref{hsmi} we have
\begin{equation}\label{sigma-gamma}
\sigma^\varphi_t\circ\gamma_s =
\gamma_{\substack{{}\\ {e^{-2\pi t} s}}}\circ
\sigma^\varphi_t\, ,\qquad s\, ,\, t\in\mathbb R\, ,
\end{equation}
so $M_{-\infty}$ is left invariant by all $\sigma^\varphi_t\,$.
Therefore, applying (\ref{transfer}) with
$\sigma =\sigma^\varphi_t\,$, we get
\begin{equation}\label{mappingb}
\sigma^\varphi_t (b) =b^{\, e^{-2\pi t}} ,\qquad t\in
\mathbb R\, .
\end{equation}

Let $\chi_{\substack{{}\\{\{ \lambda\}}}}$ be the characteristic
function of $\{ \lambda\}\subset\mathbb R\,$. For every $t\in
\mathbb R$ we have
\begin{equation*}
\sigma^\varphi_t (b) =b^{\, e^{-2\pi t}}\,\Longrightarrow\;
\sigma^\varphi_t \big(\chi_{\substack{{}\\{\{ 1\}}}} (b)\big)
=\chi_{\substack{{}\\{\{ 1\}}}} (b)\, ,
\end{equation*}
so $\chi_{\substack{{}\\{\{ 1\}}}} (b)\in M^\varphi\,$.
On the other hand, (\ref{bfixed}) implies that
$\chi_{\substack{{}\\{\{ 1\}}}} (b)$ commutes with all $U(s)\,$.
Therefore $\chi_{\substack{{}\\{\{ 1\}}}} (b)\in N$ commutes
with all $\Delta_N^{\, it} = U(1)\,\Delta_M^{\, it} U(1)^*\,$,
and so it belongs to $N^\psi\,$.

Let us denote
\begin{itemize}
\item[] $\hspace{1.2 cm}e_o =1_H -\chi_{\substack{{}\\{\{ 1\}}}} (b)
\in M^\varphi \cap N^\psi\,$,
\item[] $\hspace{1.2 cm}\varphi_o =$ the restriction of $\varphi$ to
$e_o\, M\, e_o\,$,
\item[] $\hspace{1.2 cm}\psi_o =$ the restriction of $\psi$ (hence
also of $\varphi$) to $e_o\, N\, e_o\,$,
\item[] $\hspace{1.2 cm}b_o =e_o\, b\, e_o =b\, e_o =e_o\, b
\in e_o\, M\, e_o\,$.
\end{itemize}
\smallskip
\noindent By \cite{P-T}, Proposition 4.1 and Theorem 4.6
(see also \cite{S}, Propositions 4.5 and 4.7), $\varphi_o$ is a
normal semifinite faithful weight and its modular group is the
restriction of the modular group of $\varphi$ to $e_o\, M\, e_o\,$.
Similarly, $\psi_o$ is a normal semifinite faithful weight and
its modular group is the restriction of the modular group of
$\psi$ to $e_o\, N\, e_o\,$.
In particular, by (\ref{mappingb}), we have
\begin{equation}\label{mappingb_o}
\sigma^{\varphi_o}_t (b_o)= {b_o}^{e^{-2\pi t}}\, ,\qquad
t\in\mathbb R\, .
\end{equation}

Since $0\leq b_o\leq e_o$ and the supports of both $b_o$ and
$e_o -b_o$ are equal to the unit $e_o$ of $e_o\, M\, e_o\,$,
$-\log b_o$ is a positive selfadjoint linear operator, of support
$e_o$ and affiliated with $e_o\, M\, e_o\,$. Consequently, defining
\begin{equation*}
u_s =(-\log b_o )^{i s}\, ,\qquad s\in\mathbb R\, ,
\end{equation*}
$(u_s)_{\substack{{}\\{s\in\mathbb R}}}$ is a strongly continuous
one-parameter group of unitaries in $e_o\, M\, e_o$ and
(\ref{mappingb_o}) yields
\begin{equation*}
\sigma^{\varphi_o}_t (u_s) =e^{-2\pi tsi}\, u_s\, ,\qquad s
\, ,\, t\in \mathbb R\, .
\end{equation*}
Now the characterization theorem of M. Landstad \cite{La}, Theorem 2
(see also \cite{S-V-Z}, Theorems I.3.3 and I.3.4, or \cite{S},
Theorem 19.9) implies that the von Neumann algebra $e_o\, M\, e_o$
is generated by $(e_o\, M\, e_o)^{\varphi_o} =e_o\, M^\varphi e_o$
and by $u_{\substack{{}\\ {\mathbb R}}}\,$, that is
\begin{equation}\label{Land}
e_o\, M^\varphi e_o\,\text{ and }\, b_o
\text{ generate the von Neumann algebra }e_o\, M\, e_o\, .
\end{equation}

Since $M^\varphi\subset M_{-\infty}$ and $b\in M_{-\infty}\,$,
we get that $e_o\, M\, e_o\subset M_{-\infty}\subset
\gamma_1 (M) =N\,$, that is $e_o\, M\, e_o =e_o\, N\, e_o\,$.
Consequently $\varphi_o =\psi_o\,$, and so the modular groups
$\sigma^\varphi$ and $\sigma^\psi$ have the same restriction
$\sigma^{\varphi_o} =\sigma^{\psi_o}$ on $e_o\, M\, e_o =
e_o\, N\, e_o\,$. Using (3) in Theorem \ref{hsmi}, we obtain
for every $x\in e_o\, M^\varphi\, e_o\subset M^\varphi \cap
N^\psi$ and $t\in \mathbb R\,$:
\begin{equation*}
U(1-e^{2\pi t})\, x\, U(1-e^{2\pi t})^* = \Delta_N^{-it}
\Delta_M^{it}\, x\, \Delta_M^{-it}\Delta_N^{it} =
\Delta_N^{-it}\,\sigma^\varphi_t (x)\,\Delta_N^{it} =
\sigma^\psi_{-t} (x) =x\, .
\end{equation*}
Therefore $e_o\, M^\varphi e_o\subset \big\{ x
\in M\, ;\, U(s)\, x =x\, U(s)\, ,\, s\in\mathbb R \big\}\,$,
which yields together with (\ref{bfixed}) and (\ref{Land}) :
\begin{equation*}
e_o\, M\, e_o\subset \big\{ x\in M\, ;\, U(s)\, x =
x\, U(s)\, ,\, s\in\mathbb R \big\}\, .
\end{equation*}
In other words, every $\alpha_s$ acts identically on
$e_o\, M\, e_o\subset M_{-\infty}\,$. By (\ref{bound.implem})
we conclude that $b_o$ belongs to the center of $e_o\, M\, e_o\,$.

Since $b_o\in Z(e_o\, M\, e_o)$ is invariant under the modular
automorphism group of $\varphi_o\,$, which coincides with the
restriction of the modular automorphism group of $\varphi$ to
$e_o\, M\, e_o$ as discussed above, we have $b_o\in M^\varphi\,$.
Taking into account (\ref{mappingb_o}), we obtain
\begin{equation*}
b^{\, e^{-2\pi t}} e_o ={b_o}^{e^{-2\pi t}} =
\sigma^\varphi_t (b_o) =b_o =b\, e_o\, ,\qquad t\in\mathbb R\, ,
\end{equation*}
which is possible only if $e_o=0\,$. Consequently
$\chi_{\substack{{}\\{\{ 1\}}}} (b)=1_H\,$, that is $b=1_H\,$.
But then every $\alpha_s =\text{Ad}\; U(s)\, |\, M_{-\infty}$
acts identically on $M_{-\infty}\,$, hence
\begin{equation*}
M_{-\infty}\subset \big\{ x\in M\, ;\, U(s)\, x =
x\, U(s)\, ,\, s\in\mathbb R \big\} =\big\{ x\in M\, ;\,
\gamma_s (x)=x\, ,\, s\geq 0\big\}\, .
\end{equation*}

{\it Proof of }(4). Let us first assume that
$\overline{\mathfrak M_\varphi\cap M^\varphi}^{\, so}= M^\varphi\,$.
Taking into account (8) in Theorem \ref{hsmi}, the inclusion
$M^\varphi\subset\big\{ x\in M\, ;\,\gamma_s (x) =x\, ,\, s\geq 0
\big\}$ will follow once we show that
\begin{equation*}
s>0\; ,\, x\in \mathfrak M_\varphi\cap M^\varphi\;\Longrightarrow\;
U(s)\, x_\varphi =x_\varphi\, .
\end{equation*}

For let $s>0\,$ and $x\in \mathfrak M_\varphi\cap M^\varphi$ be
arbitrary. Using (8) in Theorem \ref{hsmi}, (\ref{delta}) and
(\ref{sigma-gamma}), we get successively
\smallskip
\begin{center}
$\Delta_M^{\, it}\, U(s)\, x_\varphi =\sigma^\varphi_t\big(\gamma_s
(x)\big)_\varphi =\gamma_{\substack{{}\\ {e^{-2\pi t} s}}}
\big( x\big)_\varphi = U(e^{-2\pi t} s)\, x_\varphi
\buildrel {t\to +\infty\,}\over{\hbox to1.2cm{\rightarrowfill}}
x_\varphi\, ,$
\end{center}
\smallskip
\begin{center}
$\|\, U(s)\, x_\varphi - x_\varphi\| =\|\,\Delta_M^{\, it} \big(
U(s)\, x_\varphi - x_\varphi \big)\| =\|\,\Delta_M^{\, it}\,
U(s)\, x_\varphi -\sigma^\varphi_t(x)_\varphi\|
\buildrel {t\to +\infty\,}\over{\hbox to1.2cm{\rightarrowfill}}
0\, .$
\end{center}
\smallskip

For the second implication we prove that if
\begin{equation}\label{centr.commutes}
M^\varphi\subset\big\{ x\in M\, ;\,\gamma_s (x) =x\, ,\,
s\geq 0\big\}
\end{equation}
and $p\neq 1_H\,$, then $M (1_H-p)$ is of type III$_1\,$. Then also
$N(1_H-p) =\gamma_1\big( M(1_H-p))\big)$ will be of type III$_1\,$.

Taking into account (\ref{typeIII_1}), we have to prove that
\begin{equation*}
e\in M^\varphi\text{ projection},\, 0\neq e\leq 1_H-p\;
\Longrightarrow\;\sigma \big(\,\Delta_M\, |\, e J_M e J_M H\,\big) =
[\, 0,+\infty )\, .
\end{equation*}
For this purpose, let the projection $e\in M^\varphi\, ,\, 0\neq
e\leq 1_H-p\,$, be arbitrary. By the assumption (\ref{centr.commutes})
we have $\gamma_s (e) =e$ for all $s\geq 0\,$, hence (\ref{Ue}) holds.
Since $0\neq e\leq 1_H-p\,$, the above proved (2) entails that
(\ref{UeJeJ}) does not hold, that is $U(s)\, eJ_M eJ_M\neq
eJ_M eJ_M$ for some $s\in\mathbb R\,$. Nevertheless, by (\ref{Ue})
and by (2) in Theorem \ref{hsmi}, all $U(s)$ commute with
$eJ_M eJ_M\,$. Since $e\in M^\varphi\,$, also every $\Delta_M^{it}$
commutes with $eJ_M eJ_M\,$.

According to (\ref{generation}), $U(i)=\exp\, (-P)\in B(H)$ is
injective and $0\leq U(i)\leq 1_H\,$. By the commutation relation
(1) in Theorem \ref{hsmi}, we have
\begin{equation}\label{mappingU(i)}
\Delta_M^{it}\, U(i)\, \Delta_M^{-it}=U(i)^{\, e^{-2\pi t}} ,
\qquad t\in \mathbb R\, .
\end{equation}
Consequently, the spectral projection $f_o =
\chi_{\substack{{}\\{\{ 1\}}}} \big( U(i)\big)$ commutes with
every $\Delta_M^{it}\,$. On the other hand, $f_o$ clearly
commutes with every $U(s)\,$. Finally, the commutation of
$eJ_M eJ_M$ with all $U(s)$ implies that the projections $eJ_M eJ_M$
and $f_o$ commute.

We have already seen
that $U(s)\, eJ_M eJ_M\neq eJ_M eJ_M$ for some $s\in\mathbb R\,$.
On the other hand, since $U(s) =U(i)^{-is}\,$, we have
$U(s)\, f_o =f_o$ for every $s\in\mathbb R\,$. Therefore
$eJ_M eJ_M\not\leq f_o\,$, and so the projection
\begin{equation*}
f_1 =eJ_M eJ_M -f_o eJ_M eJ_M\leq eJ_M eJ_M
\end{equation*}
is not zero. Since all $\Delta_M^{it}$ and all $U(s)$ commute
with $eJ_M eJ_M$ and with $f_o\,$, they commute also with $f_1\,$.
Therefore, we can define an $so$-continuous one-parameter group
$(v_t)_{t\in\mathbb R}$ of unitaries on $f_1 H$ by setting
\begin{equation}\label{restr.Delta}
v_t =\Delta_M^{it}\, |\, f_1 H =(\Delta_M\, |\, f_1 H)^{it}\, .
\qquad t\in\mathbb R\, ,
\end{equation}
as well as the operator
\begin{equation*}
b_1 =U(i)\, |\, f_1 H : f_1 H\longrightarrow f_1 H\in B(f_1 H)\, ,
\end{equation*}
for which $0\leq b_1\leq f_1$ and $\, b_1\, , f_1 -b_1$ are
injective. From (\ref{mappingU(i)}) we get successively
\begin{equation*}
v_t^{}\, b_1^{} v_t^*=b_1^{\, e^{-2\pi t}} ,\qquad t\in \mathbb R\, ,
\end{equation*}
\begin{equation*}
v_t^{}\, (-\log b_1 )^{is} v_t^*=e^{-2\pi tsi}\, (-\log b_1 )^{is}
\, ,\qquad t\, ,\, s\in \mathbb R\, .
\end{equation*}

Now the Stone-von Neumann Uniqueness Theorem for canonical
commutation relations (see e.g. \cite{Ba-Ra}, Ch. 20, \S$\,$ 2 or
\cite{Sl}) entails that there exist a Hilbert space $K\neq\{ 0\}$
and a unitary operator
\smallskip
\begin{center}
$W : K\,\overline{\otimes}\, L^2(\mathbb R )\longrightarrow f_1 H$
\end{center}
\smallskip
such that
\begin{equation*}
v_t =W\circ (1_K\overline{\otimes}\, m_{-2\pi t})\circ W^*\, ,\;
(-\log b_1 )^{is}=W\circ (1_K\overline{\otimes}\,\lambda_s )\circ W^*
\, ,\qquad t\, ,\, s\in\mathbb R\, ,
\end{equation*}
where
\begin{center}
\hspace{2.5 mm}$m_t$ is the multiplication operator with
$e^{it\,\cdot}$ on $L^2(\mathbb R )\, ,\qquad t\in\mathbb R\,$,
\end{center}
\smallskip
\begin{center}
$\lambda_s$ is the translation operator $\xi\longmapsto\xi
(\,\cdot\, -s)$ on $L^2(\mathbb R )\, ,\qquad s\in\mathbb R\,$.
\end{center}
\smallskip
Using (\ref{restr.Delta}), we deduce that $\,\Delta_M\, |\, f_1 H =
W\circ (1_K\overline{\otimes}\, m_{2\pi i})\circ W^*\,$, where
$m_{2\pi i}$ is the unbounded positive selfadjoint multiplication
operator with $e^{-2\pi\,\cdot}$ in $L^2(\mathbb R )\,$.
Consequently, the spectrum of $\Delta_M\, |\, f_1 H$
is equal to the spectrum of $1_K\overline{\otimes}\, m_{2\pi i}\,$,
that is to $[0,+\infty )\,$. Since $f_1\leq eJ_M eJ_M\,$, we
conclude that also the spectrum of $\Delta_M\, |\, eJ_M eJ_M H$
is equal to $[0,+\infty )\,$.

\hfill $\square\quad$

\bigskip

\end{document}